\documentclass[12pt,a4paper]
{article}
\usepackage{amsfonts}
\usepackage{amsmath}
\usepackage{mathrsfs}
\usepackage{graphicx}
\usepackage{amsmath,amsthm,amssymb,amscd}
\newtheorem{Thm}{Theorem}[section]
\newtheorem{Lem}[Thm]{Lemma}
\newtheorem{Prop}[Thm]{Proposition}
\newtheorem{Def}[Thm] {Definition}

\newtheorem{Cor}[Thm]{Corollary}
\theoremstyle{remark}
\newtheorem{Rem} [Thm]{Remark}

\theoremstyle{claim}

\newtheorem{Que}[Thm]{Question}

\DeclareMathOperator{\Diff}{Diff}

\begin{document}

\begin{center} {\large \bf
Different Asymptotic Behavior versus Same Dynamical Complexity: Recurrence \& (Ir)Regularity}
\end{center}
\smallskip

\begin{center}
   Xueting Tian$^{*}$\\
   \medskip

  School of Mathematical Sciences, Fudan University, Shanghai 200433,\\ People's Republic of China \\
 \medskip

  E-mail: {\it xuetingtian@fudan.edu.cn; tianxt@amss.ac.cn } \\
\end{center}

\footnotetext{  Version of  January  2015.}


\footnotetext
{$^{*}$ Tian is  supported by National Natural Science Foundation of China (No. 11301088).
 }


\footnotetext{ Key words and phrases:   topological entropy;  subshifts of finite type;  $\beta-$shifts;  uniformly hyperbolic systems;  periodic, almost periodic point,  weakly almost periodic point and quasi-weakly  almost periodic point; regular, quasiregular  and irregular point;  specification property; fractal geometry.}

\medskip

\footnotetext{AMS Review:    37B10; 37B20; 37B40; 37D20; 37C45; 54H20 }

\def\abstractname{\textbf{Abstract}}

\begin{abstract}\addcontentsline{toc}{section}{\bf{English Abstract}}


For any dynamical system $T:X\rightarrow X$ of a compact metric space $X$ with $g-$almost product property and uniform  separation property, under the assumptions that   the periodic points are dense in $X$   and the periodic measures are dense in the space of invariant measures,  we distinguish various   periodic-like recurrences  and find that they all carry full topological topological entropy and so do their      gap-sets. In particular, this implies that any two kind of  periodic-like recurrences are essentially different.
  Moreover, we coordinate  periodic-like recurrences  with (ir)regularity  and obtain lots of generalized  multi-fractal analysis for all continuous observable functions.
These results are suitable for all $\beta-$shfits ($\beta>1$), topological mixing subshifts of finite type, topological mixing  expanding maps or topological mixing hyperbolic diffeomorphisms, etc.

\bigskip

 Roughly speaking, we combine  many  different ``eyes" (i.e., observable functions and periodic-like recurrences) to observe  the dynamical  complexity  and
obtain a  {\it Refined Dynamical Structure} for Recurrence Theory and Multi-fractal Analysis.

\end{abstract}

\newpage

\section{Introduction} \setlength{\parindent}{2em}

In the theory of  dynamical systems, i.e., the study of the asymptotic
behavior of orbits $\{T^n (x)\}_{n\in \mathbb{N}}$ (denoted by $Orb(x)$) when $T:X\rightarrow X$  is a continuous map of a compact
metric space $X$ and $x\in X$, one may say that two fundamental problems are to
understand how to partition different asymptotic behavior and how the points with same asymptotic behavior control or determine the complexity of system $T$.

Topological entropy is   a classical concept to describe the dynamical complexity. In this paper we are mainly to deal with a certain class of dynamical systems and show that various subsets characterized by
distinct asymptotic behavior all carry full topological entropy. To make this more precise let us introduce the following terminology. $T:X\rightarrow X$  is a continuous map of a compact
metric space $X$.

\begin{Def}\label{def-full-entropygaps} For a collection of subsets $Z_1,Z_2,\cdots,Z_k\subseteq X$ ($k\geq 2$), we say $\{Z_i\}$ has {\it full entropy gaps} with respect to $Y\subseteq X$ if
$$ h_{top} (T, (Z_{i+1}\setminus Z_i)\cap Y )=h_{top}(T, Y) \,\,\,\text{ for all } 1\leq i<k,$$
where $h_{top}(T,Z)$ denotes the topological entropy of a set $Z\subseteq X.$
\end{Def}
Often, but not always, the sets $Z_i$ are nested ($Z_i\subseteq Z_{i+1}$).
 Remark that for any system with zero topological entropy, it is obvious that any collection $\{Z_i\}$ has full entropy gaps with respect to any  $Y\subseteq X$.
Notice that  if   $X$ is a finite set, then any system on $X$ is simple and  carries zero entropy.   Thus in present paper, we always assume that $$\text{ \bf $X$ is a compact metric space with   infinite points. }$$

 In this paper, we consider following subsets of $X$ according to  different asymptotic behavior:
\begin{eqnarray*} Per(T)&:=&\{\text{ periodic points of } T\},\
  \\ A(T)&:=&\{\text{ almost periodic points of } T\}=\{\text{ points contained in minimal set}\},\
   \\ W(T)&:=&\{\text{ weakly almost  periodic points of } T\},\
    \\ QW(T)&:=&\{\text{ quasi-weakly almost periodic  points of } T\},\
     \\ Rec(T)&:=&\{\text{ recurrent points of } T \},\
 \\ \Omega (T)&:=&\{\text{ non-wandering points of } T\}.
 \end{eqnarray*}
Let $M_x(T)$ be the limit set of the empirical measures for $x$, $C_x$ be the  minimal set of attraction  for $x$, and $S_\mu$ be the support of $\mu$.  We also consider $$
    V (T):=\{x\in QW(T) |\,\exists \,\mu\in M_x(T) \text{ such that } S_\mu=C_x\}.$$
   Most of notions in above  considered sets  are well-known.
   The notions of periodic, recurrent and non-wandering can be found in \cite{Walter},
     the notion of almost periodic or minimal can be seen in
     \cite{Birkhoff,Gottschalk44-2222,Gottschalk44,Gottschalk46,Mai}  and others, for example,
      see \cite{Zhou-center-measure,Zhou93,Zhou95,ZF}.
   We will recall their  definitions and relations  in Section  \ref{section-20150108-recurrence}.
      Such sets are all $T$-invariant and they satisfy $ Per(T)\subseteq A(T)\subseteq W(T)
 \subseteq V(T) \subseteq QW(T)\subseteq Rec(T)\subseteq \Omega(T).$
  Figure \ref{Recurrence-2015} is a simple Venn diagram illustrating the containment between various sets  $\{Per(T), A(T), W(T), V(T), QW(T), $ $  Rec(T), \Omega(T)\}$ (simply, writing $\{Per, A, W, V,QW,$ $Rec,\Omega\}$ in the figure).

 \setlength{\unitlength}{1mm}
  \begin{figure}[!htbp]
  \begin{center}
  \begin{picture} (180,65) (0,0)
  \put (0,0){\scalebox{1}[1]{\includegraphics[0,0][50,30]{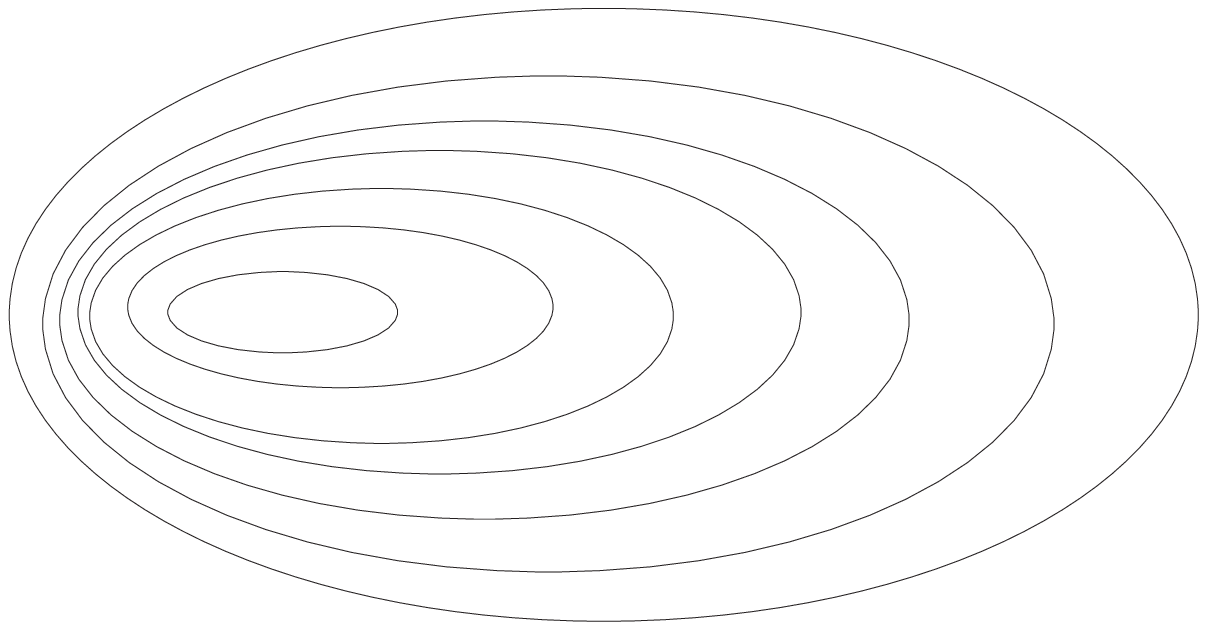}}}
   \put (43,34){$Per$}
    \put (114,34){$Rec$}  \put (131,34){$\Omega$}
\put (62,34 ){$A   $} \put (81,34 ){${\large W} $} \put (92,34 ){${\large V} $} \put (101,34 ){${\large QW} $}
  \end{picture}
  \caption{$\Omega(T)$.}
  \label{Recurrence-2015}
  \end{center}
  \end{figure}

A point  $x\in X$ is generic for some invariant measure  $\mu$ means that $M_x(T)=\{\mu\} $ (or equivalently,  Birkhoff averages of all
   continuous functions converge to the integral of $\mu$).    Let $G_\mu$   denote the set of all generic points for  $\mu $.
    Let    $M (T,X),$ $ M_{erg} (T,X)$ and $ M_{p} (T,X)$  denote the set of all $T$-invariant measures, $T$-ergodic measures and $T$-periodic measures  respectively.
We also consider \begin{eqnarray*}
 QR(T)&:=&\{\text{quasiregular points of } T\}= \cup_{\mu\in M(T,X)} G_\mu, \\
 I(T)&:=&\{\text{irregular points of } T\}=X\setminus QR(T),\\
 QR_{erg}(T) &:=&\{\text{points generic  for ergodic measures}\}
 = \cup_{\mu\in M_{erg}(T,X)} G_\mu, 
 \\ QR_{d}(T) &:=&\{\text{points of density in } QR(T)\}= \cup_{\mu\in M(T,X)} (G_\mu\cap S_\mu), 
 \\ R(T) &:=&\{\text{regular points of } T\}= QR_{d}(T)\cap QR_{erg}(T)= \cup_{\mu\in M_{erg}(T,X)} (G_\mu\cap S_\mu).
 \end{eqnarray*} Such sets are all $T-$invariant and remark that $$R(T)\subseteq QR_{d}(T)\cup QR_{erg}(T) \subseteq QR(T)=X\setminus I(T).$$  Most notions except irregular point are from \cite{Oxt} (for quasiregular point, also see \cite{DGS}) and the notion of irregular point can be found in \cite{Pesin-Pitskel1984,Bar,Takens,Barreira-Schmeling2000} etc.
  We will recall them more precisely in Section \ref{section-regular----irregular}.


\setlength{\unitlength}{1mm}
  \begin{figure}[!htbp]
  \begin{center}
  \begin{picture} (180,65) (0,0)
  \put (0,0){\scalebox{1}[1]{\includegraphics[0,0][50,30]{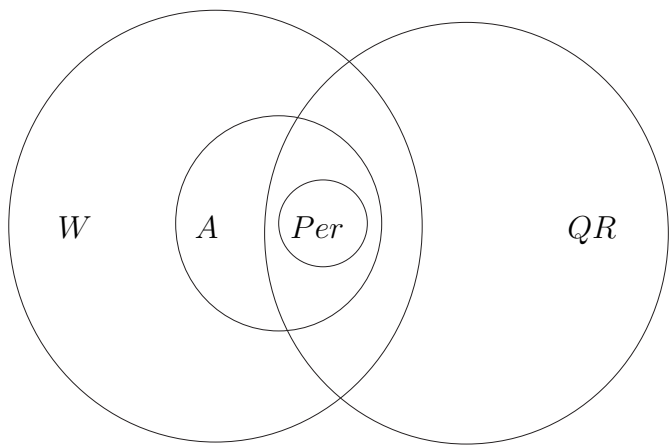}}}
   \put (45,34.5){$W$}
    \put (112,34.5){$QR$}
\put (63,34.5){$A   $} \put (75.5,34.5){${\large Per} $}
  \end{picture}
  \caption{$QR(T)\cup W(T)$.}
  \label{RecurrenceRegular-2015-1}
  \end{center}
  \end{figure}

Figure \ref{RecurrenceRegular-2015-1} is a simple Venn diagram illustrating the relations  between various sets  $\{Per(T),$ $ A(T), W(T),   QR(T)\}$ (simply, writing $\{Per, A, W, QR\}$ in the figure).   Figure \ref{RecurrenceIRRegular-2015} is a simple Venn diagram to illustrate the  relations between  following various  sets  $\{Per(T), A(T), W(T), V(T)\setminus  W(T), QW(T)\setminus  W(T), I(T)\}$ (simply, writing $\{Per, A, W, Q',$ $Q,I\}$ in the figure. Remark that $QW(T)\setminus  W(T)\subseteq I(T)$, see Theorem \ref{QW-W-in-Irregular}). Precise discussions will appear later.

\setlength{\unitlength}{1mm}
  \begin{figure}[!htbp]
  \begin{center}
  \begin{picture} (180,65) (0,0)
  \put (0,0){\scalebox{1}[1]{\includegraphics[0,0][50,30]{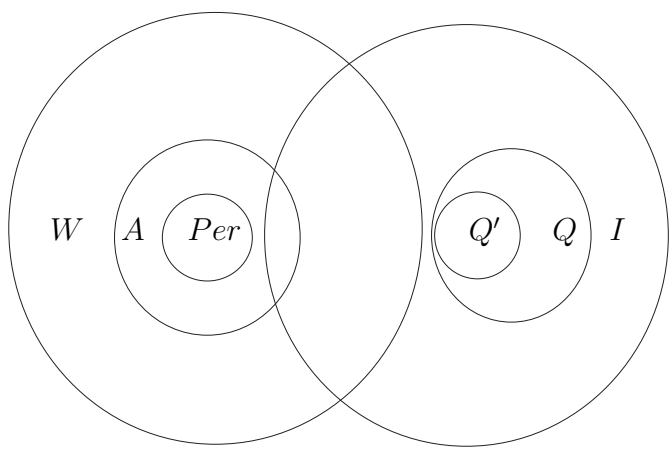}}}
   \put (44,34.5){$W$} \put (53.3,34.5){$A  $}
    \put (54,22.5){$ $}
\put (62,34.5){$Per  $} \put (117.5,34.5){${\large I} $} \put (110,34.5){${\large Q} $}\put (99,34.5){${\large Q'} $}
  \end{picture}
  \caption{$I(T)\cup W(T)$.}
  \label{RecurrenceIRRegular-2015}
  \end{center}
  \end{figure}

\subsection{Main Results}

 Now we start to state our main theorems.
   We need  two conditions called {\it $g-$almost product property}  and  {\it uniform  separation property}   which are introduced in \cite{PS} and  we will  recall them  later in Section  \ref{section-preliminaries}.

 \begin{Thm}\label{More-general-systems}
 Let  $\,T$ be a      continuous map of a compact metric space $X$ with $g-$almost product property and uniform  separation property. If the periodic points are dense in $X$ (i.e., $\overline{Per(f)}=X$) and the periodic measures are dense in the space of invariant measures (i.e., $\overline{ M_{p} (T,X) }=M (T,X)),$ then $\{A(T)\cup R(T), QR(T), W(T), V(T), QW(T), I(T)\}$ has full entropy gaps with respect to $X$.
 \end{Thm}

Theorem \ref{More-general-systems}  seem to be  required too many  very strong conditions, but they are  satisfied by many examples,
 including all topological mixing subshifts of finite type and all $\beta$-shifts, etc. (see  Section \ref{section-examples}).

\bigskip

Given a continuous function $\phi:X\rightarrow \mathbb{R}$, define $\text{  $\phi$-regular set} $ $$R_\phi(T):=\{x\in X\,| \text{ Birkhoff averages } \frac1n \sum_{i=0}^{n-1} \phi(T^i(x))   \text{ converge as } n\rightarrow +\infty\}.$$
 Define the $\phi$-irregular set $I_\phi(T)=X\setminus R_\phi(T).$ These two sets describe different asymptotic behavior under the observation of continuous functions.

 \begin{Thm}\label{More-general-systems-2}
 Let  $\,T$ be a      continuous map of a compact metric space $X$ with $g-$almost product property and uniform  separation property. If the periodic points are dense in $X$ (i.e., $\overline{Per(f)}=X$) and the periodic measures are dense in the space of invariant measures (i.e., $\overline{ M_{p} (T,X) }=M (T,X)),$ then for any continuous function $\phi:X\rightarrow \mathbb{R}$,  \\
 (1)  $\{A(T)\cup R(T), QR(T), W(T), V(T), QW(T), I(T)\}$ has full entropy gaps with respect to $R_\phi(T)$;\\
 (2) if $I_\phi(T)\neq \emptyset$, then    $\{ QR(T), W(T), V(T), QW(T), I(T)\}$  has full entropy gaps with respect to $I_\phi(T)$.
 \end{Thm}

Theorem \ref{More-general-systems} and Theorem \ref{More-general-systems-2} refine  prior multifractal results showing that $I(T)$ (or nonempty $I_\phi(T)$)  carries full topological entropy under certain conditions, for example, see  \cite{Pesin-Pitskel1984,Barreira-Schmeling2000,Bar,TV,To2010,Tho2012}.

For any continuous function $\phi:X\rightarrow \mathbb{R}$ and any $a\in\mathbb{R},$   define $\phi-$regular level set
 $$R_{\phi,a} (T):=\{x\in X|\,\,\lim_{n\rightarrow \infty}\frac1n\sum_{i=0}^{n-1}\phi (T^i (x))=a\}.$$ Remark that $$ R_\phi (T)=\bigsqcup_{a\in\mathbb{R}} R_{\phi,a} (T),$$ where $\sqcup$ denotes disjoint union.
$\phi-$regular level sets  refine the  asymptotic behavior of $ R_\phi (T).$  Define the domain of the multifractal spectrum for Birkhoff averages of $\phi,$ $$L_\phi:=[\inf\{\int\phi d\mu|\,\mu\in M (T,X)\},\,\sup\{\int\phi d\mu|\,\mu\in M (T,X)\}].$$  Let $Int (L_\phi)$ denote the interior of $L_\phi.$ That is, $$Int (L_\phi)= (\inf\{\int\phi d\mu|\,\mu\in M (T,X)\},\,\sup\{\int\phi d\mu|\,\mu\in M (T,X)\}).$$
Remark that if $I_\phi (T)\neq \emptyset,$  $Int (L_\phi)$ is a nonempty open interval,  see (\ref{IrRegular-phi-equivalent}) below.

\begin{Thm}\label{More-general-systems-3}
 Let  $\,T$ be a      continuous map of a compact metric space $X$ with $g-$almost product property and uniform  separation property. If the periodic points are dense in $X$ (i.e., $\overline{Per(f)}=X$) and the periodic measures are dense in the space of invariant measures (i.e., $\overline{ M_{p} (T,X) }=M (T,X)),$ then for any continuous function $\phi:X\rightarrow \mathbb{R}$ satisfying $I_\phi(T)\neq \emptyset$ and for any $a\in Int (L_\phi),$
    $$\{ A(T)\cup R(T), QR(T), W(T), V(T), QW(T), I(T)\}$$ has full entropy gaps with respect to $R_{\phi,a}(T)$.
 \end{Thm}

Theorem \ref{More-general-systems-3} refine  prior multifractal results  on   Birkhoff averages, for example, see \cite{TV,Ol,Barreira-Schmeling2001,Fan-Feng_Wu2001} etc.

\bigskip


Remark that Theorem \ref{More-general-systems-2}  (1) implies   Theorem \ref{More-general-systems}, since $R_\phi(T)=X$ if taking $\phi$ to be a constant function.
  We  will  prove Theorem
   \ref{More-general-systems-3} in Section \ref{Birkhoff-level-set}
 and then use 
 it  to prove Theorem \ref{More-general-systems-2} (1). We
    will prove  Theorem \ref{More-general-systems-2}  (2) in Section \ref{proof-main}.
In general the fundamental idea  for proving such results  is to construct lots of irregular points such that  you can use `so
many points'  to prove full topological entropy. However, the proof details are too long and  in the new requirements, we need to construct lots of different kind irregular points satisfying various  additional conditions.   In order to avoid writing  plenty of original techniques and proof details, inspired the results  from    \cite{PS} by Pfister
     and Sullivan, we prefer to  use  variational principle to solve the problem directly. Precise proof details will appear in Section \ref{proofs}.

\bigskip

Remark that the assumption of   density of periodic points can be replaced by existence of an invariant measure with full support.


\begin{Thm}\label{More-general-systems-201522222}
 Let  $\,T$ be a      continuous map of a compact metric space $X$ with $g-$almost product property and uniform  separation property. If  the periodic measures are dense in the space of invariant measures (i.e., $\overline{ M_{p} (T,X) }=M (T,X))$ and there exists an invariant measure with full support,  then  all results of Theorem \ref{More-general-systems}, Theorem \ref{More-general-systems-2} and Theorem  \ref{More-general-systems-3} hold.
 \end{Thm}

Let us explain why Theorem \ref{More-general-systems-201522222} holds. Under the assumption of density of periodic measures, $\overline{ Per(T) }=X \Leftrightarrow \exists \mu \in M(T,X),\,\,S_\mu=X$
 (see Proposition \ref{Prop-support-same-as-dense-periodic} below). So the assumptions of Theorem \ref{More-general-systems-201522222} are equivalent to the ones of Theorem \ref{More-general-systems}, Theorem \ref{More-general-systems-2} and Theorem  \ref{More-general-systems-3}  and thus Theorem \ref{More-general-systems-201522222} is valid.

\bigskip

Recall that from \cite{PS} $g-$almost product property is weaker than specification property  and  from \cite{DGS}, we know that for any dynamical system with Bowen's specification property, the periodic points are dense in $X$ (i.e., $\overline{Per(f)}=X$, see Proposition 21.3  \cite{DGS}) and the periodic measures are dense in the space of invariant measures (i.e., $\overline{ M_{p} (T,X) }=M (T,X),  \text{ see Proposition 21.8  \cite{DGS} or see \cite{Sig} }).$ So we have a following result as a consequence of Theorem \ref{More-general-systems}, Theorem \ref{More-general-systems-2} and Theorem  \ref{More-general-systems-3}.

\begin{Thm}\label{Thm-specification1}
 Let  $\,T$ be a      continuous map of a compact metric space $X$ with Bowen's specification  property and uniform  separation property. Then all results of Theorem \ref{More-general-systems}, Theorem \ref{More-general-systems-2} and Theorem  \ref{More-general-systems-3} hold.
 \end{Thm}


Moreover, we have a following result when the system is expansive.
\begin{Thm}\label{Thm-specification1-expansive}
 Let  $\,T$ be an    expansive   continuous map of a compact metric space $X$ with Bowen's specification  property. Then all results of Theorem \ref{More-general-systems}, Theorem \ref{More-general-systems-2} and Theorem  \ref{More-general-systems-3} hold. Moreover, Theorem \ref{More-general-systems} and  Theorem \ref{More-general-systems-2} (1)  can be stated for $$\{ A(T), A(T)\cup R(T), QR(T), W(T), V(T), QW(T), I(T)\}.$$
 \end{Thm}

 Let us explain why Theorem \ref{Thm-specification1-expansive} holds.  It is known that expansiveness  is stronger than uniform separation property, see \cite{PS}.
By Theorem \ref{Thm-specification1},  Theorem \ref{Thm-specification1-expansive} is valid except   $R(T)\setminus A(T).$
  Recall  from \cite{DGS} (see Chapter 22) that for any system with Bowen's specification  property and expansiveness, it has
  unique maximal entropy measure and this measure has full support.   So    $R(T)\setminus A(T)$ having full entropy   can be deduced from
   Theorem \ref{T11-section1} below which shows that for any dynamical  system,
 if   there is an ergodic measure with  maximal entropy and  non-minimal  support, then
 $R(T)\setminus A(T) $  has full entropy. Remark that   $R(T)\setminus A(T)\subseteq R(T)\subseteq QR(T)\subseteq  R_\phi(T)$ for any continuous function $\phi.$   So  $\{A(T), R(T)\}$ has full entropy gaps with respect to $X$ and $R_\phi(T)$.


\subsection{Applications to Standard Examples}\label{section-examples}

\subsubsection{Mixing subshifts of finite type}\label{section-subshift}

Recall from \cite{DGS} (Proposition 21.2) any topological mixing subshift of finite type satisfies Bowen's specification. As a subsystem of full shift, it is expansive. So by Theorem \ref{Thm-specification1-expansive}, we have

 \begin{Thm}\label{Thm-subshifts}
 Let  $\,T$ be a     topological mixing subshift of finite type.  Then all results of Theorem \ref{Thm-specification1-expansive}
 hold.
 \end{Thm}

\subsubsection{$\beta-$shifts}\label{section-beta-shift}

Let us recall the definition of $\beta$-shift ($\beta>1$) 
 in   \cite{Walter} (Chapter 7.3). If $\beta\geq 2$ is an integer, $\beta$-shift is the full shift of $\beta$ symbols. So we only need to recall the definition in the case that $\beta$ is not an integer.  Consider the expansion of 1 in powers of $\beta^{-1},$ i.e. $1=\sum_{n=1}^\infty a_n\beta^{-n}$ where $a_1=[\beta]$ and $a_n=[\beta^n-\sum_{i=1}^{n-1}a_i\beta^{n-i}].$ Here $[t]$ denotes the integral part of $t\in \mathbb{R}.$ Let $k=[\beta]+1.$ Then $0\leq a_n\leq k-1$ for all $n$ so we can consider $a=\{a_n\}_1^\infty$ as a point in the space $X=\prod^{+\infty}_{n=1}Y$ where $Y=\{0,1,\cdots,k-1\}$. Consider the lexicographical ordering on $X,$ i.e. $x=\{x_n\}_1^\infty<y=\{y_n\}_1^\infty$ if $x_j<y_j$ for the smallest $j$ with $x_j\neq y_j.$ Let $f:X\rightarrow X$ denote the one-sided shift transformation. Note that $f^na\leq a$ for all $n\geq 0.$  Let $$\Sigma_\beta:=\{x=\{x_n\}_1^\infty|\,x\in X\,\text{ and }\, f^n (x)\leq a\, \text{ for all }\,n\geq 0\}.$$  Then $\Sigma_\beta$ is a closed subset of $X$ and $f (\Sigma_\beta)=\Sigma_\beta.$ Let $\sigma_\beta:=f|_{\Sigma_\beta}.$ Then $ ( \Sigma_\beta, \sigma_\beta)$ is one-sided $\beta-$shift. One can obtain the two-sided $\beta-$shift by letting $$\hat{\Sigma}_\beta:=\{x=\{x_n\}_{-\infty}^\infty|\,x\in \prod^{+\infty}_{n=-\infty}Y\,\text{ and }\, (x_i,x_{i+1},\cdots)\in \Sigma_{\beta}\, \text{ for all }\,i\in \mathbb{Z}\}.$$ Then $\hat{\Sigma}_\beta$ is a closed subspace of $\prod^{+\infty}_{n=-\infty}Y$ invariant under the two-sided shift $$\hat{f}:\prod^{+\infty}_{n=-\infty}Y\rightarrow \prod^{+\infty}_{n=-\infty}Y.$$
The topological entropy of $\beta-$shift ($\beta>1$) is $\log \beta.$
 Remark that by Variational Principle, there is an ergodic measure with positive entropy.  Note that the Dirac measure supported on the fixed point $x=\{0\}_1^\infty\in\Sigma_\beta $ has zero entropy. So every $\beta-$shift is obviously not uniquely ergodic.

 \begin{Thm}\label{Thm-beta-shifts}
Every $\beta-$shift  ($\beta>1$)  satisfies all results of Theorem \ref{Thm-specification1-expansive}.
 \end{Thm}

Let us explain   why Theorem \ref{Thm-beta-shifts} holds.
 By definition     every $\beta-$shift ($\beta>1$)   is a subsystem of full shift on $[\beta]+1$ symbols and so every $\beta-$shift is expansive
 (which is stronger than uniform separation property, see \cite{PS})
  and satisfies $g-$almost product property from   \cite{PS} (see the Example on P.934).  It is known  that  the unique maximal entropy measure  of $\beta-$shifts always carries full support (see \cite{Walters78},  Theorem 13 (ii)).
   Furthermore, it was proved in \cite{Sig-76} that  the periodic measures are dense in the space of invariant measures (i.e., $\overline{ M_{p} (T,X) }=M (T,X)).$
 So the hypotheses of   Theorem \ref{More-general-systems-201522222} hold for all  $\beta-$shifts and then Theorem \ref{Thm-beta-shifts} is obtained except  $R(T)\setminus A(T)$.
 The set  $R(T)\setminus A(T)$   can be deduced from
   Theorem \ref{T11-section1}  below which shows that for any dynamical  system,
 if   there is an ergodic measure with maximal entropy and non-minimal  support, then
 $R(T)\setminus A(T) $  has full entropy. Remark that   $R(T)\setminus A(T)\subseteq R(T)\subseteq QR(T)\subseteq  R_\phi(T)$ for any continuous function $\phi.$   So  $\{A(T), R(T)\}$ has full entropy gaps with respect to $X$ and $R_\phi(T)$.  So Theorem \ref{Thm-beta-shifts} is valid.  In particular, we point out that one  can not use Theorem \ref{Thm-specification1-expansive} to prove  Theorem \ref{Thm-beta-shifts}, since  from   \cite{Buzzi} the set of parameters of $\beta$ for which Bowen's specification holds, is dense in $ (1,+\infty)$ but has Lebesgue zero measure.

In particular, for full shifts on finite symbols, we have a following result that contains more observation.

 \begin{Thm}\label{Thm-fullshifts}
 Let  $\,T$ be a     full shift on $k$ symbols ($k\geq 2$).  Then   Theorem \ref{More-general-systems}  and Theorem \ref{More-general-systems-2} (1)
   can be stated  for $$\{Per(T), A(T),  A(T)\cup R(T), QR(T), W(T), V(T), QW(T), I(T)\}.$$
 \end{Thm}

 Let us explain why Theorem \ref{Thm-fullshifts} holds. By Theorem \ref{Thm-beta-shifts}, one only needs to consider $  A (T)\setminus Per (T)$. This can be deduced from  Theorem \ref{T-almost-periodic} below which shows that for full shifts of finite symbols,  $ R(T)\cap  A (T)\setminus Per (T)$ has full entropy. Since   $R(T)\cap A(T) \setminus Per(T)\subseteq   QR(T)\cap  A (T)\setminus Per (T)\subseteq  R_\phi(T)\cap  A (T)\setminus Per (T)$ for any continuous function $\phi.$   So  $\{Per(T), A(T)\}$ has full entropy gaps with respect to $X$ and $R_\phi(T)$.

\bigskip

\subsubsection {Hyperbolic systems and Lyapunov exponents}\label{Rem-manysystemshold}

From the classical uniform hyperbolicity theory, every subsystem restricted on a
  topological mixing locally maximal hyperbolic    set (called basic set or elementary set)
satisfies specification property (for example, see \cite{Sig}) and satisfies expansiveness. So by Theorem \ref{Thm-specification1-expansive}, we have

 \begin{Thm}\label{Thm-hyperbolic}  Let $f:M\rightarrow M$ be a $C^1$  
  diffeomorphism of a compact Riemanian manifold $M$.
 Let  $\,T$ be a   subsystem   restricted on a  topological mixing locally maximal hyperbolic    set.   Then all results of Theorem \ref{Thm-specification1-expansive}
 hold.
 \end{Thm}
In particular, this result can be applicative to all topological mixing Anosov diffeomorphisms.
Remark that similar results can be stated for topological mixing expanding maps (for example, $T:S^1\rightarrow S^1, x\mapsto k x\,\,mod\, \,1$ for an integer $k\geq 2$).

Moreover, we can use Lyapunov exponents to observe the periodic-like recurrence. Let $f:M\rightarrow M$ be a $C^1$ diffeomorphism of a compact Riemanian manifold $M$. Let $E\subset TM$ be a $Df-$invariant subbundle. Define the (maximal) Lyapunov exponent of $E$ at a point $x\in M$ by $$\lim_{n\rightarrow \infty} \frac 1n\log \|Df^n|_{E_x}\|$$ if the limit exists. Such points are called {\it Lyapunov-regular.} Otherwise, the points are called  {\it Lyapunov-irregular.}  Denote the sets of all Lyapunov-regular and Lyapunov-irregular points by $R_{Lya} (f)$ and $I_{Lya} (f)$ respectively. Let $\Lambda\subseteq M$ be a compact invariant set.   The subbundle $E$  is called {\it conformal } on $\Lambda,$  if for any $x\in \Lambda,\,n\geq 1$ $$ \|Df^n|_{E_x}\|=\prod_{j=0}^{n-1}\|Df|_{E_{f^j (x)}}\|.$$ Let $\phi (x)=\log \|Df|_{E_x}\|$ , if $E$ is a continuous subbundle then it is a continuous function. Thus, for subsystem $T=f|_\Lambda,$ using this $\phi$ in Theorem \ref{More-general-systems-2} and  Theorem \ref{More-general-systems-3} we have

\begin{Thm}\label{Main-Thm-Lyapunov-exponents} Let $f:M\rightarrow M$ be a $C^1$  
 diffeomorphism of a
compact Riemanian manifold $M$.  Let $T:\Lambda \rightarrow \Lambda$ be a  subsystem   restricted on a
topological mixing locally maximal hyperbolic    set  $\Lambda$  and let $E\subset T_\Lambda M$ be a continuous conformal $Df-$invariant
subbundle on $\Lambda$. Then Theorem \ref{More-general-systems-2} and  Theorem \ref{More-general-systems-3} hold  for the function $\phi(x)=\log\|Df|_{E(x)}\|$
(replacing $R_{\phi} (T)$ and $I_{\phi} (T)$ by $R_{Lya} (T)$ and $I_{Lya} (T)$ respectively).
\end{Thm}

In other words, we can use the `eyes' of Lyapunov exponents to  distinguish different `periodic-like' recurrence. Remark that similar results can be stated for topological mixing conformal  expanding maps.

 \subsection{An answer for  Zhou and Feng's question}\label{zhoufeng-OpenProblem}

 There is an open problem in   \cite{ZF} by Zhou and Feng that whether the set

 {\it  $$\{QW (T)\setminus W (T)|\exists\, \mu\in M_x (T) \,s.t. \,S_\mu=C_x\}\neq \,\emptyset\,\,?$$}
This is the set $V\setminus W(T)$ according to the definition of $V$  at the beginning of the paper. It has been solved positively by constructing  examples, see   \cite{OS,HYZ,WHH} etc. From  Theorem \ref{More-general-systems}, for a certain class of dynamical systems (including  topological mixing  subshifts of finite type, all $\beta-$shifts,   systems restricted on mixing locally maximal hyperbolic  sets), $V\setminus W(T)$ is not only nonempty but also has full topological entropy (and so does its  complementary  set in  $QW (T)\setminus W (T)$).  In other words,   $V\setminus W(T)$  has  very strong dynamical complexity which reaches the complexity of dynamical system itself. In particular, we know that positive topological entropy implies $V\setminus W(T)$ has  uncountable elements. So our  Theorem \ref{More-general-systems} can be as a  strong answer for Zhou and Feng's open problem, provided that the given system has positive entropy.




\subsection{ Layout of the Paper}

The remainder of this paper is organized as follows.  In Section \ref{section-preliminaries}
we will recall the notions of entropy, $g$-almost product property, uniform separation and recall
some classical  results including saturated property and entropy-dense property.
In Section \ref{section-20150108-recurrence} we will recall the notions of various
`perioidic-like' recurrence and introduce some simple observation. In Section \ref{section-regular----irregular}
we will recall the notions of regularity and irregularity and introduce some simple facts. In
Section \ref{section-overlap-regular-recurrence} we will coordinate `perioidic-like' recurrence and (ir)regularity
and give some basic  discussion.  In Section \ref{section-usefulfacts-lemmas} we recall and introduce some useful
facts and lemmas.
 In Section \ref{proofs} we divide our main theorems  into several propositions to prove and in this process we state every proposition for possibly  applicative to  more general dynamical systems,   in particular some results can be applied for time-$t$ maps of mixing hyperbolic flows. In Section \ref{section-transitive} we give some similar results involving the set of transitive points. Finally, in Section \ref{section-geometric-Bowen-specification} for generic systems or  systems with Bowen's specification, we give some topological or geometric characterization of various subsets with distinct asymptotic behavior.

\section{Entropy, $g-$almost product property \& uniform  separation}\label{section-preliminaries}

\subsection{Entropy}


Let $T:X\rightarrow X$  be a continuous map of a compact
metric space $X$.
Now let us to recall the definition of topological entropy  in \cite{Bowen1} by Bowen.

Let $x\in X$. The dynamical ball $B_n (x,\varepsilon)$ is the set $$B_n (x,\varepsilon):=\{y\in X|\,\max\{d (T^j (x),T^j (y))|\,0\leq j\leq n-1\}\leq \varepsilon\}.$$
Let $E\subseteq X,$ and $\mathfrak{F}_n (E,\epsilon)$ be the collection of all finite or countable covers of $E$ by sets of the form $B_m (x,\epsilon)$ with $m\geq n.$ We set $$C (E;t,n,\epsilon,T):=\inf\{\sum_{B_m (x,\epsilon)\in \mathcal{C}}2^{-tm}:\mathcal{C}\in \mathfrak{F}_n (E,\epsilon)\},$$ and $$C (E;t,\epsilon,T):=\lim_{n\rightarrow \infty}C (E;t,n,\epsilon,T).$$ Then $$h_{top} (E,\epsilon,T):=\inf\{t:C (E;t,\epsilon,T)=0\}=\sup\{t: C (E;t,\epsilon,T)=\infty\}$$ and the {\it topological entropy} of $E$ is defined as $$h_{top} (T,E);=\lim_{\epsilon\rightarrow 0}h_{top} (E,\epsilon,T).$$ In particular, if $E=X,$ we also denote $h_{top} (T,X)$ by $h_{top} (T)$.
It is known   from   \cite{Bowen1} that if $E$ is an invariant compact subset, then the topological entropy $h_{top} (T,E)$ is same as the classical definition (for classical definition of topological entropy, see Chapter 7 in \cite{Walter}).

Let us recall some basic facts about topological entropy.
From   \cite{Bowen1}
for any subsets $Y_1\subseteq Y_2\subseteq X,$
\begin{eqnarray}\label{Y1containY2}
h_{top} (T,Y_1)\leq h_{top} (T,Y_2).
\end{eqnarray}
 If one considers a collection of subsets of $X$:   $\{Y_i\}_{i=1}^n\,(n\in \mathbb{N}\cup\{+\infty\}),$    from   \cite{Bowen1} we know that    the topological entropy satisfies \begin{eqnarray}\label{Bowen-disjointSet-entropy}
 h_{top} (T,\cup_{i=1}^nY_i)=\sup_{1\leq i\leq n} h_{top} (T,Y_i).
 \end{eqnarray}

 Let $M(X)$ denote the space of all Borel  probability  measures supported on $X$.
 Let $\xi=\{V_i|\,i=1,2,\cdots,k\},$ be a finite partition of measurable sets of $X$. The entropy of $\nu\in M(X)$ with respect to $\xi$
 is $$H(\nu,\xi):=-\sum_{V_i\in\xi}\nu(V_i)\log \nu(V_i).$$
 We write $T^{\vee n}\xi:=\vee_{k\in \Lambda }T^{-k}\xi.$ The entropy of $\nu\in M(T,X)$ with respect to $\xi$ is $$h(T,\nu,\xi):=\lim_{n\rightarrow \infty}\frac 1n H(\nu, T^{\vee n}\xi),$$ and the {\it metric entropy} of $\nu$ is $$h_\nu(T):=\sup_{\xi} h(T,\nu,\xi).$$
 More information of metric entropy, see Chapter 4 of \cite{Walter}.

Let us recall some relations on metric entropy and topological entropy.
By classical  Variational Principle (Theorem 8.6 and Corollary 8.6.1 in \cite{Walter}), we know that \begin{eqnarray}\label{eq-classical-Variational-Principle}
 h_{top} (T)=\sup_{\mu\in M (T,X)}  h_\mu (T)
=\sup_{\mu\in M_{erg} (T,X)}  h_\mu (T).
\end{eqnarray}
 From Theorem 1 in \cite{Bowen1}
for any ergodic measure $\mu$ and any subset $Z\subseteq X$, if $\mu (Z)=1$, then
\begin{eqnarray}\label{FullMeasure-larger}
h_\mu (T)\leq h_{top} (T,Z).
\end{eqnarray} Moreover,  every set with totally full
 measure (i.e., being full measure for all invariant measures) carries  full topological entropy, that is,

 \begin{Thm}\label{thm-fullmeasure-fullentropy}
 Let  $\,T$ be a      continuous map of a compact metric space $X$. If $Y\subseteq X$ is a set with totally full measure, then  \begin{eqnarray}\label{FullMeasure-FullEntropy}
 h_{top} (T,Y) =h_{top} (T).\end{eqnarray}
 \end{Thm}


  {\bf Proof.}
   Let $Y\subseteq X$ be a set with totally full measure.
 By (\ref{Y1containY2}), (\ref{eq-classical-Variational-Principle}) and (\ref{FullMeasure-larger}), \begin{eqnarray*}\label{FullMeasure-FullEntropy-proof}
 h_{top} (T,Y)&\leq& h_{top} (T,X)=h_{top} (T)=\sup_{\mu\in M (T,X)}  h_\mu (T)\nonumber\\
&=&\sup_{\mu\in M_{erg} (T,X)}  h_\mu (T)\leq h_{top} (T,Y).\end{eqnarray*}
This means that $Y$ carries full entropy. \qed


\subsection{$g-$almost product property}

Firstly we recall  the definition of   specification property which is stronger than  $g-$almost product property, see \cite{DGS,Sig,Bow,Bowen2,Bowen71-trans,To2010}. Let  $\,T$ be a      continuous map of a compact metric space $X$.

\begin{Def}\label{specification}  We say that the dynamical system $T$  satisfies {\it specification property}, if the following holds:  for any $\epsilon>0$ there exists an integer $M(\epsilon)$ such that for any $k
\geq 2,$ any $k$ points $x_1,\cdots,x_k$, any  integers $$a_1\leq b_1<a_2\leq b_2\cdots<a_k\leq b_k$$ with $a_{i+1}-b_i\geq M(\epsilon)\,(2\leq i\leq k),$   there exists a point $x\in X$ such that \begin{eqnarray}\label{specification-inequality}
 d(T^j(x),T^j(x_i))<\epsilon,\,\,\,\,for \,\,a_i\leq j\leq b_i,\,1\leq i\leq k.\end{eqnarray}

\end{Def}

The original definition of specification, due to Bowen, was stronger.

\begin{Def}\label{Bowen-specification} We say that the dynamical system $T$  satisfies {\it Bowen's  specification property}, if  under the assumptions of Definition \ref{specification} and   for any  integer $p\geq M(\epsilon)+b_k-a_1,$ there exists a point $x\in X$ with $T^p(x)=x$ satisfying (\ref{specification-inequality}).

\end{Def}

Now we start to recall the concept $g-$almost product property in \cite{PS}  (there is a slightly weaker variant,  called almost specification, see \cite{Tho2012}).  It is weaker than specification property(see Proposition 2.1 in \cite{PS}).  A striking and typical example of   $g-$almost product property (and almost specification) is that it applies to every $\beta-$shift \cite{PS,Tho2012}. In sharp contrast, the set of $\beta$ for which
the $\beta-$shift has specification property has zero Lebesgue measure  \cite{Buzzi,Schmeling}.

Let $\Lambda_n=\{0,1,2,\cdots,n-1\}.$ The cardinality of a finite set $\Lambda$ is denoted by $\# \Lambda.$ Let $x\in X$. The dynamical ball $B_n(x,\varepsilon)$ is the set $$B_n(x,\varepsilon):=\{y\in X|\,\max\{d(T^j(x),T^j(y))|\,j\in\Lambda_n\}\leq \varepsilon\}.$$

\begin{Def}\label{blowup-function} Let $g:\mathbb{N}\rightarrow \mathbb{N}$  be a given nondecreasing unbounded map with the properties $$g(n)<n\,\,\text{ and } \lim_{n\rightarrow \infty}\frac{g(n)}n=0.$$ The function $g$ is called {\it blowup function.} Let $x\in X$ and $\varepsilon>0.$ The $g-$blowup of $B_n(x,\varepsilon)$ is the closed set
$$B_n(g;x,\varepsilon):=\{y\in X|\, \exists \Lambda\subseteq\Lambda_n ,\#(\Lambda_n\setminus\Lambda)\leq g(n)\,\text{ and }\,\max\{d(T^j(x),T^j(y))|\,j\in\Lambda\}\leq \varepsilon\}.$$

\end{Def}

\begin{Def}\label{product-property} We say that the dynamical system $T$  satisfies {\it $g-$almost product property}  with blowup function $g$, if  there is a nonincreasing function $m:\mathbb{R}^+\rightarrow \mathbb{N},$ such that
for any $k
\geq 2,$ any $k$ points $x_1,\cdots,x_k\in X$, any positive $\varepsilon_1,\cdots,\varepsilon_k$ and any integers $n_1\geq m(\varepsilon_1),\cdots, n_k\geq m(\varepsilon_k),$ $$\bigcap_{j=1}^k T^{-M_{j-1}}B_{n_j}(g;x_j,\varepsilon_j)\neq \emptyset,$$ where $M_0:=0,M_i:=n_1+\cdots+n_i,i=1,2,\cdots,k-1.$

\end{Def}

\subsection{Uniform  separation}

Now we recall the definition of uniform separation property \cite{PS}. For $x\in X,$ define $$\Upsilon_n (x):=\frac1n\sum_{j=0}^{n-1}\delta_{T^j (x)}$$ where
  $\delta_y$ is the Dirac  probability measure supported
  at $y\in X$.   
For $\delta>0$ and  $\varepsilon>0$, two points $x$ and $y$ are
$(\delta,n,\varepsilon)-$separated if $$\#\{j:d(T^jx,T^jy)>\varepsilon,\,j\in\Lambda_n\} \geq \delta n.$$ A subset
$E$ is  $(\delta,n,\varepsilon)-$separated if any pair of different points of $E$ are  $(\delta,n,\varepsilon)-$separated.      Let $F\subseteq M(X)$ be a neighborhood of $\nu\in M(T,X)$.
 Define $$ X_{n,F}:=\{x\in X|\, \Upsilon_n(x)\in F\},$$
and define
 \begin{eqnarray*}N(F;\delta,n,\varepsilon):=\text{maximal cardinality of a } (\delta, n,\varepsilon)-\text{separated subset of } X_{n,F}.\end{eqnarray*}

\begin{Def}\label{uniform-separation-property} We say that the dynamical system $T$  satisfies {\it uniform separation  property}, if  following holds. For any $\eta>0,$ there exist $\delta^*>0,\epsilon^*>0$ such  that for $\mu$ ergodic and any neighborhood $F\subseteq M(X)$ of $\mu$, there exists $n^*_{F,\mu,\eta},$ such that for $n\geq n^*_{F,\mu,\eta},$ $$N(F;\delta^*,n,\epsilon^*)\geq 2^{n(h_\mu(f)-\eta)}.$$

\end{Def}

Now let us recall a basic relation between expansiveness and uniform separation in \cite{PS}.
\begin{Thm}\label{Thm-expansive} (Theorem 3.1 in \cite{PS}) Let  $\,T$ be a      continuous map of a compact metric space $X$. If $\,T$ is expansive (or asymptotically $h$-expansive),   $T$ satisfies uniform separation. \end{Thm}


\subsection{Variational Principle for saturated sets}

Let  $\,T$ be a      continuous map of a compact metric space $X$. Recall the definition that  $\Upsilon_n (x):=\frac1n\sum_{j=0}^{n-1}\delta_{T^j (x)}$.  Let $M_x(T)$ be the limit set of the empirical measures for $x$, i.e., the set of all limits of $\Upsilon_n(x)$
 in weak$^*$ topology.

Now we recall  a result from   \cite{PS}.
   The system $T$ is said to be {\it saturated} (or $T$ has saturated property), if  for any  compact connected nonempty set $K \subseteq M (T,X ),$
$$h_{top} (T,G_K)=\inf\{h_\mu (T)\,|\,\mu\in K\},$$ where $G_K=\{x\in X|\,M_x (T)=K\}.$

\begin{Lem}\label{lem-PS} (Variational Principle, Theorem 1.1 in \cite{PS}) \\
 Let  $\,T$ be a      continuous map of a compact metric space $X$ with $g-$almost product property and uniform  separation property.
Then $T$ is saturated.



\end{Lem}

On the other hand, from \cite{PS} if one does not have uniform separation property, then the saturated property just holds for any singleton $K$.  For convience to compare saturated property, we give a following notion called single-saturatd property. We say $T$ is {\it single-saturated}, if  $h_{top} (T,G_\mu)= h_\mu (T) $ holds for any $\mu\in M(T,X),$ where $G_\mu=\{x\in X|\,M_x (T)=\{\mu\}\}.$

\begin{Lem}\label{lem-PS-etds2007-no-uniform-separation} (Variational Principle, Theorem 1.2 in \cite{PS}) \\
 Let  $\,T$ be a      continuous map of a compact metric space $X$ with $g-$almost product property.
Then $T$ is single-saturated.



\end{Lem}

Remark that for any continuous map $\,T$ of a compact metric space $X$, there is a general fact (see Theorem 4.1 (3) in \cite{PS}):
for any  compact connected nonempty set $K \subseteq M (X, T ),$
\begin{eqnarray}\label{eq-generalentropyestimate-cptconnected-K}
h_{top} (T,G_K)\leq \inf\{h_\mu (T)\,|\,\mu\in K\},\,\, \text { where }  \,\,G_K=\{x\in X|\,M_x (T)=K\}.
\end{eqnarray}
In particular,  for  any $\mu \in  M (X, T ),$ we have
\begin{eqnarray}\label{eq-generalentropyestimate-singlemeasure}
h_{top} (T,G_\mu)\leq  h_\mu (T) ,\,\, \text { where }  \,\,G_\mu=\{x\in X|\,M_x (T)=\{\mu\}\}.
\end{eqnarray}


\subsection{Entropy-dense property}

Now  let's recall the entropy-dense property of Theorem 2.1 in   \cite{PS2005} (or see \cite{PS}, also  see \cite{EKW} for similar discussion).
  Roughly speaking, any invariant probability measure $\mu$ is the  limit of a sequence of ergodic measures
$\{\mu_n\}_{n=1}^{\infty}$ in weak$^*$ topology such that the entropy of $\mu$ is the limit of the entropies of  $\mu_n$. Recall $M (X)$ and $M(T,X)$ denote the space of all Borel probability measures and the space of invariant measures, respectively.
 Here we further  require $S_{\mu_n}\neq X$ in this property (whose statement is a little stronger than usual entropy-dense property).  More precisely, we say $T$ has {\it entropy-dense} property, if for any $\nu\in M (T,X)$, any neighborhood $G\subseteq M (X)$ of $\nu$ and any $ h_* < h_\nu (T),$ there exists an ergodic
measure $\mu\in G\cap M (T,X)$ such that $S_\mu\neq X$ and $h_\mu (T) > h_*.$

\begin{Lem}\label{lem-entropy-dense-Ps-2005} 
Let  $\,T$ be a      continuous map of a compact metric space $X$ with $g-$almost product property. Suppose that $M (T,X)$ is not a singleton. Then $T$ has entropy-dense property.

\end{Lem}

{\bf Proof.}  Here this lemma is   slightly stronger than Theorem 2.1 in   \cite{PS2005},  since we further add $S_\mu\neq X$. Now let us explain this  more precisely. Let  $\nu\in M (T,X)$ and  $G\subseteq M (X)$ be a neighborhood  of $\nu.$
 Since $M (T,X)$ is not a singleton, then  we can take an open ball  $G'\subseteq M(X)$ such that   $\nu\in G'\subseteq  \overline{G'}\subset G$ and $M (T,X)\setminus \overline{G'} \neq \emptyset.$ Note that $M (T,X)\setminus \overline{G'}$ is open in $M (T,X)$.

 From the proof of Proposition 2.3 (1) of   \cite{PS2005}, one construct  a closed invariant set $Y$ and there exists $n_{G'} \in\mathbb{N}$ such that $ h_{top} (T,Y)>h_*$ and for
any $y\in Y$ and any $n\geq n_{G'},$    $\Upsilon_n (y)\in G'.$ So for
 any $m\in M_{erg} (T,Y),$ by Birkhorff ergodic theorem there is $y\in Y$ such that $ \Upsilon_n (y)$ converge to $m$ in weak$^*$ topology and thus $m\in \overline{G'}.$ In other words, $M_{erg} (T,Y)\subseteq \overline{ G'}.$  By convex property of the ball $G'$ and Ergodic Decomposition theorem,  $M (T,Y)\subseteq \overline{ G'}.$
 Then  $Y\neq X,$ since $M (T,X)\setminus \overline{G'} \neq \emptyset.$
   By Variational Principle, for $0<\varepsilon<h_{top} (T,Y)- h_*$,  take a $\mu\in M_{erg} (T,Y)$ such that $h_\mu (T)>h_{top} (T,Y)-\varepsilon>h_*.$ Then $\mu$ is the measure we need.  For more details, see    \cite{PS2005}.
\qed


\section{Periodic and Periodic-like Recurrence}\label{section-20150108-recurrence}
 One important way to partition points with different asymptotic behavior is according to the recurrence property.

 In the classical study of   dynamical systems, an important concept is non-wandering point. A point $x\in X$ is called {\it wandering}, if there is a neighborhood $U$ of $x$ such that the sets $T^{-n}U,\,\,n\geq 0$, are mutually disjoint. Otherwise, $x$ is called {non-wandering.}   Let $\Omega (T)$ denote the set of all non-wandering points, called {\it non-wandering set}.
The interesting action of $T$ takes place in $\Omega (T)$ and recall that  from  Theorem 5.6, Theorem 6.15  and  Corollary 8.6.1 in \cite{Walter} $\Omega (T)$ is always invariant, compact, carries totally full measure and owns the whole complexity of the system $$h_{top} (T)= h_{top} (T,\Omega (T)).$$ So in general  one can always consider the subsystem $T:\Omega (T)\rightarrow \Omega (T)$ to replace the original  system $T:X\rightarrow X$.  
  It is   interesting to ask for general dynamical systems,  how about the dynamical complexity of $X\setminus\Omega (T)$?  Unfortunately, in this paper the systems studied is  the case $X=\Omega (T)$, since we require that the periodic points are dense in $X$ in our main theorems (Theorem \ref{More-general-systems}-\ref{More-general-systems-3}).

The set $\Omega (T)$ consists of those points with a weak recurrence property. Now let us recall the concept of recurrent point.  It is known that recurrent points play important roles in the ergodic theory of dynamical systems. Given $x\in X,$ let $\omega_T (x)$ denote the $\omega$-limit set. We call $x\in X$ to be {\it recurrent}, if $$x\in \omega_T (x).$$
 Let  $Rec (T)$ denote the set of all recurrent points. By definition, it is obvious that $Rec (T)\subseteq \Omega (T).$

\begin{Thm}\label{T1} {\it For any continuous   map $T:X\rightarrow X$    of a compact
metric space $X$, the recurrent set $Rec (T)$ has full topological entropy.}

   \end{Thm}

{\bf Proof.}
It is known that from Poincar$\acute{\text{e}}$ recurrent theorem recurrent set $Rec(T)$ has totally full measure (see Remark on Page 157 of \cite{Walter})  so that by Theorem \ref{thm-fullmeasure-fullentropy}  we complete the proof of Theorem \ref{T1}.  \qed

\bigskip

In the study of   (smooth or topological) dynamical systems, many people pay attention to refine recurrent set according to the `recurrent frequency'.
 A standard and important kind of recurrent point with same asymptotic behavior is periodic point, which returns itself through  finite iterates.
Let $Per (T)$ denote the set of periodic points. Then $$Per (T)\subseteq Rec (T).$$ A fundamental question in dynamical systems is to search the existence of periodic points.  For    system with Bowen's specification (such as topological mixing subshifts of finite type,  topological mixing uniformly hyperbolic and topological mixing expanding systems), it is well known that the set of periodic points is dense in the whole space (for example, see Proposition 21.3 in \cite{DGS}) and moreover, if the system is expansive, then the periodic set  is  countable  and thus has zero entropy.
 Here we give a basic observation on the topological  entropy of periodic set for general dynamical systems, admitting existence of uncountable periodic points.

\begin{Thm}\label{T2} {\it For any continuous   map $T:X\rightarrow X$    of a compact
metric space $X$, the periodic  set $Per (T)$ either is empty or has zero topological entropy.}
\end{Thm}

{\bf Proof.} 
  Notice that $Per (T)=\cup_{n\geq 1} P_n(T),$ where $P_n(T)=\{x\in X|\,T^n(x)=x\}.$ By (\ref{Bowen-disjointSet-entropy}), we only need to show that for any $n\geq 1,$ $P_n(T)$ carries zero entropy. Note that for any fixed $n\geq 1,$  $P_n(T)$ is an invariant closed (compact)  set, every ergodic measure $\mu$ supported on $P_n(T)$ is periodic and so $\mu$ has zero metric entropy.  Applying  classical Variational Principle (Corollary 8.6.1 in \cite{Walter}) for subsystem $T|_{P_n(T)},$ we have
$$ h_{top}(T,P_n(T))=\sup_{\mu\in M_{erg}(T, P_n(f))} h_\mu(f)=0.$$  This ends the proof of Theorem \ref{T2}. \qed

\bigskip

From Theorem \ref{T2}, the periodic set has no dynamical complexity in the sense of topological entropy.
  However, a classical interesting  result states that for expansive systems with Bowen's specification (which implies topological mixing),  the topological entropy can be characterized by the exponential growth of periodic points with same period (for example, see Proposition 22.7 in \cite{DGS}). More precisely, $$h_{top} (T)=\lim_{n\rightarrow \infty}\frac1n \log \# P_n (T)$$ where $P_n (T)=\{x|\,T^n (x)=x\}$ and $\#A$ denotes the cardinality of the set $A $ (more better characterization of entropy and the growth of periodic points, see Proposition 22.6 in \cite{DGS}).  Roughly speaking,  periodic points have  same complexity as the system itself from the viewpoint of  two different ways to  interpret complexity. This result holds for Axiom
A diffeomorphisms in any dimension (see   \cite{Bowen2}), but this is a somewhat special situation. It is well known that Axiom A
diffeomorphisms are not dense in $\Diff^1 (M)$ (the space of all diffeomorphisms on a compact Riemanian manifold $M$). In 2004 Kaloshin   \cite{Kal} showed that in general
$\# P_n (T)$ can grow much faster than entropy. Moreover, it is well known that for $C^1$ generic diffeomorphisms, all periodic points are hyperbolic so that countable and they form a dense subset of the non-wandering set (by  classical  Kupka-Smale theorem, Pugh's  or Ma\~{n}$\acute{\text{e}}$'s   Ergodic Closing lemma from   Smooth Ergodic Theory, for example,  see \cite{Kupka1,Smale1,Pugh1,Pugh2,Mane1}). For nonuniformly hyperbolic systems, many people studied the existence of periodic point by showing closing lemma,  for example, see   \cite{K3}. In particular, an interesting result from  \cite{K3} is that   any surface diffeomorphism with positive entropy always carries a lot of periodic points.    All in all, periodic points have been studied more and more in the research of modern  dynamical systems.


In general, it is well known that there are lots of topological dynamical systems without  periodic points. The standard example is irrational rotation.  So   many  generalizations of periodic points
are introduced.  One such kind `periodic-like' point is almost periodic point. A point $x\in X$ is {\it almost periodic}, if for every open neighborhood $U$ of $x$, there exists $N>0$ such that $f^k (x)\in U$ for some $k\in [n,n+N]$ and every integer $n\geq 1.$ Let $A (T)$ denote the set of all almost periodic points, called almost periodic set for convenience.   It is well-known  from   \cite{Birkhoff}   that a point $x$ is almost periodic if and only if $x$ is minimal, i.e., $x\in\omega_T (x)$ and $\omega_T (x)$ is minimal (see
\cite{Gottschalk44-2222,Gottschalk44,Gottschalk46,Mai} for more related discussion in the sense that  the space $X$ is more  general, not necessarily being compact metric space). Here an invariant  set $E\subseteq X$ is called  minimal, if for every point $y\in E,$
$\omega_T (y)=E.$ In particular, if $X$ is minimal, we say the system  $T$ to be minimal.  Remark that the almost periodic set  $A(T)$ can be written as the union of all minimal sets.


Remark that for any uniquely ergodic system, the support of the unique invariant measure must be minimal. However,  there are minimal invariant sets which are
not uniquely ergodic  \cite{Oxt}. Note that   $E\subseteq X$ is a minimal invariant subset if and only if the support of any invariant measure supported on $E$ coincides with $E$.  By Zorn's lemma, one can show that any dynamical system contains at least one
minimal invariant subset.    So different from periodic points, almost periodic points naturally exist   and thus played important roles in the study of  all topological  dynamical systems. Moreover, constructions of minimal examples are studied a lot by many researchers. For homeomorphisms, there are many examples of subshifts which are strictly ergodic and
has positive entropy   \cite{Gri,HaKa}. For systems on manifolds, from   \cite{Rees} there are minimal homeomorphisms on 2-torus with positive entropy and it was proved in   \cite{BCorL} that any compact manifold of dimension $d \geq 2 $ which carries a minimal uniquely
ergodic homeomorphism also carries a minimal uniquely ergodic homeomorphism with positive topological
entropy.  M. Herman asked whether, for diffeomorphisms, positive topological
entropy was compatible with minimality or strict ergodicity. It was  constructed in   \cite{Herman} a 4-dimensional example of a minimal
 (but not strictly ergodic) diffeomorphism with positive topological entropy. For any $C^{1+\alpha}$ surface diffeomorphism with positive entropy, by classical Pesin theory
 there are  lots of periodic points, see    Corollary 4.4 in   \cite{K3},   so that the system is not minimal. For general  non-minimal systems, we have a following observation.

 \begin{Thm} \label{Thm-22-minimal-supportnot} {\it For any non-minimal continuous   map $T:X\rightarrow X$    of a compact
metric space $X$, if there is an ergodic measure $\mu$ with full support, then $\mu(A(T))=0$.}
\end{Thm}



{\bf Proof.}
    By contradiction,   $\mu (A (T))>0$.  By ergodicity and invariance of $A (T),$ $\mu (A (T))=1.$ Let $S_\mu$ be the support of $\mu$ and by assumption $X=S_\mu$.  By ergodicity and full support of $\mu$, it is known that  the set  $D:=\{x\in X|\,\overline{Orb(x)}=X\}$  has $\mu$ full measure (for example,   see Theorem 1.7 or  Theorem 5.16 in \cite{Walter}).
So  $\mu(A (T)\cap D)=1$ and thus we can take  $z\in A (T)$ such that the orbit of $z$ is dense in   $X.$   $z\in A (T)$ implies that $z\in \omega_T(z)$ and   the closed invariant set $\omega_T (z)$ is minimal. $z\in \omega_T(z)$ implies that $\omega_T(z)=\overline{Orb(z)}$.  Then    $X= \overline{Orb(z)}= \omega_T (z)$  is minimal,    contradicting that $T$ is not minimal.\qed

\bigskip





Notice that for systems such as topological mixing Anosov systems, mixing subshifts of finite type,
the unique measure with maximal entropy is ergodic  and its support is not minimal. Thus, by Theorem \ref{Thm-22-minimal-supportnot}  if one wants  to find some kind of periodic-like points with totally full measure, we need to generalize almost periodic point to be more general.
 Zhou etc. (see   \cite{Zhou93,Zhou95,ZF}) introduced such two more general  concepts of `periodic-like'  points which have totally full measure for all dynamical systems.

One is called {\it weakly almost periodic} and another is more weaker called {\it quasi-weakly almost periodic}.
 Different `recurrent frequency' determines different asymptotic behavior.
 If $E\subseteq X$ is nonempty and $x\in X$, define $$\,\,\,\,\,\,\,\,\underline{P}_x (E):=\liminf_{n\rightarrow \infty}\frac1n \sum_{i=0}^{n-1} \chi_E (T^i (x))
 \,\,\,\,\,\text{and}\,\,\,\overline{P}_x (E):=\limsup_{n\rightarrow \infty}\frac1n \sum_{i=0}^{n-1}  \chi_E (T^i (x)).$$
 In other words, recalling  the definition of $\Upsilon_n (x)=\frac1n\sum_{j=0}^{n-1}\delta_{T^j (x)}$,  $$ \underline{P}_x (E)=\liminf_{n\rightarrow \infty}\Upsilon_n (x) (E),\,\,\,\overline{P}_x (E)=
 \limsup_{n\rightarrow \infty}\Upsilon_n (x) (E).$$
 If $\underline{P}_x (E)=\overline{P}_x (E),$ we denote by $P_x (E),$ which means the probability of the orbit of $x$ enters in $E.$  Let $V_\varepsilon (x)$ denote  $\varepsilon$-neighborhood of $x$, i.e., $ V_\varepsilon (x)=\{y\in M\,|\, d (x,y)<\varepsilon\}$.

 \begin{Def} \label{Weakly-periodic} ({\bf (quasi-)weakly almost periodic}) We call $x$ to be a weakly almost periodic point, if for any $\varepsilon>0$, $$\underline{P}_x (V_\varepsilon (x))>0.$$
$x$ is called to be a quasi-weakly almost periodic point, if for any $\varepsilon>0$, $$\overline{P}_x (V_\varepsilon (x))>0.$$

\end{Def}

  Let $QW (T)$ and $W (T)$ denote the sets of  all weakly almost periodic points and  all
  quasi-weakly almost periodic points, respectively.
   Remark that
  \begin{eqnarray}\label{weakperiodic-in-reccurent}
 \Omega (T)\supseteq Rec (T)\supseteq QW (T)\supseteq W (T)\supseteq A (T)  \supseteq Per (T).
 \end{eqnarray}
 We have known  $\Omega (T), Rec (T)$ both have   full measure
  for any invariant measure. Now we show that   $ QW (T)$ and $ W (T)$ also both have   full measure
  for any invariant measure.

  \begin{Thm}\label{T3} {\it For any continuous   map $T:X\rightarrow X$    of a compact
metric space $X$, $ QW (T)$ and $ W (T)$   both have   full measure
  for any invariant measure.  In particular,  each one of 
$QW (T)$ and $ W (T)$ carries full topological entropy.}

\end{Thm}

  {\bf Proof.}   For any ergodic measure $\mu$,  let $G (\mu)$ denote the set of all points satisfying that $ \Upsilon_n (x) $ converges
  to   $\mu$ in weak$^*$ topology. By Birkhoff ergodic theorem  $G (\mu)$ is of $\mu$ full measure. Let $S_\mu$ denote the support of  $\mu$,
  meaning that $S_\mu$ is the smallest closed  invariant set with $\mu$ full measure.  So  $S_\mu\cap G (\mu)$ is of $\mu$ full measure and every
  $x\in S_\mu\cap G (\mu)$ satisfies $\Upsilon_n (x)\rightarrow \mu$ in weak$^*$ topology. By weak$^*$ topology for open sets (see Remarks (3) (iii) on Page 149 of \cite{Walter}),
  we have  for any $x\in S_\mu\cap G(\mu),\,\varepsilon>0,$
   $$\underline{P}_x (V_\varepsilon (x))=\liminf_{n\rightarrow \infty}\Upsilon_n (x) (V_\varepsilon (x))\geq \mu (V_\varepsilon (x))>0.$$ That is,
   for any ergodic measure $\mu$, $\mu$ a.e. $x$ belongs to $W (T)$. Thus by Ergodic Decomposition theorem,  $W (T)$ has full measure for any invariant measure.
   By (\ref{weakperiodic-in-reccurent})    $QW (T)$ also has full measure for any invariant measure.
    By Theorem \ref{thm-fullmeasure-fullentropy} we complete the proof.  
\qed

Remark that by (\ref{weakperiodic-in-reccurent})
  the proof of Theorem \ref{T3} also can be as  an  alternative proof to show  that $\Omega (T), Rec (T)$ both have   full measure
  for any invariant measure.

\bigskip

There are some equivalent statements of weakly almost periodic point and quasi-weakly almost periodic point.  Firstly we need to recall  some notions (see \cite{Zhou-center-measure}, for the case of  flow  see \cite{N-book}).
A set $E$ is said to be a {\it center of attraction} of $T$ with respect to $X_0\subseteq X$, if
$E$ is a closed invariant set, and for any $x \in X_0$ and $\varepsilon > 0,$ the limit $$P_x (V (E,\varepsilon))=\lim_{n\rightarrow \infty}\frac1n \sum_{i=0}^{n-1}  \chi_{_{V (E,\varepsilon)}} (T^i (x))$$ exists and equals to 1, where $V (E,\varepsilon)$ denotes the
$\varepsilon$-neighborhood of $E$, i.e., $V (E,\varepsilon) = \{y\in X | d (E, y) < \varepsilon\}.$
A set $E$ is called {\it minimal  center of attraction}  relative to $X_0$ if $E$ is a center of attraction relative to $X_0$ but no proper subset of $E$ also is. Denote by
$C (X_0)$ the minimal center of attraction relative to $X_0$. In particular, let
 $$C_x:=C (\{x\}),\,x\in X.$$ Let $M_x (T)$ be the set of all limits of $\Upsilon_n (x)$ in weak$^*$ topology. Recall that $S_\mu$ denotes the support of a measure $\mu$ (that is, the smallest closed subset of $X$ with $\mu$ full measure). Let $M_{X_0}(T)=\cup_{x\in X_0} M_x(T).$
 From   \cite{Zhou-center-measure,Zhou93,Zhou95,ZF}  we know  some basic facts:
\begin{Lem}
Let $T:X\rightarrow X$  be  a  continuous   map  of a compact
metric space $X$.   For any $X_0\subseteq X,\,x\in X,$
\begin{eqnarray}\label{C_x=UnionSupport}
C({X_0})=\overline{\cup_{m\in M_{X_0} (T)}S_m},  \text{ in particular  } \,\,  C_x=\overline{\cup_{m\in M_x (T)}S_m};
\end{eqnarray}
 \begin{eqnarray}\label{C_xinOmega}
\forall\,\, x\in X,\,\,\,C_x\subseteq \omega_T (x);
\end{eqnarray}
For any recurrent point $ x\in Rec (T),$ \begin{eqnarray}\label{WT}
\,x\in W (T)\Leftrightarrow C_x=S_\mu,\,\forall \mu\in M_x (T)\Leftrightarrow \omega_T (x)=S_\mu,\,\forall \mu\in M_x (T);
\end{eqnarray}
\begin{eqnarray}\label{QW}
\,x\in QW (T)\Leftrightarrow x\in C_x\Leftrightarrow \omega_T (x)=C_x\Leftrightarrow\omega_T (x)=\overline{\cup_{\mu\in M_x (T)}S_\mu}.
\end{eqnarray}

\end{Lem}

\bigskip

If $C_x=X,$ then by (\ref{C_xinOmega}) $\omega_T (x)=X\ni x$ so that  $x\in Rec (T)\cap C_x.$ By (\ref{QW}), $x\in QW (T).$  That is, we have  following corollary.
\begin{Cor} Let $T:X\rightarrow X$  be  a  continuous   map  of a compact
metric space $X$. For any $x\in X,$
\begin{eqnarray}\label{C_x=X}
C_x=X\Rightarrow x\in QW (T).
\end{eqnarray}

\end{Cor}

For convenience of readers, we   prefer to introduce some classical equivalent definitions.
Let $S\subseteq \mathbb{N}$, define $$\bar{d} (S):=\limsup_{n\rightarrow\infty}\# (S\cap \{0,1,\cdots,n-1\})$$
and $$ \underline{d} (S):=\liminf_{n\rightarrow\infty}\# (S\cap \{0,1,\cdots,n-1\}).$$ These two concepts are called {\it lower density} and {\it upper density} of $S$, respectively. If $\bar{d} (S)=\underline{d}=d,$ we call $S$ to have density of $d.$ A set $S\subseteq \mathbb{N}$ is called {\it syndetic}, if there is $N>0$ such that for any $n\geq 1,$ $$S\cap \{n,n+1,\cdots,n+N\}\neq \emptyset.  $$
Let $U,V\subseteq X$ be two nonempty open subsets and $x\in X.$ Define sets of recurrent time $$N (U,V):=\{n\geq 1|\,U\cap f^{-n} (V)\neq \emptyset\}$$ and $$N (x,U):=\{n\ge 1|\,f^n (x)\in U\}.$$ Then it is easy to check that  for any $x\in X$,
\begin{eqnarray*}
 x \text{ is almost periodic } &\Leftrightarrow& \,\forall \,\epsilon>0,\,N (x,V_\epsilon (x))\text{ is syndetic, }\\
x  \text{ is weakly almost periodic } &\Leftrightarrow&\,\forall \,\epsilon>0,\,N (x,V_\epsilon (x))\text{ has positive lower} \\
   & &\text{   density,} \,\,\,\\
x \text{ is quasi-weakly almost periodic } &\Leftrightarrow& \,\forall \,\epsilon>0,\,N (x,V_\epsilon (x))\text{ has positive upper}\\
   & &\text{ density,} \,\, \\
x    \text{ is recurrent } &\Leftrightarrow& \,\forall \,\epsilon>0,\,N (x,V_\epsilon (x))\neq \emptyset,\\
x  \text{ is non-wandering } &\Leftrightarrow& \, \forall \,\epsilon>0,\,N (V_\epsilon (x),V_\epsilon (x))\neq \emptyset.
\end{eqnarray*}
Remark that from these equivalent statements,   for recurrent  and non-wandering points  it is not required positive upper or lower density and for almost periodic points it is  required not only lower density but also some uniformly good order.  Thus,   these periodic-like recurrences  essentially reflect different `recurrent frequency'. A natural question arises:

 {\it How much difference are there between these periodic-like recurrences? }

\bigskip

One fundamental way is to consider non-emptiness of gap-sets.
Recall that $A (T)$ maybe not have totally full measure (for example, mixing subshifts of finite type and $\beta$-shifts) but $W (T)$ has totally full measure, thus the concepts of weakly and quasi-weakly almost periodic  is more general than  almost periodic from the probabilistic  perspective.  Moreover,
 many people pay attention to which system has  {\it nonempty} gap between two periodic-like level-sets (see   \cite{ZF,OS,HYZ,WHH} etc.), $QW (T)\setminus W (T)$.
  For generic diffeomorphisms and systems with Bowen's specification such as topological mixing subshifts of finite type and  topological mixing hyperbolic systems, we will prove that $QW (T)\setminus W (T)$   always contains a dense $G_\delta$ subset (see  Proposition \ref{Main-Thm-6} and Theorem \ref{Main-Thm-6-zhou-open-problem} below). This is a characterization from topological or geometric perspective so that the concept of quasi-weakly almost periodic is more general than weakly almost periodic point.
   However, note that the sets $QW (T)\setminus W (T)$ and $Rec (T)\setminus QW (T)$ have zero measure for any invariant measure.  So we need to find another way to characterize the difference.

   As said in the beginning of the paper, it is known that  topological entropy is a better and deeper tool  to study dynamical complexity than non-emptiness.
Now let us firstly consider some simple observation for general dynamical systems.

 \begin{Thm}\label{T4} {\it For any continuous   map $T:X\rightarrow X$    of a compact
metric space $X$ and a  set $Y\subseteq X$ with totally full measure, \begin{eqnarray}\label{FullMeasure-minus-periodic-FullEntropy}
Y\setminus Per (T) \text{ either is empty or carries  full topological entropy. }
\end{eqnarray}
In particular,   $ W (T)\setminus Per (T)$ either is empty or  carries  full topological entropy and by (\ref{weakperiodic-in-reccurent}) so does
  $ QW (T)\setminus Per (T)$ and $ Rec (T)\setminus Per (T)$.}

   \end{Thm}

{\bf Proof.} Suppose $Y\setminus Per (T)\neq \emptyset.$  For $Y_1=Per(T)$ and $Y_2=Y\setminus Per(T)$, by (\ref{Bowen-disjointSet-entropy})  and  Theorem \ref{T2}   $Y\setminus Per (T)$ should carry  full topological entropy.  In particular,  $ W (T)\setminus Per (T)$ carries  full topological entropy and by (\ref{weakperiodic-in-reccurent}) so do $ QW (T)\setminus Per (T)$, $ Rec (T)\setminus Per (T)$ and $ \Omega (T)\setminus Per (T)$. \qed

\bigskip

Note that $$W (T)\setminus Per (T)= (W (T)\setminus A (T))\sqcup (A (T)\setminus Per (T)),$$  where $\sqcup$ denotes the disjoint union.
 By the   discussion of $W (T)\setminus Per (T)$ in Theorem \ref{T4} and (\ref{Bowen-disjointSet-entropy}),

\begin{Thm}\label{T5}  {\it For any continuous   map $T:X\rightarrow X$    of a compact
metric space $X$,   at least one of $\,W (T)\setminus A (T)$ and  $ A (T)\setminus Per (T)$ carries  full topological entropy.}

   \end{Thm}




\section{Regularity, Quasiregularity and Irregularity} \label{section-regular----irregular}

Another important way  to observe points is from Birkhorff Ergodic theorem, called {\it quasiregular points} (for example, see \cite{DGS,Oxt}), {\it regular points} (see \cite{Oxt})   and {\it irregular points} (for example, see Chapter 8 in \cite{Bar}).

Firstly let us recall the definition of  generic point and   quasiregular point (see Chapter 4 in \cite{DGS}).
A point $x\in X$ is said to be {\it generic}  for a measure $\mu\in M(T,X)$, if for any continuous function $\phi:X\rightarrow \mathbb{R}$,   the limit $$\phi^*(x):=\lim_{n\rightarrow \infty}\frac1n\sum_{i=0}^{n-1}\phi (T^i (x)) $$ exists and equals to $\int \phi d\mu.$ Let $G_\mu$ (or $G(\mu)$) denote the set of all generic points for  $\mu,$   called generic set of $\mu$ for convenience.  By weak$^*$ topology, \begin{eqnarray}\label{eq-generic-equivalent-2015}  x \in G_\mu\,\,\Leftrightarrow\,\, \lim_{n\rightarrow \infty}\Upsilon_n (x)=\mu\Leftrightarrow M_x(T)=\{\mu\}. \end{eqnarray}  Remark that for different $\mu\neq \nu \in M(T,X),$
$ G_\mu\cap G_\nu=\emptyset.$ Recall a classical result (see Theorem 3 in \cite{Bowen1})   that for any continuous map $T:X\rightarrow X$, every ergodic measure
$\mu\in M_{erg}(T,X)$ satisfies that
\begin{eqnarray} \label{eq-Bowen-ergodicentropy-equal-generic}
h_\mu(f)=h_\mu(G_\mu).
\end{eqnarray}
By Birkhoff ergodic theorem, for any ergodic measure $\mu\in M_{erg}(T,X),$  $G_\mu$ is of $\mu$ full measure.
 Thus, by Ergodic Decomposition theorem, for any invariant measure $\mu\in M(T,X),$ $$\mu(\cup_{\nu\in M_{erg}(T,X)} G_\nu)=1.$$ This implies for any invariant measure $\mu\in M(T,X)\setminus M_{erg}(T,X),$ $G_\mu$ either is empty or satisfies that $$\mu( G_\mu)=0.$$
  A point
$x\in X$ is called   {\it quasiregular} with respect to $T$,  if  it is generic with respect to some invariant  measure.
Denote  by $QR(T)$ the set of all quasiregular points  with respect to $T$, called quasiregular set for convenience.
 By definition, 
 \begin{eqnarray}\label{eq-regular-equivalent-2015}
 x\in QR(T) \Leftrightarrow \,\,\exists\,\,\mu\in M (T,X),\,\,x\in G_\mu \Leftrightarrow M_x (T) {\text { is   a singleton}} .
 \end{eqnarray}
For convenience, for any $x\in QR(T)$,   denote  by $\mu_x$ the invariant measure for which $x$ is generic.
 Now let us introduce another concept, called $\phi$-regular,  and then we use it to give more equivalent statements of $QR(T)$.
Let    $\phi:X\rightarrow \mathbb{R}$ be  a continuous observable  function. A point  $x\in X$ is called to be {\it quasiregular for $\phi$}  with respect to $T$ (simply,  {\it $\phi$-regular}), if the limit  $\lim_{n\rightarrow \infty}\frac1n\sum_{i=0}^{n-1}\phi (T^i (x)) $ exists.  Define the $\phi$-regular set   to be the set of all $\phi$-regular points,
 that is, $$R_\phi (T):=\{x\in X| x\,\,\text{ is } \,\,\phi\text{-regular}\}.$$ Let $C^0 (X)$ denote 
  the space of all continuous functions on $X$.
 Then by definition and weak$^*$ topology,  

\begin{Thm}\label{T-8-2015-Regular-equal}
  For any continuous   map $T:X\rightarrow X$    of a compact
metric space $X$,   $$QR(T)=\bigcup_{\mu\in M(T,X)} G_\mu= \bigcap_{\phi\in C^0 (X)}R_\phi (T). $$

\end{Thm}
{\bf Proof.} The first equality is from definition. For the second equality, by definition, the relation ``$ \subseteq$" is obvious and so we only need to prove the relation ``$\supseteq$".   By contradiction, there is some $x\in X$ such that $x$ is  $\phi$-regular  point for   any continuous function $\phi:X\rightarrow\mathbb{R}$, but $x$ is not in $QR(T).$  From (\ref{eq-regular-equivalent-2015})   $\Upsilon_n (x)$ do not converge to a unique measure, then this sequence has at least two different limit points $\mu$ and $\nu$. By weak$^*$ topology and the definition of $\phi$-regular point,  for any continuous function $\phi:X\rightarrow\mathbb{R}$, the limit of   $$\frac1n\sum_{i=0}^{n-1}\phi (T^i (x))$$ equals to $\int \phi (x)d\mu$ and $\int \phi (x)d\nu$. In other words,  $\int \phi (x)d\mu=\int\phi (x)d\nu$  holds for any continuous functions.  By weak$^*$ topology this implies $\mu=\nu$, which is a contradiction. \qed

 \bigskip

Now we start to recall the concept of regular point (see \cite{Oxt}). A point $x\in QR(T)$ is called a {\it point of density}, if $\mu_x(U)>0$ for every open set $U\subseteq X$ containing $x$.  Let $QR_d(T)$ denote the set of all points of density in  $QR(T)$ and  for convenience in present paper $QR_d(T)$ is called density set.  It is easy to check that  for any $x\in QR(T),$
\begin{eqnarray}\label{eq-density-quasiregular}
 x\in QR_d(T) \Leftrightarrow x\in  S_{\mu_x}.
 \end{eqnarray} Thus \begin{eqnarray}\label{eq-density-quasiregular222222} QR_d(T) =\bigcup_{\mu\in M(T,X)} (G_\mu\cap S_\mu).\end{eqnarray}  Let $ QR_{erg}(T):= \cup_{\nu\in M_{erg}(T,X)} G_\nu.$  In \cite{Oxt} the point in $ QR_{erg}(T)$ is called transitive, but in present paper transitive point means that its orbit is dense in the whole space $X$. So for the sake of confuse, in this paper the point in $ QR_{erg}(T)$ is called {\it ergodic-transitive} and the set   $ QR_{erg}(T)$ is called ergodic-transitive set.   A point $x\in X$ is called  {\it regular}, if it belongs to the set $ R(T)=QR_{d}(T)\cap QR_{erg}(T) $ (called regular set).  Remark that
\begin{eqnarray}\label{eq-regular-in-quasiregular}
 R(T)=\bigcup_{\mu\in M_{erg}(T,X)} (G_\mu\cap S_\mu)\subseteq QR_{d}(T)\cup QR_{erg}(T) \subseteq QR(T).
 \end{eqnarray}

By  Birkhoff Ergodic theorem and Ergodic Decomposition theorem, $R(T)$  has  totally full measure (see \cite{Oxt} for a proof) and so does $QR(T)$  and every $\phi$-regular set  $R_\phi (T)$. Thus,  by Theorem \ref{thm-fullmeasure-fullentropy}  regular points (resp., quasiregular points  or $\phi$-regular points) essentially determine the dynamical  complexity of any system $T:X\rightarrow X$. That is,

\begin{Thm}\label{T-8}
  For any continuous   map $T:X\rightarrow X$    of a compact
metric space $X$,   $  R(T)$ carries full topological entropy  and  so does $QR_{d}(T)$, $ QR_{erg}(T)$,  $ QR(T)$ and every  $R_\phi (T)$  $ (\forall \phi\in C^0 (X)).$

\end{Thm}

Let $I (T)$ denote  the complementary set of quasiregular set, that is $I(T)=X\setminus QR(T).$ This set is called {\it irregular set} and its element is called {\it irregular point}
  (for example, see Chapter 8 in \cite{Bar}, also called {\it point with historic behavior} in \cite{Ruelle,Takens} and called  `non-typical' point in \cite{Barreira-Schmeling2000}).
  In other words, a point $x\in X$ is irregular if and only if there is some continuous function  $\phi:X\rightarrow \mathbb{R}$
    such that the limit  $\lim_{n\rightarrow \infty}\frac1n\sum_{i=0}^{n-1}\phi (T^i (x)) $ does not  exist. Note that every irregular point  is observed by at least one continuous function  and so let us recall a concept  called irregular point for the Birkhoff averages of  a function (for example, see \cite{Barreira-Schmeling2000,Tho2012,To2010}). Let
    $\phi:X\rightarrow \mathbb{R}$ be  a continuous    function. A point  $x\in X$ is called to be {\it irregular    for the Birkhoff averages of  $\phi$}  (simply,  {\it $\phi$-irregular}), if the limit  $\lim_{n\rightarrow \infty}\frac1n\sum_{i=0}^{n-1}\phi (T^i (x)) $ does not exist.
  Remark that $I_\phi(T)= X\setminus R_\phi (T).$
 By Theorem \ref{T-8-2015-Regular-equal},
 \begin{eqnarray}\label{eq-irregularset}
 I(T)=X\setminus QR(T)=\bigcap_{\mu\in M(T,X)} (X\setminus G_\mu)= \bigcup_{\phi\in C^0 (X)}X\setminus R_\phi (T)=\bigcup_{\phi\in C^0 (X)}I_\phi (T).
  \end{eqnarray}
Thus by  (\ref{eq-regular-equivalent-2015}) we have   \begin{eqnarray}\label{Irregular-equivalent}
x\in I (T) \Leftrightarrow M_x (T) {\text { is not a singleton}} \Leftrightarrow \Upsilon_n(x) \text{ does not converge.}
\end{eqnarray}
From Theorem \ref{T-8},  $QR(T)$ has totally full measure. Thus by (\ref{eq-irregularset})  $I(T)$ has zero measure for any invariant measure and so does every $I_\phi(T).$

However, in recent several
   years many people have focused on studying the dynamical complexity of  irregular set from different sights, for example,
    in the sense of dimension theory and topological entropy (or pressure) etc. Pesin and Pitskel   \cite{Pesin-Pitskel1984} are the first to notice the phenomenon of the irregular set carrying full topological entropy in the case of the full shift on two symbols. Barreira, Schmeling, etc. studied   the irregular set  in the setting
of subshifts of finite type and beyond, see   \cite{Barreira-Schmeling2000,Bar,TV} etc.
     Recently,  Thompson shows in   \cite{To2010,Tho2012} that  every $\phi-$irregular set  $I_\phi (T)$ either is empty or carries full topological entropy
 (or pressure) when the system satisfies (almost) specification, which is inspired from   \cite{PS} by Pfister
     and Sullivan and   \cite{TV} by Takens and Verbitskiy. For convenience of readers, let us recall a result from \cite{Tho2012}.




\begin{Thm}\label{T-7} (\cite{Tho2012}) {\it Let  $\,T$ be a      continuous map of a compact metric space $X$ with  almost  specification (slightly weaker than $g$-almost product property). Then for any continuous function $\phi:X\rightarrow \mathbb{R},$  the  $\phi-$irregular set $I_\phi(T)$ either is empty or carries full topological entropy.}

\end{Thm}

  $I (T)$ contains any $\phi-$irregular set so that it has similar results in a certain class of dynamical systems such as    \cite{Barreira-Schmeling2000,Bar,TV,To2010,Tho2012} etc. For example, by Theorem \ref{T-7} we have

\begin{Thm}\label{T-9} {\it Let  $\,T$ be a      continuous map of a compact metric space $X$ with almost  specification. Then the   irregular set $I (T)$ either is empty or carries full topological entropy.}

\end{Thm}


For the case of specification, a separable proof of this result can be found in \cite{CTS} (irregular point is called {\it divergence point} there).
 Moreover, if further assuming  that the system is not uniquely ergodic, we have

\begin{Thm}\label{T-9-1} {\it Let  $\,T$ be a      continuous map of a compact metric space $X$ with  almost  specification. If $\,T$ is not uniquely ergodic (equivalently, $M(T,X)$ is not a singleton), then the   irregular set $I (T)$  is not empty and   carries full topological entropy.}

\end{Thm}

 We need following fact from \cite{Tho2012} that

\begin{Lem}\label{Lem-IC-notempty} Let $T$ be a     continuous map of a compact metric space $X$ with    almost   specification (which is a little weaker than  $g$-almost product  property).  Let  $\psi:X\rightarrow \mathbb{R}$ be a continuous function. Then    $$\inf_{\mu\in M(T,X)}\int  \psi(x)d\mu<\sup_{\mu\in  M(T,X)}\int \psi(x)d\mu\Leftrightarrow I_\psi(T)\neq \emptyset.$$

\end{Lem}

{\bf Proof.}  For the case of `$\Leftarrow$', it is the fact (\ref{Ir-phi-nonempty-imply-differentmeasure}).  For the case of `$\Rightarrow$', see the paragraph behind of Lemma 2.1 in \cite{Tho2012}, as a corollary of Lemma 2.1 and Theorem 4.1 there on Page 5397 (for the case of specification, see Lemma 1.6 of \cite{To2010}).
 \qed

\bigskip

{\bf Proof of Theorem \ref{T-9-1}.} By assumption, there are two different invariant measures $\mu_1,\mu_2$. By weak$^*$ topology, there is a continuous function $\phi:X\rightarrow \mathbb{R}$ such that $$\int \phi d\mu_1 \neq \int \phi d\mu_2.$$ Then  $$\inf_{\mu\in M (T,X)}\int \phi (x)d\mu<\sup_{\mu\in M (T,X)}\int \phi (x)d\mu.$$
 By Lemma \ref{Lem-IC-notempty} this implies   $I_\phi (T)\neq \emptyset$ and so does $I (T)$.
Thus by Theorem \ref{T-9} we complete the proof.  \qed


\bigskip

Some classical results are known on multi-fractal analysis of Birkhoff averages, for example see   \cite{PS,TV,Ol}. More precisely, firstly  let's recall  $R_{\phi,a} (T)$ which denote the set of points whose Birkhoff average by $\phi$ equal to $a$, that is, $$R_{\phi,a} (T):=\{x\in X|\,\,\lim_{n\rightarrow \infty}\frac1n\sum_{i=0}^{n-1}\phi (T^i (x))=a\}.$$ Remark that $$ R_\phi (T)=\bigsqcup_{a\in\mathbb{R}} R_{\phi,a} (T),$$
 where $\sqcup$ denotes the disjoint union.   Let us state some basic facts as follows.
For $\phi\in C^0 (X)$ and $a\in \mathbb{R},$ by the continuity of $\phi$ and weak$^*$ topology, \begin{eqnarray}\label{Regular-equalto-a} x\in R_{\phi,a} (T) \Leftrightarrow\lim_{n\rightarrow \infty}\frac1n\sum_{i=0}^{n-1}\phi (f^i (x))=a\Leftrightarrow \,M_x (T)\subseteq \{\rho|\,\,\int \phi d\rho=a\}.\end{eqnarray}
 So by (\ref{Regular-equalto-a})
\begin{eqnarray}\label{Regular-phi-equivalent} x\in R_{\phi} (T) \Leftrightarrow \exists\,\,a\in \mathbb{R},\,\, \,M_x (T)\subseteq \{\rho|\,\,\int \phi d\rho=a\}.\end{eqnarray}
Thus one has \begin{eqnarray}\label{IrRegular-phi-equivalent} x\in I_{\phi} (T) \Leftrightarrow \exists\,\ \mu_1,\mu_2\in M_x (T)\, \text{ such that }   \int\phi d\mu_1\neq \int\phi d\mu_2.\end{eqnarray}
In particular, \begin{eqnarray}\label{Ir-phi-nonempty-imply-differentmeasure}
 I_\phi (T)\neq\emptyset\Rightarrow \inf\{\int\phi d\mu|\,\mu\in M (T,X)\}<\sup\{\int\phi d\mu|\,\mu\in M (T,X)\}.\end{eqnarray}
Recall  $ L_\phi:=[\inf\{\int\phi d\mu|\,\mu\in M (T,X)\},\,\sup\{\int\phi d\mu|\,\mu\in M (T,X)\}],$ and let $Int (L_\phi)$ denote the interior of $L_\phi.$ That is, $$Int (L_\phi)= (\inf\{\int\phi d\mu|\,\mu\in M (T,X)\},\,\sup\{\int\phi d\mu|\,\mu\in M (T,X)\}).$$
Remark that if $I_\phi (T)\neq \emptyset,$ by (\ref{IrRegular-phi-equivalent}) $Int (L_\phi)$ is a nonempty and open interval.


\bigskip

Now let us recall a result of \cite{PS} about variational principle on $R_{\phi,a}(T)$.

\begin{Prop}\label{3Main-Thm-0} (Proposition 7.1 of \cite{PS})
 Let  $\,T$ be a     continuous map of a compact metric space $X$. Suppose that $T$ is single-saturated.
 Let  $\phi:X\rightarrow \mathbb{R}$ be a continuous function.   Then for any real number $a\in  L_\phi ,$
$$ h_{top} (T,R_{\phi,a} (T))=\sup \{h_{\rho} (T)|\,\,\rho\in M (T,X)\,\,and\,\,\int \phi d\rho=a\} .$$

\end{Prop}

For uniformly hyperbolic maps and H$\ddot{\text{o}}$lder continuous functions,  Barriera and Saussol
established similar result as Proposition \ref{3Main-Thm-0}.
 The study of
multifractal analysis for arbitrary (that is, non-H$\ddot{\text{o}}$lder) continuous functions was initiated in
the symbolic dynamics setting by Olivier  \cite{Ol},  Fan and Feng \cite{Fan-Feng2000,Fan-Feng_Wu2001}.
Similar results for maps with specification can be found in \cite{TV} (for pressure, see \cite{Thompson2009}).

\bigskip


Let us consider the  entropy map
$\Psi : L_\phi \rightarrow\mathbb{R}:$ $$a \mapsto \sup \{h_{\rho} (T)|\,\,\rho\in M (T,X)\,\,and\,\,\int \phi d\rho=a\}. $$

\begin{Prop}\label{Prop-remark-3Main-Thm-0}
 Let  $\,T$ be a     continuous map of a compact metric space $X$ with positive entropy.
 Let  $\phi:X\rightarrow \mathbb{R}$ be a continuous function   with $I_\phi (T)\neq \emptyset.$    Then  $\Psi $ is a positive function on the interval $ Int(L_\phi) $.

\end{Prop}

{\bf Proof.}
Take a $\omega\in M (T,X)$ with positive entropy. If $\int \phi d\omega=a,$ then $$\sup \{h_{\rho} (T)|\,\,\rho\in M (T,X)\,\,and\,\,\int \phi d\rho=a\}\geq h_\omega(f) >0.$$ Otherwise, without loss of generality, we assume $\int \phi d\omega>a.$  By definition of $L_\phi$ and connectedness   of $M(T,X)$,  we can  take a measure $\sigma\in M(T,X)$ such that $\int \phi d\sigma<a.$ Then we can choose suitable $\xi\in (0,1)$ such that $\nu=\xi \omega+ (1-\xi)\sigma$  satisfies $\int \phi d\nu=a.$ Remark that $h_\nu (T)\geq \xi h_\omega (T) >0.$  So  $\sup \{h_{\rho} (T)|\,\,\rho\in M (T,X)\,\,and\,\,\int \phi d\rho=a\}\geq h_\nu (T)>0.$  \qed

Moreover,  we point out a result on the continuity and concave property of entropy function.

\begin{Prop}\label{Fact-concave}  Let  $\,T$ be a     continuous map of a compact metric space $X$.  Let  $\phi:X\rightarrow \mathbb{R}$ be a continuous function with $I_\phi (T)\neq \emptyset.$ Then the  entropy map
$\Psi : L_\phi \rightarrow\mathbb{R}$   is a   concave function and thus continuous on $Int (L_\phi)$.  In particular,
  $h_{top}(T)=\sup_{a\in Int (L_\phi)} \Psi(a).$

  \end{Prop}

{\bf Proof.} For any $a_1,a_2\in L_\phi,\theta\in[0,1],$  \begin{eqnarray*} & &\theta\sup\{h_\mu (f):\int\phi (x)d\mu=a_1\}+ (1-\theta)\sup\{h_\mu (f):\int\phi (x)d\mu=a_2\}\\
&=& \sup\{\theta h_{\mu_1} (f)+ (1-\theta) h_{\mu_2} (f):\int\phi (x)d\mu_i=a_i,\,i=1,2\}\\
&=& \sup\{ h_{\theta\mu_1+ (1-\theta)\mu_2} (f):\int\phi (x)d\mu_i=a_i,\,i=1,2\}\\
&\leq& \sup\{h_\mu (f):\int\phi (x)d\mu=\theta a_1+ (1-\theta) a_2\}
\end{eqnarray*}
  By classical convex analysis theory, convex function is always (locally  Lipshitz) continuous over  interior subset of the domain.

    Fix $\varepsilon>0.$  By Variational Principle, we can take one $\mu\in M (T,X)$ such that $$h_\mu (T)> h_{top} (T)-\varepsilon.$$
  Then by (\ref{Ir-phi-nonempty-imply-differentmeasure}) there is another $\nu\in M (T,X)$ such that $$\int \phi d\mu \neq \int \phi d\nu.$$ Take $\theta\in (0,1)$
  close to 1 such that $$h_\omega (T)> h_{top} (T)-\varepsilon,$$   where $\omega=\theta \mu + (1-\theta)\nu.$ Remark that the value of $\int \phi d\omega$ is
  between $\int \phi d\mu$ and $ \int \phi d\nu$ so that if $a=\int \phi d\omega$, then $a\in Int (L_\phi)$ and  
  $$\Psi(a)=\sup \{h_{\rho} (T)|\,\,\rho\in M (T,X)\,\,and\,\,\int \phi d\rho=a\}\geq h_\omega (T)>h_{top} (T)-\varepsilon.$$
   Now we complete the proof.   \qed

  \medskip

However, it is unknown  the higher smoothness of the concave function $\Psi_\xi$.

\section{Overlap of Regularity and Periodic-like Recurrence}\label{section-overlap-regular-recurrence}

Regularity and recurrence are two different ``eyes" to study asymptotic behavior.
 Inspired from above analysis,
  a natural idea is to consider the recurrence and (ir)regularity  simultaneously. Roughly speaking, under the observation of two ``eyes" , we aim to obtain more deeper results.
 In this section let us first deal with some simple relations.

\begin{Thm}\label{QW-W-in-Irregular}
For any     continuous   map $T:X\rightarrow X$    of a compact
metric space $X$,
\begin{eqnarray} R(T)\subseteq   W (T)\cap QR(T),\,\,\,\,\,Rec (T)\cap QR(T)\subseteq W (T),\,\,\,\,\,\,\, QW (T)\setminus W (T)\subseteq I (T).\nonumber
\end{eqnarray}

\end{Thm}

{\bf Proof.} We   prove the three relations respectively as follows:

(1) If $x\in R(T)$, by definition $x\in QR(T).$ Now we start to prove $x\in W(T).$ Recall  $\mu_x$ denotes the ergodic measure for which $x$ is generic.
Then by definition of regular point, it is a point of density so that for any $\varepsilon>0,$ $\mu (V_\varepsilon (x))>0.$
 By weak$^*$ topology for open sets (see Remarks (3) (iii) on Page 149 of \cite{Walter}),
 $$\underline{P}_x (V_\varepsilon (x))=\liminf_{n\rightarrow \infty}\Upsilon_n (x) (V_\varepsilon (x))\geq \mu (V_\varepsilon (x))>0.$$
 So $x\in W(T)$. 

(2) If $ x\in Rec (T)\cap QR(T),$ then there is $\mu\in M (T,X)$ such that $x\in G_\mu$ (the set of all generic points of $\mu$). So $M_x (T)=\{\mu\}$ and by (\ref{C_x=UnionSupport}) $C_x=\overline{\cup_{m\in M_x (T)}S_m}=S_\mu.$ Then by (\ref{WT}) $ x\in Rec (T)$ and $C_x=S_\nu,\,\forall\,\nu\in M_x(T)$ imply   $x\in W (T).$

(3) Let $x\in QW (T)\setminus W (T)$. Then by (\ref{weakperiodic-in-reccurent}) $ x\in Rec (T).$ If $x\in QR(T),$ by the second
  statement $x\in W (T),$ which is a contradiction to $x\in QW (T)\setminus W (T)$.
\qed

\bigskip

 Figure \ref{weak-periodic-2015-1} is to illustrate the relations between $R(T),W(T),Per(T),A(T)$ (simply, writing $R,W,Per,A$ respectively in the figure).
\setlength{\unitlength}{1mm}
  \begin{figure}[!htbp]
  \begin{center}
  \begin{picture} (180,65) (0,0)
  \put (0,0){\scalebox{1}[1]{\includegraphics[0,0][50,30]{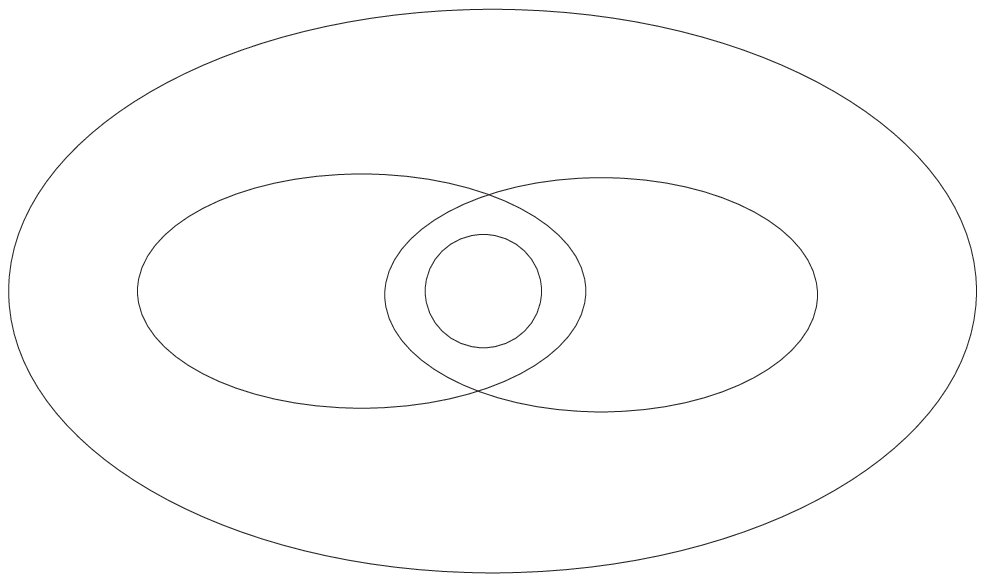}}}
   \put (80,13.5){$W$}
    \put (102,34.5){$R$}
\put (63,34.5){$A   $} \put (79.5,34.5){${\large Per} $}
  \end{picture}
  \caption{$ W(T)$.}
  \label{weak-periodic-2015-1}
  \end{center}
  \end{figure}

\bigskip

Note that obviously  one has  $Per (f)\subseteq  R(T)$.
From (\ref{FullMeasure-minus-periodic-FullEntropy}) for any dynamical system, $ R(T)\setminus Per (T)$ always carries full topological entropy, since $R(T)$ always has totally measure.
 Now we deal with $ R(T)$ and  $A (T)$.
  For full shifts on finite symbols, $ R(T)\cap  A (T)$ has  full topological entropy.

\begin{Thm}\label{T-almost-periodic}
 {\it For full shift  on $k$ symbols ($k\geq 2$),  $ R(T)\cap  A (T)\setminus Per (T)$ carries  full topological entropy and so does   $ A (T)\setminus Per (T)$. In particular, the almost periodic set $A (T)$ carries full topological entropy.}

   \end{Thm}

{\bf Proof.}
It is not difficult to prove. Denote the set of all ergodic measures of all uniquely ergodic minimal subshifts with {\it positive} entropy by  $M^*_{erg} (T,X).$  Then $$\bigcup_{\mu\in M^*_{erg} (T,X)}S_\mu =\bigcup_{\mu\in M^*_{erg} (T,X)}G_\mu\cap S_\mu \subseteq  R(T)\cap A (T) \setminus Per (T)$$ and for any $\mu \in M^*_{erg} (T,X),$  by   classical Variational Principle (Theorem 8.6 and Corollary 8.6.1 in \cite{Walter}) for $T|_{S_\mu}$,
 $$  h (T,S_\mu)=\sup_{\nu\in M_{erg} (T,S_\mu)}h_\nu(f)= h_\mu (T).$$  Recall a result from   \cite{Gri,HaKa} that for any full shift  on  finite symbols, there exist uniquely ergodic minimal subshifts with any given entropy.  This implies that
 $$h_{top} (T)=\sup_{\mu\in M^*_{erg} (T,X)}  h_\mu (T)=\sup_{\mu\in M^*_{erg} (T,X)}  h (T,S_\mu).$$
 Thus  
  by (\ref{Y1containY2}),   \begin{eqnarray}\label{AlmostPeriodic-fullentropy}
 & &h_{top} (T, R(T)\cap A (T)\setminus Per (T))\leq h_{top} (T,X)=h_{top} (T)\nonumber\\
 &=& \sup_{\mu\in M^*_{erg} (T,X)}  h (T,S_\mu)\leq h_{top} (T, R(T)\cap A (T)\setminus Per (T)).\end{eqnarray}
We complete the proof. \qed

\bigskip

Remark that even for full shifts on finite symbols, we still do not know  the entropy estimate  of   $A(T)\setminus   R(T)$.

\bigskip

 Now we consider
$  R(T)\setminus A (T) .$  Firstly let us give some  simple observation as follows.

\begin{Lem}\label{lem-non-minimalsupport-disjoint-almostperiod}
For any     continuous   map $T:X\rightarrow X$    of a compact
metric space $X$,  if there is an invariant measure $\mu$ with non-minimal support, then    $G_\mu\cap A(T)=\emptyset.$
\end{Lem}
{\bf Proof.}  By contradiction, there is $x\in G_\mu \cap A(T)$. Then $\omega_T(x)$ is a minimal compact set and $x\in\omega_T(x)$. Thus $\omega_T(x)=\overline{Orb(x)}. $  We claim that   $\mu(\omega_T(x))=1$.  More precisely, by weak$^*$ topology for closed
 sets  (see Remarks (3) (ii) on Page 149 of \cite{Walter}), $\Upsilon_{n}(x)\rightarrow \mu$ implies that
$$ 1=\limsup_{n\rightarrow \infty}\Upsilon_{n}(x)(  Orb(x) )\leq \limsup_{n\rightarrow \infty}\Upsilon_{n}(x)( \overline{Orb(x)})\leq \mu (\overline{Orb(x)})=\mu(\omega_T(x)).$$ So $\mu(\omega_T(x))=1$ implies  $S_\mu\subseteq \omega_T(x)$, contradicting that   $S_\mu$ is not minimal. \qed

\bigskip


\begin{Thm}\label{T11-section1}    For any  continuous   map $T:X\rightarrow X$    of a compact
metric space $X$, if  there is an ergodic measure $\mu$ with maximal entropy and  non-minimal  support, then    $   R(T)\setminus A (T) $  carries  full topological entropy.

   \end{Thm}

{\bf Proof.}
 Since $\mu$ is ergodic, by Birkhoff ergodic theorem  $G_\mu$ is of $\mu$ full measure.  Then
 $G_\mu\cap  S_\mu $ also has full measure, since $S_\mu$ has full measure.
 Since $\mu$  has  maximal entropy, then
  by  (\ref{FullMeasure-larger})  $h_{top}(G_\mu\cap S_\mu)\geq  h_\mu(f)=h_{top}(T).$  By   non-minimal assumption of $S_\mu$  and Lemma \ref{lem-transipoint-disjoint-AlmostPeriodc-non-minimalsystem}   $A(T)\cap G_\mu=\emptyset$.
   By ergoidcity,  $G_\mu\cap  S_\mu \subseteq   R(T)$ and thus $ G_\mu\cap  S_\mu \subseteq   R(T)\cap G_\mu \subseteq   R(T) \setminus A (T).$
    Then by (\ref{Y1containY2})
         $ h_{top}( R(T)\setminus A (T))\geq h_{top}(G_\mu\cap  S_\mu)\geq h_{top}(T). $  We complete the proof.
\qed

   \bigskip

From Theorem \ref{QW-W-in-Irregular},   $R(T) \subseteq QR(T)\cap W (T)$ and thus
$   R(T)\setminus    A (T)\subseteq     QR(T)\cap  W(T) \setminus A (T) .$
So by Theorem \ref{T11-section1} we have a following consequence.

\begin{Thm}\label{T12-section1}   For any   continuous   map $T:X\rightarrow X$    of a compact
metric space $X$, if  there is an ergodic measure $\mu$ with maximal entropy and non-minimal  support, then    $   QR(T)\cap  W (T)\setminus A (T) $ carries  full topological entropy and so does $   QR(T) \setminus A (T) $ and $      W (T)\setminus A (T) $.
\end{Thm}



   Remark  that Theorem \ref{T11-section1} and Theorem \ref{T12-section1} are applicative to   all the examples in Section \ref{section-examples}, since each example    is non-minimal and  the  unique maximal entropy measure  is ergodic and has full support.


\section{Some   useful facts and lemmas}\label{section-usefulfacts-lemmas}

Firstly let us recall two basic properties of invariant measures with full support from \cite{DGS}.  Let $M_{full}(T,X)$ denote the set of all invariant measures with full support.

\begin{Lem}\label{lem-DGS-fullsupport-empty-2015} (Proposition 21.11 in \cite{DGS}) Let $T:X\rightarrow X$ be a continuous map on a  compact metric space $X.$
  Then the set $M_{full}(T,X)$  is either empty or a dense $G_\delta$ subset in $M(T,X).$


\end{Lem}

\begin{Lem}\label{lem-DGS-full-support-2015-22222} (Proposition 21.12 in \cite{DGS}) Let $T:X\rightarrow X$ be a continuous map on a  compact metric space $X.$
 If  the periodic points are dense in $X$, then
    the set $M_{full}(T,X)$ is  a dense $G_\delta$ subset in $M(T,X).$


\end{Lem}

A direct corollary of Lemma \ref{lem-DGS-fullsupport-empty-2015} is following.

\begin{Lem}\label{lem-DGS-fullsupport-2015}   Let $T:X\rightarrow X$ be a continuous map on a  compact metric space $X.$
 If  there is an invariant measure with full support,  then
    the set $M_{full}(T,X)$  is   a dense $G_\delta$ subset in $M(T,X).$


\end{Lem}

Moreover, we discuss the  relation between periodic points, periodic measures  and measures with full support.

\begin{Prop}\label{Prop-supportimplies-dense-periodic}  Let $T:X\rightarrow X$ be a continuous map on a  compact metric space $X.$
 If  there is an invariant measure $\mu$  with full support and a sequence of periodic measures $\mu_n$   such that $\mu_n$ converges to $\mu$ in weak$^*$ topology.
 Then the periodic points are dense in $X$.


\end{Prop}

{\bf Proof.} By weak$^*$ topology for closed
 sets  (see Remarks (3) (ii) on Page 149 of \cite{Walter}),
$$ 1=\limsup_{n\rightarrow \infty}\mu_n( \overline{Per(T)})\leq \mu (\overline{Per(T)}).$$
 So  $S_\mu \subseteq \overline{Per(T)}$ and then $X = \overline{Per(T)}$, since $\mu$ has full support. \qed

\bigskip

By Proposition \ref{Prop-supportimplies-dense-periodic} and Lemma \ref{lem-DGS-full-support-2015-22222}, we have a following consequence.

\begin{Prop}\label{Prop-support-same-as-dense-periodic}  Let $T:X\rightarrow X$ be a continuous map on a  compact metric space $X$. Suppose that the periodic measures are dense in the space of invariant measures. Then    $$\overline{ Per(T)}=X \Leftrightarrow \exists \mu \in M(T,X),\,\,S_\mu=X.$$

\end{Prop}

\medskip

For non-ergodic measures,  we have a following result.

\begin{Lem}\label{lem-non-ergodic-2015}   Let $T:X\rightarrow X$ be a not uniquely ergodic continuous map on a  compact metric space $X.$
     Then  the set of non-ergodic measures  $M(T,X)\setminus M_{erg}(T,X)$  is   a dense   subset in $M(T,X).$


\end{Lem}

 {\bf Proof.} One can use Ergodic Decomposition theorem to prove, here we give a constructed  proof. In fact, we only need to prove for any given $\mu\in M_{erg}(T,X),$
     $\mu $ can be approximated by non-ergodic measures.  More precisely, by the assumption of not uniquely ergodic, take another ergodic measure $\nu\neq \mu.$ Let $\nu_n=\frac1n \nu +(1-\frac1n)\mu.$ Then $\nu_n$ converges to $\mu$ in weak$^*$ topology. Remark that obviously  $\nu_n$ are all not ergodic. \qed

     \bigskip

Now we start to consider the number of periodic points  in open sets.

\begin{Prop}\label{prop-number-periodic-in-openset}
Let $T:X\rightarrow X$ be a continuous map on a  compact metric space $X.$
 If the periodic points are dense in $X$ (i.e., $\overline{Per(f)}=X$) and the periodic measures are dense in the space of invariant measures (i.e., $\overline{ M_{p} (T,X) }=M (T,X)),$ then for any nonempty open set $U\subseteq X,$ there are  infinite periodic points in $U.$ In particular, every nonempty set $U$ is not finite.

\end{Prop}

{\bf Proof.}   By contradiction, there is some nonempty open set $U\subseteq X$ such that $U$   contains at most  finite periodic points. This implies that  $$\Delta:=\{x\in Per(T)|\,Orb(x)\cap U\neq \emptyset\}$$ is an  invariant set with at most   finite  elements and thus closed.
 By assumption of density of periodic points, the open set $U$ contains some  periodic point. Thus $\Delta$ is a nonempty  finite,  closed invariant set.
 Take $\mu:=\frac 1{\#\Delta} \sum_{x\in \Delta}\delta_x.$  It is a finite convex sum of all periodic measures supported on  $\Delta$ and  thus $\mu(U)>0$.


  Since in present paper $X$ is assumed infinite, then by assumption,  the number of periodic points are infinite.
 Then we can take a periodic measure $\nu$ whose orbit is contained in $X\setminus U.$  So $\nu(\Delta)=0$.  Let $\omega=\frac12 \mu+\frac12 \nu$. Then by assumption, there is a sequence of periodic measures $\omega_n$ such that $\omega_n$ converges to $\omega$ in weak$^*$ topology. By weak$^*$ topology for open sets  (see Remarks (3) (iii) on Page 149 of \cite{Walter}), $$\liminf_{n\rightarrow \infty} \omega_n(U)\geq \omega(U)=\frac12\mu(U)>0.$$ Then we can take a large $N_1$ such that for any $n\geq N_1,$  $\omega_n(U)>0.$
 So the periodic orbit $S_{\omega_n}$  is contained in $\Delta$ and thus $\omega_n(\Delta)=1.$   On the other hand, Using weak$^*$ topology for the  open set $X\setminus \Delta$,
 $$\liminf_{n\rightarrow \infty} \omega_n(X\setminus \Delta)\geq \omega(X\setminus \Delta)\geq \frac12\nu(X\setminus \Delta)=\frac12>0.$$ Then we can take a large integer $N>N_1$ such that   $\omega_N(X\setminus \Delta)>0.$  This contradicts $\omega_N(\Delta)=1.$  \qed

\bigskip

By weak$^*$ topology and the continuity of $\phi\in C^0 (X)$, it is easy to see that if $\mathcal{R}$ is a dense subset of $M (T,X)$, then $$\{\int \phi (x) d \nu |\,\,\nu\in \mathcal{R}\}$$ is dense in $L_\phi.$ This implies that

\begin{Lem}\label{lem-phi-different-000000} Let $T:X\rightarrow X$ be a continuous map on a  compact metric space $X.$   Let $\phi\in C^0 (X)$ with $I_\phi (T)\neq \emptyset$. Then for any $\nu\in M (T,X)$,  there is $\varrho\in \mathcal{R}$ such that $\int\phi d\varrho\neq \int\phi d\nu.$ Moreover, for any $n\in \mathbb{Z}^+$ and $a\in Int (L_\phi),$  there is  $\mu_1, \mu_2, \cdots, \mu_n\in \mathcal{R}$ and $\nu_1, \nu_2, \cdots, \nu_n\in\mathcal{R}$ such that $$\int \phi (x) d\mu_n<\cdots<\int \phi (x) d\mu_1<a<\int \phi (x) d\nu_1   <\cdots<\int \phi (x) d\nu_n.$$





\end{Lem}

 In particular, by Lemma \ref{lem-phi-different-000000}   we have a following consequence.

 \begin{Lem}\label{lem-phi-different} Let $T:X\rightarrow X$ be a continuous map on a  compact metric space $X.$   Let $\phi\in C^0 (X)$ with $I_\phi (T)\neq \emptyset$. Then \\
 (1) if $M_p(T,X)$  is dense in $M(T,X),$ then for any $\nu\in M (T,X)$,  there is $\varrho\in M_p(T,X)$ such that $\int\phi d\varrho\neq \int\phi d\nu.$  Moreover, for any $n\in \mathbb{Z}^+$ and $a\in Int (L_\phi),$  there is  $\mu_1, \mu_2, \cdots, \mu_n$ $\in M_p(T,X)$ and $\nu_1, \nu_2, \cdots, \nu_n\in M_p(T,X)$ such that $ \int \phi (x) d\mu_n<\cdots<\int \phi (x) d\mu_1<a<\int \phi (x) d\nu_1 < \cdots<\int \phi (x) d\nu_n.$  Moreover, $$ Int (L_\phi)\subseteq \{ \int \phi d\mu\,|\,\mu\in M(T,X),\,\, S_\mu \text{ is composed of  at most two periodic orbits}\}.$$

 (2) if there is some invariant measure    with full support, then for any $\nu\in M (T,X)$,  there is $\varrho\in M_{full}(T,X)$ such that $\int\phi d\varrho\neq \int\phi d\nu.$  Moreover, for any $n\in \mathbb{Z}^+$ and $a\in Int (L_\phi),$  there is  $\mu_1, \mu_2, \cdots, \mu_n\in M_{full}(T,X)$ and $\nu_1, \nu_2, \cdots, \nu_n\in M_{full}(T,X)$ such that $\int \phi (x) d\mu_n<\cdots<\int \phi (x) d\mu_1<a<\int \phi (x) d\nu_1 <\cdots<\int \phi (x) d\nu_n.$
 Moreover, $$ Int (L_\phi)\subseteq \{ \int \phi d\mu\,|\,\mu\in M(T,X),\,\, S_\mu=X\}.$$

(3) if $M_{erg}(T,X)$  is dense in $M(T,X),$  then for any $\nu\in M (T,X)$,  there is $\varrho\in M_{erg}(T,X)$ such that $\int\phi d\varrho\neq \int\phi d\nu.$  Moreover, for any $n\in \mathbb{Z}^+$ and $a\in Int (L_\phi),$  there is  $\mu_1, \mu_2, \cdots, \mu_n $ $\in M_{erg}(T,X)$ and $\nu_1, \nu_2, \cdots, \nu_n\in M_{erg}(T,X)$ such that  $\int \phi (x) d\mu_n<\cdots<\int \phi (x) d\mu_1<a<\int \phi (x) d\nu_1   <\cdots<\int \phi (x) d\nu_n.$ \\

 (4) if there is some invariant measure    with full support and $M_{erg}(T,X)$  is dense in $M(T,X),$  then for any $\nu\in M (T,X)$,  there is $\varrho\in M_{erg}(T,X) \cap M_{full}(T,X)$ such that $\int\phi d\varrho\neq \int\phi d\nu.$  Moreover, for any $n\in \mathbb{Z}^+$ and $a\in Int (L_\phi),$  there is  $\mu_1, \mu_2, \cdots, \mu_n\in M_{erg}(T,X)\cap M_{full}(T,X)$ and $\nu_1, \nu_2, \cdots, \nu_n\in M_{erg}(T,X)\cap M_{full}(T,X)$ such that  $\int \phi (x) d\mu_n$ $<\cdots<\int \phi (x) d\mu_1<a<\int \phi (x) d\nu_1   <\cdots<\int \phi (x) d\nu_n.$ \\

 (5)  if there is some invariant measure    with full support and $M_{p}(T,X)$  is dense in $M(T,X),$  then for any $\nu\in M (T,X)$,  there is $\varrho\in M_{erg}(T,X) \cap M_{full}(T,X)$ such that $\int\phi d\varrho\neq \int\phi d\nu.$  Moreover, for any $n\in \mathbb{Z}^+$ and $a\in Int (L_\phi),$  there is  $\mu_1, \mu_2, \cdots, \mu_n\in M_{erg}(T,X)\cap M_{full}(T,X)$ and $\nu_1, \nu_2, \cdots, \nu_n\in M_{erg}(T,X)\cap M_{full}(T,X)$ such that  $\int \phi (x) d\mu_n$ $<\cdots<\int \phi (x) d\mu_1<a<\int \phi (x) d\nu_1   <\cdots<\int \phi (x) d\nu_n.$ \\

(6)    for any $\nu\in M (T,X)$,  there is $\varrho\in M(T,X)\setminus M_{erg}(T,X)$ such that $\int\phi d\varrho\neq \int\phi d\nu.$  Moreover, for any $n\in \mathbb{Z}^+$ and $a\in Int (L_\phi),$  there is  $\mu_1, \mu_2, \cdots, \mu_n\in M(T,X)\setminus M_{erg}(T,X)$ and $\nu_1, \nu_2, \cdots, \nu_n\in M(T,X)\setminus M_{erg}(T,X)$ such that  $\int \phi (x) d\mu_n<\cdots<\int \phi (x) d\mu_1<$ $a<\int \phi (x) d\nu_1   <\cdots<\int \phi (x) d\nu_n.$
Moreover, $$ Int (L_\phi)\subseteq \{ \int \phi d\mu\,|\,\mu\in M(T,X)\setminus M_{erg}(T,X)\}.$$

\end{Lem}

 {\bf Proof.}  (1)     Letting $\mathcal{R}=M_p(T,X)$, by   Lemma \ref{lem-phi-different-000000}  we only need to prove $$ Int (L_\phi)\subseteq \{ \int \phi d\mu\,|\,\mu\in M(T,X),\,\, S_\mu \text{ is composed of  at most two periodic orbits}\}.$$  Take  two periodic measures $\mu_1,\nu_1$ such that $\int \phi (x) d\mu_1<a<\int \phi (x) d\nu_1$. Then we can choose some suitable
 $\theta\in (0,1)$ such that $\int \phi d\omega=a,$ where $\omega=\theta \mu_1+(1-\theta)\nu_1.$ Remark that $S_\omega$ is composed of at most two periodic orbits.

 (2) By  Lemma \ref{lem-DGS-fullsupport-2015} $M_{full}(T,X)$ is dense in $M(T,X)$.  So by   Lemma \ref{lem-phi-different-000000}  we only need to prove
 $$ Int (L_\phi)\subseteq \{ \int \phi d\mu\,|\,\mu\in M(T,X),\,\, S_\mu=X\}.$$ Take  two invariant  measures $\mu_1,\nu_1$ with full support such that $\int \phi (x) d\mu_1<a<\int \phi (x) d\nu_1$. Then we can choose some suitable
 $\theta\in (0,1)$ such that $\int \phi d\omega=a,$ where $\omega=\theta \mu_1+(1-\theta)\nu_1.$ Remark that $S_\omega$ also has full support.

 (3)  It is obvious just letting  $\mathcal{R}=M_{full}(T,X)$ in   Lemma \ref{lem-phi-different-000000}.

 (4)   It is well-known that $M_{erg}(T,X)$ is a non-empty $G_\delta$ subset of $M(T,X)$ (see Proposition 5.7 in \cite{DGS}). Thus,
 by assumption and Lemma \ref{lem-DGS-fullsupport-2015}, $M_{erg}(T,X)$ and  $M_{full}(T,X)$ both are dense  $G_\delta$ subsets of $M(T,X)$. So
 $M_{erg}(T,X)\cap M_{full}(T,X)$ is also   dense  $G_\delta$ in $M(T,X)$. Then  Lemma \ref{lem-phi-different-000000} implies that (4) is true.

 (5) Note that $M_{p}(T,X)\subseteq M_{erg}(T,X).$ So by assumption,  $M_{erg}(T,X)$  is dense in $M(T,X) $ and thus (4) implies (5).

 (6)
From (\ref{Ir-phi-nonempty-imply-differentmeasure})  $I_\phi (T)\neq \emptyset $   implies that the system is naturally not uniquely ergodic.
By Lemma \ref{lem-non-ergodic-2015}, $M(T,X) \setminus M_{erg}(T,X)$ is dense in $M(T,X)$. Let $\mathcal{R}=M(T,X) \setminus M_{erg}(T,X)$ and then by
Lemma \ref{lem-DGS-fullsupport-2015},  we only need to show $$ Int (L_\phi)\subseteq \{ \int \phi d\mu\,|\,\mu\in M(T,X)\setminus M_{erg}(T,X)\}.$$
Take  two non-ergodic  measures $\mu_1,\nu_1$ such that $\int \phi (x) d\mu_1<a<\int \phi (x) d\nu_1$. Then we can choose some suitable
 $\theta\in (0,1)$ such that $\int \phi d\omega=a,$ where $\omega=\theta \mu_1+(1-\theta)\nu_1.$ Remark that $\omega$ is not ergodic, since it is known that every ergodic measure is extremal point  in $M(T,X)$ (for example, see Proposition 5.6 in \cite{DGS}).
 \qed



\section{Proofs of main results}\label{proofs}

\subsection{ Proof  of    Theorem \ref{More-general-systems-2} (2)}\label{proof-main}

In this section we divide Theorem \ref{More-general-systems-2} (2) into several propositions to prove. Here we will state every proposition  for  possible more general systems (for example, it is possible to apply partial  results in some partially hyperbolic systems, see Section \ref{section-partial-hyperbolic}).

\begin{Prop}\label{Main-Thm-1}  
 Let  $\,T$ be a     continuous map of a compact metric space $X$. Suppose that $T$ is saturated and there is some invariant measure $\mu$ with full support (i.e., $S_\mu=X $).
 Then for any $\phi\in C^0 (X)$, either $I_\phi (T)=\emptyset$ or
$$h_{top} (T, (W (T)\setminus QR(T))\cap I_\phi (T))=h_{top} (T, I_\phi (T))=h_{top} (T).$$

\end{Prop}

{\bf Proof.} Suppose $I_\phi (T)\neq\emptyset$ and
fix $\varepsilon>0.$ By    classical Variational Principle (Theorem 8.6   in \cite{Walter}), we can take  $\mu\in M (T,X)$ such that $$h_\mu (T)>h_{top} (T)-\varepsilon.$$ By Lemma \ref{lem-phi-different} (2), we can take $\nu\in M (T,X)$  such that $S_\nu=X$ and $\int\phi d\nu\neq\int\phi d\mu$. Then we can choose two different numbers $0<\theta_1<\theta_2<1$ close to 1 enough such that for $\omega_i=\theta_i\mu + (1-\theta_i)\nu,\,i=1,2,$ one has
\begin{eqnarray}\label{entropy-1} \,h_{\omega_i} (T)=\theta_ih_\mu (T) + (1-\theta_i)h_\nu (T)\geq \theta_ih_\mu (T)>h_{top} (T)-\varepsilon,\,i=1,2 .\,\,\,\,\,\,\,\,\,\,\,\,\end{eqnarray}
Remark that $\theta_1\neq\theta_2$ and $\int\phi d\mu\neq \int\phi d\nu$ imply
\begin{eqnarray}\label{phi-1}\int\phi d\omega_1\neq\int\phi d\omega_2;\, \end{eqnarray} and $S_\nu=X$ implies
\begin{eqnarray} \label{support1} S_{\omega_i}=S_\mu\cup S_\nu=X,\,i=1,2. \,\,\,\,\,\,\,\end{eqnarray}

 Let $$K=\{\tau \omega_1 + (1-\tau)\omega_2|\,\tau\in [0,1]\}.$$ Then by (\ref{support1}) and (\ref{entropy-1})
  for any $m=\tau \omega_1 + (1-\tau)\omega_2\in K,$
$$S_m=X,\,\,\,h_m (T)\geq \min\{h_{\omega_1} (T),h_{\omega_2} (T)\}>h_{top} (T)-\varepsilon.$$
  Since $T$ is saturated, then
   $$h_{top} (T,G_K)=\inf\{h_m (T)\,|\,m\in K\}\geq h_{top} (T)-\varepsilon,$$ where $G_K=\{x\in X|\,M_x (T)=K\}.$

To complete the proof, we only need to prove $G_K\subseteq I_\phi (T)\cap W (T).$
On one hand, for any $x\in G_K$, notice that $\omega_1,\omega_2\in K= M_x (T)$ and thus by (\ref{IrRegular-phi-equivalent}), (\ref{phi-1}) implies  $x\in I_\phi (T).$ On the other hand, by (\ref{C_x=UnionSupport}) the equality  $$C_x=\overline{\cup_{\kappa\in M_x (T)}S_\kappa}=\overline{\cup_{\kappa\in K}S_\kappa}=X=S_m$$ hold for all $m\in K=M_x (T).$ By (\ref{C_xinOmega}) $x\in X=C_x\subseteq \omega_T (x)$ and thus $x\in Rec (T).$  So from (\ref{WT}) one has $x\in W (T).$
\qed

\bigskip

\begin{Prop}\label{Main-Thm-2}
Let  $\,T$ be a      continuous map of a compact metric space $X$. 
Suppose that  $T$ is saturated,  has entropy-dense property and
there is some invariant measure with full support.
Then for any $\phi\in C^0 (X)$, either $I_\phi (T)=\emptyset$ or
$$h_{top} (T, I_\phi (T)\cap V(T)\setminus W(T))=h_{top} (T, I_\phi (T))=h_{top} (T).$$  

\end{Prop}


{\bf Proof.} Suppose $I_\phi (T)\neq\emptyset$ and
fix $\varepsilon>0.$   
By    classical Variational Principle (Theorem 8.6   in \cite{Walter}), we can take  $\mu_0\in M (T,X)$ such that $$h_{\mu_0} (T)>h_{top} (T)-\varepsilon.$$ By entropy-dense property  we can choose $\mu\in  M_{erg} (T,X)$ (close to $\mu_0$) such that $S_\mu\subsetneqq X$ and  $$h_\mu (T)>h_{top} (T)-\varepsilon.$$
 By Lemma \ref{lem-phi-different} (2),  $I_\phi (T)\neq \emptyset$ implies that we can take $\nu\in M (T,X)$   such that $S_\nu=X$ and $\int\phi d\nu\neq\int\phi d\mu$. Then we can take $\theta\in (0,1)$ close to 1 such that
 $h_\omega (T)\geq \theta h_\mu (T) >h_{top} (T)-\varepsilon$  where  $\omega=\theta\mu+ (1-\theta)\nu.$ Remark that $S_\omega=S_\mu\cup S_\nu=X$ and \begin{eqnarray}\label{different-phi2}\int\phi d\omega\neq\int\phi d\mu.\end{eqnarray}

 Let $$K=\{\tau \omega + (1-\tau)\mu|\,\tau\in [0,1]\}.$$ Then
  for any $m=\tau \omega + (1-\tau)\mu\in K\setminus\{\mu\},$
$S_m=X$ and for any $m=\tau \omega + (1-\tau)\mu\in K$ $$h_m (T)\geq \min\{h_{\omega} (T),h_{\mu} (T)\}>h_{top} (T)-\varepsilon.$$
  Since $T$ is saturated, then
   $$h_{top} (T,G_K)=\inf\{h_m (T)\,|\,m\in K\}>h_{top} (T)-\varepsilon,$$ where $G_K=\{x\in X|\,M_x (T)=K\}.$


To complete the proof, we  only need to prove $G_K\subseteq   I_\phi (T)\cap\{x\in  QW (T)\setminus W (T)|\exists\, \omega\in M_x (T) \,s.t. \,S_\omega=C_x\}.$  In fact, fix
 $x\in G_K$.  on one hand, notice that  $\omega,\mu\in K= M_x (T)$   and thus by (\ref{IrRegular-phi-equivalent}), (\ref{different-phi2}) implies $x\in I_\phi (T).$ On the other hand, by (\ref{C_x=UnionSupport}) one has the following equality $$C_x=\overline{\cup_{m\in M_x (T)}S_m}=\overline{\cup_{m\in K}S_m}=X.$$  By (\ref{C_x=X}) $x\in QW (T).$  Thus by (\ref{weakperiodic-in-reccurent}) $x\in  Rec (T).$ Notice that $\mu\in M_x (T)$ and  $C_x=X\neq S_\mu$ so that from (\ref{WT})   $x\in X\setminus W (T)$. Recall that $\omega\in K=M_x (T)$ and $S_\omega=X.$ So $x\in I_\phi (T)\cap\{x\in  QW (T)\setminus W (T)|\exists\, \omega\in M_x (T) \,s.t. \,S_\omega=C_x\}.$ 
\qed

\bigskip

\begin{Prop}\label{Main-Thm-3} 
 Let  $\,T$ be a      continuous map of a compact metric space $X$.  
  Suppose that  $T$ is saturated and  has entropy-dense property,
 the periodic points are dense in $X$ (i.e., $\overline{Per(f)}=X$) and
 the periodic measures are dense in the space of invariant measures (i.e., $\overline{ M_{p} (T,X) }=M (T,X)).$
Then for any $\phi\in C^0 (X)$, either $I_\phi (T)=\emptyset$ or
$$h_{top} (T,I_\phi (T)\cap    QW (T)\setminus V (T))=h_{top} (T, I_\phi (T))
=h_{top} (T).
$$

\end{Prop}


{\bf Proof.} Suppose $I_\phi (T)\neq\emptyset$ and
fix $\varepsilon>0.$ 
By    classical Variational Principle (Theorem 8.6   in \cite{Walter}), we can take  $\mu_0\in M (T,X)$ such that $$h_{\mu_0} (T)>h_{top} (T)-\varepsilon.$$
 By entropy-dense property,  we can choose $\mu\in  M_{erg} (T,X)$ (close to $\mu_0$) such that $S_\mu\subsetneqq X$ and  $$h_\mu (T)>h_{top} (T)-\varepsilon.$$

 Since $X$ and $M(T,X)$ are compact metric spaces, by assumption  we can take   a countable dense subset $P_1\subseteq Per(T)$ and a countable dense subset $M_1\subseteq M_p(T,X)$. Then   $  \cup_{\mu\in M_1} S_\mu\cup P_1$ is still a countable dense subset  of $Per(T)$, denoted by $\{x_i\}_{i=1}^\infty $.
  Moreover,  $\bigcup_{x\in P_1}m_x \cup M_1$  is countable and dense in $M (T,X),$ where $m_x$ denote the $T$-invariant measure supported on the periodic orbit of $x$.   So $\overline{\{m_i\}_{i=1}^\infty}=M (T,X)$,  where $m_i $  denotes the $T$-invariant measure supported on the periodic orbit of $x_i$.



  Take a strictly  increasing sequence of $\{\theta_i|\,\,\theta_i\in (0,1)\}_{i=1}^\infty$ such that  $$\lim_{i\rightarrow +\infty}\theta_i=1$$ and  $$h_{\nu_i} (T)\geq \theta_i h_\mu (T) >h_{top} (T)-\varepsilon$$ where $\nu_i=\theta_i\mu+ (1-\theta_i)m_i,\,i=1,2,3,\cdots.$ Remark that for any $i,$ $S_{\nu_i}=S_\mu\cup S_{m_i}.$

  By Lemma \ref{lem-phi-different-000000}, if let $\mathcal{R}=\{m_i\}_{i=1}^\infty,$  then $I_\phi (T)\neq \emptyset$ implies that we can take one periodic measure $m_{i_0}$  such that $\int\phi dm_{i_0}\neq\int\phi d\mu$.    Without loss of generality, we can assume $m_1=m_{i_0}$. Then \begin{eqnarray}\label{different-phi3}\int\phi d\nu_1\neq\int\phi d\mu.\end{eqnarray}

Now we consider $$K=\bigcup_{i\geq 1} \{\tau \nu_i+ (1-\tau)\nu_{i+1}|\,\tau\in[0,1]\}\,\,\,\bigcup\,\,\,\{\mu\}.$$
Remark that $K$ is a nonempty connected compact subset of $M (T,X)$ because $\nu_i\rightarrow \mu$ in weak$^*$ topology.
 Since $T$ is saturated, then  $$h_{top} (T,G_K)=\inf\{h_m (T)\,|\,m\in K\}=\min\{\inf_{i\geq 1}\{h_{\nu_i} (T) \},h_\mu (T)\}\geq h_{top} (T)-\varepsilon,$$ where $G_K=\{x\in X|\,M_x (T)=K\}.$

To complete the proof, we only need to prove $G_K\subseteq I_\phi (T)\cap \{x\in  QW (T)\setminus W (T)|\forall\, m\in M_x (T) \,s.t. \,S_m\neq C_x\}. $ 
Fix $x\in G_K$. Recall that $\nu_1,\mu\in K= M_x (T)$  and   thus by (\ref{IrRegular-phi-equivalent}), (\ref{different-phi3}) implies
   $x\in I_\phi (T).$ Clearly by (\ref{C_x=UnionSupport}) $$C_x=\overline{\cup_{m\in M_x (T)}S_m}=\overline{\cup_{m\in K}S_m}\supseteq\overline{\cup_{i\geq 1}S_{\nu_i} } =\overline{\cup_{i\geq 1} (S_{m_i}\cup S_\mu)}\supseteq\overline{\cup_{i\geq 1}\{x_i\}}=X.$$  So $C_x=X.$ By (\ref{C_x=X}) $x\in QW (T).$ By (\ref{weakperiodic-in-reccurent}) $x\in  Rec (T).$ Notice that  $C_x=X\neq S_\mu$ and $\mu\in M_x (T)$ so that from (\ref{WT}) $x\in X\setminus W (T)$.  For any $ m\in M_x (T)=K,$ by definition of $K$ there is some $i$ such that $S_m\subseteq S_\mu\cup S_{m_i}\cup S_{m_{i+1}}.$ Note that $X\setminus S_\mu$ is a nonempty open set so that by Proposition \ref{prop-number-periodic-in-openset}   $X\setminus S_\mu$ is not a finite set. But $S_{m_i}\cup S_{m_{i+1}}$ is just composed of two periodic orbits and thus is a finite set. Hence $ S_m \subsetneqq X=C_x$. I.e., $S_m\neq C_x.$\qed

\bigskip

\begin{Prop}\label{Main-Thm-4}  
Let  $\,T$ be a      continuous map of a compact metric space $X$. 
 Suppose that  $T$ is saturated and  has entropy-dense property,
 the periodic points are dense in $X$ (i.e., $\overline{Per(f)}=X$) and
 the periodic measures are dense in the space of invariant measures (i.e., $\overline{ M_{p} (T,X) }=M (T,X)).$  Then for any $\phi\in C^0 (X)$, either $I_\phi (T)=\emptyset$ or
$$h_{top} (T,I_\phi (T)\cap I (T)\setminus QW (T))=h_{top} (T, I_\phi (T))=h_{top} (T).$$

\end{Prop}


{\bf Proof.} Suppose $I_\phi (T)\neq\emptyset$ and
fix $\varepsilon>0.$
By    classical Variational Principle (Theorem 8.6   in \cite{Walter}), we can take  $\mu_0\in M (T,X)$ such that $$h_{\mu_0} (T)>h_{top} (T)-\varepsilon.$$  By entropy-dense property,   we can choose $\mu\in  M_{erg} (T,X)$ (close to $\mu_0$) such that $S_\mu\subsetneqq X$
and  $$h_\mu (T)>h_{top} (T)-\varepsilon.$$ Since $ X\setminus S_\mu$
 is nonempty and open, by density of periodic points 
 $$B:=\{\nu\in M_p (T,X)|\,\,S_{\nu}\setminus S_\mu\neq \emptyset\}\neq \emptyset.$$

Now we will construct an invariant measure $\kappa$ such that the set  $S_{\kappa}\setminus S_\mu$ is composed of one periodic orbit and $$\int\phi d\kappa \neq\int\phi d\mu.$$
More precisely,  if there is a periodic measure $\nu\in B$ such that   $\int\phi d\nu\neq\int\phi d\mu,$ then take $\kappa=\nu$ and remark that $S_{\kappa}\setminus S_\mu=S_\nu.$ Otherwise, for any $\nu\in B,$   $\int\phi d\nu=\int\phi d\mu.$ Take such a measure $\nu.$
By Lemma \ref{lem-phi-different} (1),  $\overline{ M_{p} (T,X) }=M (T,X)$ and $I_\phi (T)\neq \emptyset$ imply $$Y:=\{\tau|\,\,\int\phi d\tau\neq\int\phi d\mu,\,\tau\in M_p (T,X)\}\neq \emptyset.$$  Then we can take $\nu'\in Y$ such that $\int\phi d\nu'\neq\int\phi d\mu.$   Remark that in this case $Y\cap B=\emptyset$ so that $S_{\nu'} \setminus S_\mu=\emptyset.$ Then  $S_{\nu'}\subseteq S_\mu.$
So if we take $\kappa=\frac12 (\nu+\nu'),$ then   $\int\phi d\kappa\neq\int\phi d\mu.$ Note that $S_\kappa=S_\nu\cup S_{\nu'} $ and $S_{\kappa}\setminus S_\mu=S_\nu. $





\bigskip

 Take $\theta\in (0,1)$ close to 1 such that  $\omega=\theta \mu+ (1-\theta)\kappa$ satisfies $h_\omega (T)\geq \theta h_\mu (T)> h_{top} (T)-\varepsilon.$
 Then $\omega$ also satisfies that $\int\phi d\omega\neq\int\phi d\mu.$ Remark that  $S_\omega=S_\mu\cup S_\nu\cup S_{\nu'}= S_\mu\sqcup S_\nu.$ Let $K=\{\tau\mu+ (1-\tau)\omega|\,\tau\in[0,1]\}.$
 Since $T$ is saturated, then  
   $$h_{top} (T,G_K)=\inf\{h_m (T)\,|\,m\in K\}\geq \min\{h_\mu (T),h_\omega (T)\}>h_{top} (T)-\varepsilon,$$ where $G_K=\{x\in X|\,M_x (T)=K\}.$

To complete the proof, we only need to prove $$G_K\subseteq I_\phi (T)\setminus QW (T).$$
 For any $x\in G_K$, recall that  $\omega,\mu\in K= M_x (T)$ and $\int\phi d\omega\neq\int\phi d\mu.$ Thus by (\ref{IrRegular-phi-equivalent}) $x\in I_\phi (T).$ If $x\in QW (T)$, then  by (\ref{weakperiodic-in-reccurent}) $x\in Rec (T)$ so that  by (\ref{QW})   $x\in C_x=\omega_T (x).$ Then by (\ref{C_x=UnionSupport}) $$x\in \omega_T (x)=C_x=\overline{\cup_{m\in M_x (T)}S_m}=\overline{\cup_{m\in K}S_m}=S_\mu\sqcup S_\nu,$$ which implies    $x \in S_\mu$ or $S_\nu.$     So by invariance of $S_\mu$ and $S_\nu,$ one has $Orb (x)\subseteq S_\mu$ or $Orb (x)\subseteq S_\nu$ and thus by compactness of $S_\mu$ and $S_\nu,$ we have $\omega_T (x)\subseteq S_\mu\subsetneqq S_\mu\sqcup S_\nu=C_x=\omega_T (x)$ or $\omega_T (x)\subseteq S_\nu\subsetneqq S_\mu\sqcup S_\nu=C_x=\omega_T (x).$   That $\omega_T (x)\subsetneqq \omega_T (x)$ is a contradiction.  Hence, $x$ is not in $QW (T)$.
\qed

\bigskip



{\bf Proof of Theorem \ref{More-general-systems-2} (2) }  By Lemma \ref{lem-DGS-full-support-2015-22222}, 
 there is some invariant measure with full support. By Lemma \ref{lem-PS}, $T$ is saturated. 
 By Lemma \ref{lem-entropy-dense-Ps-2005}, $T$ has entropy-dense property. 
  So Theorem \ref{More-general-systems-2} (2) can be deduced from Propositions \ref{Main-Thm-1}-\ref{Main-Thm-4}. \qed


\bigskip



\subsection{Proofs of  Theorem \ref{More-general-systems-2} (1) and Theorem  \ref{More-general-systems-3}}\label{Birkhoff-level-set}

We    divide  Theorem \ref{More-general-systems-3}
into several propositions  and then use it to prove
 Theorem  \ref{More-general-systems-2} (1). Here we will state every proposition for  possible more general systems.


\begin{Prop}\label{3Main-Thm-2015-R-minus-A00000}
 Let  $\,T$ be a
  continuous map of a compact metric space $X$. Suppose that $T$ is    single-saturated.
 Let  $\phi:X\rightarrow \mathbb{R}$ be a continuous function  with $I_\phi (T)\neq \emptyset.$     Then for any real number $a\in  Int( L_\phi)  ,$  \\

 (1) $\{    QR_{erg}(T), QR(T)\}$ has full entropy gaps with respect to $ R_{\phi,a} (T)$.\\
 

 (1')  $\{    R(T), QR(T)\}$  has full entropy gaps with respect to $ R_{\phi,a} (T)$. \\


 (2)  If there is some invariant measure   such that its support  is not minimal, then
 $\{   A(T)\cup QR_{erg}(T), QR(T)\}$ has full entropy gaps with respect to $ R_{\phi,a} (T)$. \\ 
 
 (2')  If there is some invariant measure   such that its support  is not minimal, then 
  $\{    A(T)\cup R(T), QR(T)\}$   has full entropy gaps with respect to $ R_{\phi,a} (T)$.  

\end{Prop}

{\bf Proof.}  Since $ R(T)\subseteq QR_{erg}(T)  $,  (1) implies (1') and (2) implies (2') so that we only need to prove (1) and (2). 
   Remark that from (\ref{Ir-phi-nonempty-imply-differentmeasure})  $I_\phi (T)\neq \emptyset $   implies that the system is naturally not uniquely ergodic.

Fix $a\in Int (L_\phi) $  and
let $t=\sup \{h_{\rho} (T)|\,\,\rho\in M (T,X)\,\,and\,\,\int \phi d\rho=a\}.$  For any $\epsilon>0,$ we can take $\rho\in M(T,X)$ such that $\int \phi d\rho=a,\,\,h_\rho(f)>t-\epsilon.$
\medskip

(1)   By  Proposition \ref{3Main-Thm-0}, we only need to show that  $$ h_{top} (T,R_{\phi,a} (T)\cap  QR(T)\setminus   QR_{erg}(T) )\geq t.$$
If   $\rho$ is non-ergodic, then take $m:=\rho.$ Otherwise, $\rho$ is ergodic,  by Lemma \ref{lem-phi-different} (6), take a non-ergodic measure $\sigma$ such that
  $ \int \phi d\sigma=a.$     Take $\theta'\in (0,1)$ close to 1 such that  $m:=\theta'\rho+(1-\theta')\sigma$ satisfies $h_m(f)\geq \theta' h_\rho(f) >t-\epsilon.$ Remark that $\int \phi dm=a$ and $m$ is not ergodic, since $\sigma$ is not ergodic.  Since  $T$ is single-saturated, then $h(T,G_m)=h_m(f) >t-\epsilon.$

  We only need to prove that     $G_m\subseteq R_{\phi,a} (T)\cap QR(T)\setminus QR_{erg}(T).$  By construction of $m,$  it is easy to see that    $G_m\subseteq     QR(T)\setminus QR_{erg}(T).$
 For any $x\in G_{m},$ $M_x(T)=\{m\}$ so that $M_x(T)\subseteq \{\vartheta |
 \int \phi d\vartheta=a\}.$  By (\ref{Regular-equalto-a}) $x\in R_{\phi,a}(T).$ Now we complete the proof of (1).
 \medskip

(2) Take same $m$ as in (1).  If $S_m$  is not minimal, then we claim that  (2) can be deduced from  $G_m\subseteq R_{\phi,a} (T)\cap QR(T)\setminus (QR_{erg}(T)\cup A(T)).$  Otherwise, there is $x\in G_m\cap A(T)$. Then $\omega_T(x)$ is a minimal compact set and $x\in\omega_T(x)$.  By (\ref{C_x=UnionSupport}) and (\ref{C_xinOmega}) $S_m=C_x\subseteq \omega_T(x)$,    contradicting that   $S_m$ is not minimal.

 Now we face the case that $S_m$ is minimal.
By  assumption there is some   invariant measure  $\mu$   with non-minimal support. Then $m\neq \mu.$
Let $b:=\int \phi d\mu .$  If $b=a,$ take $\omega:=\mu.$
If $b \neq a,$ without loss of generality, we can assume $b>a$. By definition of $L_\phi$ and connectedness   of $M(T,X)$, we can choose some $\nu$ such that $c:=\int \phi d\nu<a.$
 In this case we can take suitable $\theta\in(0,1)$ such that $\omega=\theta \mu +(1-\theta)\nu$ satisfies that $ \int \phi d\omega=a.$ So in any case, $\omega$ satisfies that $\int \phi d\omega=a$ and $S_\omega \supseteq S_\mu$ is not minimal.   Take $\theta'\in (0,1)$ close to 1 such that  $\tau:=\theta'm+(1-\theta')\omega$ satisfies $h_\tau(f)\geq \theta' h_m(f) >t-\epsilon.$ Remark that $\int \phi d\tau=a$ and $S_\tau=S_\omega\cup S_m\supseteq S_\omega$ is not minimal,  and $S_m$ is minimal implies that $m\neq \omega $  and so $\tau$ is non-ergodic. Since  $T$ is single-saturated, then $h(T,G_\tau)=h_\tau(f) >t-\epsilon.$

 We need to prove $G_\tau\subseteq R_{\phi,a} (T)\cap QR(T)\setminus (QR_{erg}(T)\cup A(T)).$  Similar as the proof of (1), it is easy to check that $G_\tau\subseteq R_{\phi,a} (T)\cap QR(T)\setminus  QR_{erg}(T).$ Now we start to prove $G_\tau \cap A(T)=\emptyset.$  Otherwise, there is $x\in G_\tau \cap A(T)$. Then $\omega_T(x)$ is a minimal compact set and $x\in\omega_T(x)$. 
By (\ref{C_x=UnionSupport}) and (\ref{C_xinOmega}) $S_\tau=C_x\subseteq \omega_T(x)$, contradicting that    $S_\tau$ is not minimal. Now we complete the proof of (2). \qed


\begin{Prop}\label{3Main-Thm-2015-R-minus-A11111}
 Let  $\,T$ be a
  continuous map of a compact metric space $X$. Suppose that $T$ is    single-saturated and    there is some invariant measure with full support.
 Let  $\phi:X\rightarrow \mathbb{R}$ be a continuous function  with $I_\phi (T)\neq \emptyset.$     Then for any real number $a\in  Int( L_\phi)  ,$  \\

 (1)   $\{    QR_{erg}(T), QR_d(T)\}$ has full entropy gaps with respect to $ W(T)\cap R_{\phi,a} (T)$.\\
 
 (1')   $\{    R(T), QR_d(T)\}$   has full entropy gaps with respect to $ W(T)\cap  R_{\phi,a} (T)$. \\

  (2)  If further the system $T$   is not minimal, then
 $\{   A(T)\cup QR_{erg}(T), QR_d(T)\}$ has full entropy gaps with respect to $ W(T)\cap R_{\phi,a} (T)$.\\
 
 (2')    $\{    A(T)\cup R(T), QR_d(T)\}$  has full entropy gaps with respect to $ W(T)\cap  R_{\phi,a} (T)$.  

\end{Prop}

{\bf Proof.}   Since $ R(T)\subseteq QR_{erg}(T)  $,  (1) implies (1') and (2) implies (2') so that we only need to prove (1) and (2).

Remark that from (\ref{Ir-phi-nonempty-imply-differentmeasure})  $I_\phi (T)\neq \emptyset $   implies that the system is naturally not uniquely ergodic.
By Lemma \ref{lem-phi-different} (2) there is some invariant measure $\xi$ with full support such that $\int \phi d\xi=a.$  Fix $a\in Int (L_\phi) $  and
let $t=\sup \{h_{\rho} (T)|\,\,\rho\in M (T,X)\,\,and\,\,\int \phi d\rho=a\}.$ Fix $\epsilon>0,$

 (1) Take $m$ same as in the proof of Proposition \ref{3Main-Thm-2015-R-minus-A00000} (1)  such that $m$ is not ergodic, $\int \phi dm=a,\,\,h_m(f)>t-\epsilon.$  Take $\alpha\in (0,1)$ close to 1 enough such that the measure $m'=\alpha m+(1-\alpha) \xi$ such that $h_{m'}(f)\geq \alpha h_m(f)>t-\epsilon.$ Remark that $m'$ is not ergodic, has full support and
 $\int \phi dm'=a.$ Since  $T$ is single-saturated, then $h(T,G_{m'})=h_{m'}(f) >t-\epsilon.$

 We only need to prove that     $G_{m'} \subseteq W(T)\cap R_{\phi,a} (T)\cap QR_d(T)\setminus QR_{erg}(T).$  By construction of $m',$  it is easy to see that    $G_{m'}\subseteq     QR(T)\setminus QR_{erg}(T).$      Now we start to show that $G_{m'}\subseteq QR_d(T).$ In fact,  $S_{m'}=X$ implies $G_{m'}= G_{m'}\cap S_{m'}$ and thus by (\ref{eq-density-quasiregular222222}) $G_{m'}\subseteq QR_d(T).$
 For any $x\in G_{m'},$ $M_x(T)=\{m'\}$ so that $M_x(T)\subseteq \{\rho|
 \int \phi d\rho=a\}.$
 Thus  by (\ref{Regular-equalto-a}) $x\in R_{\phi,a} (T).$  Moreover, for $x\in G_{m'},$ $M_x(T)=\{m'\}$ implies that  $C_x=S_{m'}=X\ni x.$ By (\ref{C_xinOmega}) $x\in Rec(T)$ and then  by (\ref{WT})  $x\in W(T)$.   Now we complete the proof of (1).
 \medskip

(2) By non-minimal assumption the measure $m'$ in (1) satisfies that $S_{m'}=X$ is not minimal. Similar as   the proof of Proposition \ref{3Main-Thm-2015-R-minus-A00000} (2) it is easy to check that  $G_{m'}\cap A(T)=\emptyset$.
 Thus  by (1)  $G_{m'}\subseteq W(T)\cap R_{\phi,a}(T)\cap QR_d(T)\setminus( A(T)\cup QR_{erg}(T)).$
  Then we can follow the proof of (1) to  complete the proof of (2).
\qed



\bigskip

\begin{Prop}\label{3Main-Thm-1} 
 Let  $\,T$ be a     continuous map of a compact metric space $X$. Suppose that $T$ is saturated, there is some invariant measure    with full support and $M_{erg}(T,X)$  is dense in $M(T,X).$    Let  $\phi:X\rightarrow \mathbb{R}$ be a continuous function with $I_\phi (T)\neq \emptyset.$  Then for any real number $a\in Int (L_\phi),$
$$h_{top} (T,R_{\phi,a} (T)\cap W (T)\setminus QR(T))=h_{top} (T,R_{\phi,a} (T)).$$

\end{Prop}

{\bf Proof.} Fix $a\in Int (L_\phi)$ and
let $t=\sup \{h_{\rho} (T)|\,\,\rho\in M (T,X)\,\,and\,\,\int \phi d\rho=a\}.$  By  Proposition \ref{3Main-Thm-0}, we only need to show that $$ h_{top} (T,R_{\phi,a} (T)\cap W (T)\setminus QR(T))\geq t.$$  Fix $\varepsilon>0.$  We need to construct two measures as follows, which are also useful to prove other propositions.

 \begin{Lem}\label{3Main-Thm-1-TwoFullsupport} Let  $\,T$ be a     continuous map of a compact metric space $X$.
  Suppose that $T$ is saturated, there is some invariant measure    with full support and $M_{erg}(T,X)$  is dense in $M(T,X).$
   Then there are two different measures $\omega,\omega'\in M (T,X)$$ (\omega\neq \omega')$ such that $$
\min\{ h_{\omega} (T)  ,\,\, h_{\omega'} (T) \}>t-\varepsilon$$ and $\int \phi d \omega=\int \phi d \omega'=a,\,\, S_{\omega}=S_{\omega'}=X.$
\end{Lem}

{\bf Proof.} By Lemma \ref{lem-phi-different} (4)  we can take three   different  ergodic measures of $\mu_i\, (i=0,1,2 )$ with support $X$ such that $$\int \phi d\mu_0<a<\int \phi d\mu_1<\int \phi d\mu_2.$$ Then we can choose suitable $\theta_i\in (0,1) (i=1,2)$ such that $\nu_i=\theta_i\mu_0+ (1-\theta_i)\mu_{i}$  satisfy $$\int\phi d\nu_i=a, i=1,2.$$ Remark that by ergodicity of $\mu_i,$ $\nu_1\neq \nu_2$ and $S_{\nu_i}=S_{\mu_0}\cup S_{\mu_i}=X, i=1,2.$


 By  definition of $t$,
 we can take  $\mu\in M (T,X)$ such that $\int \phi d\mu=a$ and  $h_\mu (T)>t-\varepsilon.$
 Then we can choose  $0<\theta<1$ close to 1 such that $\omega=\theta\mu + (1-\theta)\nu_1,\,\omega'=\theta\mu + (1-\theta)\nu_2$ satisfy $$
h_{\omega} (T)=\theta h_\mu (T) + (1-\theta)h_{\nu_1} (T)\geq \theta h_\mu (T)>t-\varepsilon,$$ $$h_{\omega'} (T)=\theta h_\mu (T) + (1-\theta)h_{\nu_2} (T)\geq \theta h_\mu (T)>t-\varepsilon.$$   Remark that $\int \phi d \omega=\int \phi d \omega'=a,\,\, S_{\omega}=S_{\omega'}=X$ and  $\nu_1\neq \nu_2$  implies $\omega\neq \omega'.$ \qed

\bigskip

Now we continue the proof of Proposition \ref{3Main-Thm-1}.
  Let $K=\{\tau \omega + (1-\tau)\omega'|\,\tau\in [0,1]\},$ then for any $m=\tau \omega + (1-\tau)\omega'\in K,$
$$S_m=X,\,\,\,h_m (T)\geq \min\{h_{\omega} (T),h_{\omega'} (T)\}>t-\varepsilon,\,\,\int\phi d m=a.$$
 Since $T$ is saturated, then
   $$h_{top} (T,G_K)=\inf\{h_m (T)\,|\,m\in K\}>t-\varepsilon,$$ where $G_K=\{x\in X|\,M_x (T)=K\}.$
We only need to prove $G_K\subseteq R_{\phi,a} (T)\cap W (T)\setminus QR(T).$


Fix $x\in G_K$. Then $ M_x (T)=K $ so that for any $ m\in M_x (T),$    $\int \phi dm=a$. Thus  $M_x (T)\subseteq \{\rho|\,\,\int \phi d\rho=a\}$ and so by (\ref{Regular-equalto-a}) $x\in R_{\phi,a} (T).$ Notice that for any $m\in M_x (T),$ by (\ref{C_x=UnionSupport}) $$C_x=\overline{\cup_{m\in M_x (T)}S_m}=\overline{\cup_{m\in K}S_m}=X=S_m.$$  By (\ref{C_xinOmega}) $x\in X=C_x\subseteq \omega_T (x)$ and thus $x\in Rec (T).$  So from (\ref{WT}) one has $x\in W (T).$ Since $M_x (T)=K$ is not a singleton, by (\ref{Irregular-equivalent}) $x\in I (T).$
\qed

\bigskip

\begin{Prop}\label{3Main-Thm-2} 
Let  $\,T$ be a      continuous map of a compact metric space $X$.  Suppose that  $T$ is saturated and has entropy-dense property,
 the periodic points are dense in $X$ (i.e., $\overline{Per(f)}=X$) and
the periodic measures are dense in the space of invariant measures (i.e., $\overline{ M_{p} (T,X) }=M (T,X)).$
Let  $\phi:X\rightarrow \mathbb{R}$ be a continuous function with $I_\phi (T)\neq \emptyset.$  Then for any real number $a\in Int (L_\phi),$
$$h_{top} (T, R_{\phi,a} (T)\cap V(T) \setminus W (T) )=h_{top} (T,R_{\phi,a} (T)).$$


\end{Prop}

{\bf Proof.} Fix $a\in Int( L_\phi (T))$ and
let $t=\sup \{h_{\rho} (T)|\,\,\rho\in M (T,X)\,\,and\,\,\int \phi d\rho=a\}.$  By  Proposition \ref{3Main-Thm-0}, we only need to show that $$ h_{top} (T,R_{\phi,a} (T)\cap V (T)  \setminus W (T))\geq t.$$
Fix $\varepsilon\in (0,t).$  We need to construct a measure as follows, which is also useful to prove other propositions.

 \begin{Lem}\label{3Main-Thm-2-NOT-Fullsupport} Let  $\,T$ be a      continuous map of a compact metric space $X$. 
  Suppose that  $T$ is saturated and has entropy-dense property, the periodic points are dense in $X$ (i.e., $\overline{Per(f)}=X$) and
the periodic measures are dense in the space of invariant measures (i.e., $\overline{ M_{p} (T,X) }=M (T,X)).$
Then there is a measure $\omega \in M (T,X)$ such that $$
h_{\omega} (T) >t-\varepsilon  \,\,\, and \,\,\,\int \phi d \omega =a,\,\, S_{\omega} \subsetneqq X.$$

\end{Lem}

 {\bf Proof.} Take a $\nu\in M (T,X)$ such that $h_{\nu} (f)>t-\frac{\varepsilon}3$ and $\int \phi d\nu=a.$ By Lemma \ref{lem-phi-different} (1), take two periodic measures $\nu_i (i=1,2)$ such that $b_1:=\int \phi d\nu_1>a>\int \phi d\nu_2=:b_2.$ Let $\delta>0$ small enough such that $$\min{\{\frac{b_1-a}{b_1-a+\delta},\frac{a-b_2}{a-b_2+\delta}\}}>\frac{t-\varepsilon}{t-\frac{2\varepsilon}3}.$$
 Then by entropy-dense property, we can take one ergodic measure $\mu$  close to $\nu$ enough (in weak$^*$ topology) such that $S_\mu\subsetneqq X$ and
$$|\int \phi d\mu-a|=|\int \phi d\mu-\int \phi d\nu|<\delta,\,\, h_\mu (f)> t-\frac{2\varepsilon}3.$$
If $\int \phi d\mu=a,$ then take $\omega=\mu.$ Otherwise, $\int \phi d\mu\neq a.$ Without loss of generality, we assume $\int \phi d\mu<a$. Take $$\omega=\frac{b_1-a}{b_1-\int \phi d\mu}\mu+ (1-\frac{b_1-a}{b_1-\int \phi d\mu})\nu_1.$$ Then $\int \phi d\omega =a,\,\,h_\omega (f)\geq \frac{b_1-a}{b_1-\int \phi d\mu}h_\mu (f)>\frac{b_1-a}{b_1-a+\delta}h_\mu (f)>t-\varepsilon.$ Recall that $S_{\nu_1}$ is a finite closed set but $X\setminus S_\mu$ is nonempty,  open  and thus by Proposition \ref{prop-number-periodic-in-openset}  it is not a finite set. So $S_\omega=S_\mu\cup S_{\nu_1}\neq X.$ \qed

\bigskip

Now we continue the proof of Proposition \ref{3Main-Thm-2}. By density of periodic points and Lemma \ref{lem-DGS-full-support-2015-22222}, there is some invariant measure    with full support and by density of periodic measures,  $M_{erg}(T,X)$  is dense in $M(T,X).$    Then one can construct  $\omega'$ same  as in Lemma \ref{3Main-Thm-1-TwoFullsupport} such that $$h_{\omega'} (T)>t-\varepsilon,\int\phi d\omega'=a\,\,\,{\text {and} }\,\,\,S_{\omega'}=X.$$
Clearly $S_\omega\neq S_{\omega'}$ and thus $\omega\neq \omega'.$

 Let  $K=\{\tau \omega + (1-\tau)\omega'|\,\tau\in [0,1]\},$ then    for any $m=\tau \omega + (1-\tau)\omega'\in K\setminus\{\omega\},$
$S_m=X$ and for any $m=\tau \omega + (1-\tau)\omega'\in K,$ $$h_m (T)\geq \min\{h_{\omega} (T),h_{\omega'} (T)\}>t-\varepsilon,\,\,\int\phi d m=a.$$
   Since $T$ is saturated, then
   $$h_{top} (T,G_K)=\inf\{h_m (T)\,|\,m\in K\}>t-\varepsilon,$$ where $G_K=\{x\in X|\,M_x (T)=K\}.$
We only need to prove $$G_K\subseteq R_{\phi,a} (T)\cap V (T)  \setminus W (T).$$

Fix $x\in G_K$. Then $ M_x (T)=K $ so that for any $ m\in M_x (T),$    $\int \phi dm=a$. Thus  $$M_x (T)\subseteq \{\rho|\,\,\int \phi d\rho=a\}$$
and so by (\ref{Regular-equalto-a}) $x\in R_{\phi,a} (T).$ Notice that by (\ref{C_x=UnionSupport}) $$C_x=\overline{\cup_{m\in M_x (T)}S_m}=\overline{\cup_{m\in K}S_m}=X.$$
So by (\ref{C_x=X}) $C_x=X$ implies $x\in QW (T).$ By (\ref{weakperiodic-in-reccurent}) $x\in Rec (T).$ Then from (\ref{WT}) $C_x=X\neq S_{\omega}$
and ${\omega}\in M_x (T)$ imply $x\in X\setminus W (T)$.   Recall $S_{\omega'}=X$ and $\omega'\in K=M_x (T).$
So  $x\in   R_{\phi,a} (T)\cap V (T)  \setminus W (T).$ \qed 

\bigskip

\begin{Prop}\label{3Main-Thm-3} Let  $\,T$ be a      continuous map of a compact metric space $X$.  Suppose that  $T$ is saturated and has entropy-dense property,
 the periodic points are dense in $X$ (i.e., $\overline{Per(f)}=X$) and
the periodic measures are dense in the space of invariant measures (i.e., $\overline{ M_{p} (T,X) }=M (T,X)).$
Let  $\phi:X\rightarrow \mathbb{R}$ be a continuous function with $I_\phi (T)\neq \emptyset.$  Then for any real number $a\in Int (L_\phi),$
$$h_{top} (T, R_{\phi,a} (T)\cap    QW (T)\setminus V(T) )=h_{top} (T,R_{\phi,a} (T)).$$

\end{Prop}

{\bf Proof.} Fix $a\in L_\phi (T)$ and
let $t=\sup \{h_{\rho} (T)|\,\,\rho\in M (T,X)\,\,and\,\,\int \phi d\rho=a\}.$ By  Proposition \ref{3Main-Thm-0}, we only need to show that $$ h_{top} (T,R_{\phi,a} (T)\cap QW (T)  \setminus V (T))\geq t.$$
Fix $\varepsilon>0.$ Firstly, by assumption we can take  same $\omega$ as in Lemma \ref{3Main-Thm-2-NOT-Fullsupport} such that  $\int \phi d\omega =a,\,\,h_\omega (f)>t-\varepsilon$ and $S_\omega\neq X.$

By density of periodic points and density of periodic measures, we can choose
  $\{x_i\}_{n=1}^\infty$ and $\{m_i\}_{n=1}^\infty$  same as in the proof Proposition \ref{Main-Thm-3}. They are composed of periodic points  and periodic measures and $\overline{\{x_i\}_{i=1}^\infty}=X,\,\,\overline{\{m_i\}_{i=1}^\infty}=M (T,X),\,\,\cup_{i\geq 1}S_{m_i}=\{x_i\}_{i=1}^\infty.$

  Let $$K_1:=\{m|\int \phi dm>a, \,m\in \{m_i\}_{n=1}^\infty \},$$  $$K_2:=\{m|\int \phi dm<a, \,m\in \{m_i\}_{n=1}^\infty \}$$ and $$K_3:=\{m|\int \phi dm=a, \,m\in \{m_i\}_{n=1}^\infty \}.$$ By Lemma \ref{lem-phi-different} if let $\mathcal{R}=\{m_i\}_{n=1}^\infty$, then  it is easy to see that $K_1$ and $K_2$ are countable. Remark that $K_3$ may be empty, finite or countable. Without loss of generality,  we can assume $K_i=\{m^{ (i)}_{j}\}_{j=1}^\infty, \,i=1,2,3.$   Then we can choose suitable $\theta_{j,k}\in (0,1)$ such that $m_{j,k}=\theta_{j,k}m^{ (1)}_{ j}+ (1-\theta_{j,k})m^{ (2)}_{k}$ satisfies $\int\phi dm_{j,k}=a.$ For any $n\geq 1$, let $$l_n=\frac{\sum_{j+k=n}m_{j,k}+m^{ (3)}_{n}}n,$$ then $\int\phi dl_n=a.$
Remark that every $S_{l_n}$ is composed of finite periodic orbits and  $\bigcup_{n\geq 1} S_{l_n}=\{x_i\}_{n=1}^\infty$ is dense in $X.$

    Take an increasing sequence of  $\{\theta_i|\,\,\theta_i\in (0,1)\}_{i=1}^\infty$ convergent to 1 such that $h_{\omega_i} (T)>t-\varepsilon$ where $\omega_i=\theta_i\omega+ (1-\theta_i)l_i.$ Remark that $S_{\omega_i}=S_\omega\cup S_{l_i}.$ In particular,  for all $i$,  $\int\phi d\omega_i=\int\phi d\omega=a.$

Now we consider $$K=\{\omega\}\cup\bigcup_{i\geq 1} \{\tau \omega_i+ (1-\tau)\omega_{i+1}|\,\tau\in[0,1]\}.$$ Then $K$ is nonempty connected compact subset of $M (T,X)$ because $\omega_i\rightarrow \omega$ in weak$^*$ topology. Since $T$ is saturated, then  $$h_{top} (T,G_K)=\inf\{h_\nu (T)\,|\,\nu\in K\}=\min\{\inf_{i\geq 1}\{h_{\omega_i} (T) \},h_\omega (T)\}\geq  t-\varepsilon,$$ where $G_K=\{x\in X|\,M_x (T)=K\}.$ We only need to prove $$G_K\subseteq \{x\in R_{\phi,a} (T)\cap QW (T)\setminus W (T)|\forall\, m\in M_x (T) \,s.t. \,S_m\neq C_x\}.$$

Fix  $x\in G_K$. Note that for any $m \in  M_x (T)=K$,  $\int\phi dm=a$. Thus  $M_x (T)\subseteq \{\rho|\,\,\int \phi d\rho=a\}$ and so by (\ref{Regular-phi-equivalent})  $x\in R_{\phi,a} (T).$ Notice that by (\ref{C_x=UnionSupport}) $$C_x=\overline{\cup_{m\in M_x (T)}S_m}=\overline{\cup_{m\in K}S_m}\supseteq\overline{\cup_{i\geq 1}S_{\omega_i}} =\overline{\cup_{i\geq 1} (S_{l_i}\cup S_\omega)}\supseteq\overline{\cup_{i\geq 1}\{x_i\}}=X.$$ So $  X=C_x$. From (\ref{C_x=X}) one has $x\in QW (T).$ By (\ref{weakperiodic-in-reccurent}) $x\in Rec (T).$ Notice that from (\ref{WT}) $C_x=X\neq S_\omega$ and $\omega\in M_x (T)$ imply $x\in X\setminus W (T)$.  For any $ m\in M_x (T)=K,$ there is some $i$ such that $S_m\subseteq S_\omega\cup S_{l_i}\cup S_{l_{i+1}} \subsetneqq X=C_x$. This is because $S_{l_i}\cup S_{l_{i+1}}$ is a finite set but $X \setminus S_\omega$ is nonempty,  open and thus by Proposition \ref{prop-number-periodic-in-openset} it is an infinite set. Hence   $S_m\neq C_x.$
\qed

\bigskip

\begin{Prop}\label{3Main-Thm-4} Let  $\,T$ be a      continuous map of a compact metric space $X$.  Suppose that  $T$ is saturated and has entropy-dense property,
 the periodic points are dense in $X$ (i.e., $\overline{Per(f)}=X$) and
the periodic measures are dense in the space of invariant measures (i.e., $\overline{ M_{p} (T,X) }=M (T,X)).$
Let  $\phi:X\rightarrow \mathbb{R}$ be a continuous function with $I_\phi (T)\neq \emptyset.$  Then for any real number $a\in Int (L_\phi),$
$$
h_{top} (T,R_{\phi,a} (T)\cap I (T)\setminus QW (T)) =h_{top} (T,R_{\phi,a} (T)).$$

\end{Prop}

{\bf Proof.}
Fix $a\in L_\phi (T)$ and
let $t=\sup \{h_{\rho} (T)|\,\,\rho\in M (T,X)\,\,and\,\,\int \phi d\rho=a\}.$ By  Proposition \ref{3Main-Thm-0}, we only need to show that $$ h_{top} (T,R_{\phi,a} (T)\cap I (T)  \setminus QW (T))\geq t.$$
Fix $\varepsilon>0.$  Firstly, by assumption we can take same $\omega$ as in Lemma \ref{3Main-Thm-2-NOT-Fullsupport} such that  $\int \phi d\omega =a,\,\,h_\omega (f)>t-\varepsilon$ and $S_\omega\neq X.$

Let $B:=\{m\in M_p (T,X)|\,S_m\setminus S_\omega\neq \emptyset\}.$ Since $X\setminus S_\omega$ is  nonempty, open and invariant, by density of periodic points,   $B\neq \emptyset.$ If there is $m\in B$ such that $\int \phi d m=a,$ take $\mu=m.$
Otherwise, for any $m\in B,$ $\int \phi dm\neq a.$ Take one $\mu_1\in B.$ Then  $\int \phi d\mu_1\neq a.$ Without loss of generality, we assume $\int \phi d\mu_1<a.$ By Lemma \ref{lem-phi-different} (1) we can take a periodic measure $\mu_2$
 such that $\int \phi d\mu_2>a.$  Then we can choose suitable $\theta\in (0,1)$ such that $\mu=\theta\mu_1+ (1-\theta)\mu_2$ satisfies  $ \int \phi d \mu =a. $ Remark that 
  $ S_\mu\setminus S_\omega  $    is composed of one periodic orbit or two periodic orbits containing $S_{\mu_1}$.

Take $\theta'\in (0,1)$ close to 1 such that   $\omega'=\theta'\omega+ (1-\theta')\mu$  satisfies $h_{\omega'} (T)>t-\varepsilon.$ Remark that $\int \phi d \omega'=a$ and $S_{\omega'}\setminus S_\omega=S_\mu\setminus S_\omega$ is composed of one periodic orbit or two periodic orbits. So  $S_{\omega'}\setminus  S_{\omega}$ is nonempty, finite,   invariant and compact and thus $\omega\neq \omega'.$

 Let $K=\{\tau \omega + (1-\tau)\omega'|\,\tau\in [0,1]\},$ then  for any $m\in K,$ $$h_m (T)\geq \min\{h_\omega (T),h_{\omega'} (T)\}>t-\varepsilon,\,\,\int\phi dm=a.$$
  Since $T$ is saturated, then
   $$h_{top} (T,G_K)=\inf\{h_m (T)\,|\,m\in K\}>t-\varepsilon,$$ where $G_K=\{x\in X|\,M_x (T)=K\}.$
We only need to prove $G_K\subseteq R_{\phi,a} (T)\cap I (T)\setminus QW (T).$

Fix  $x\in G_K$. Since $K$ is not singleton, by (\ref{Irregular-equivalent}) $x\in I (T).$ Note that for any $m \in  M_x (T)=K$,  $\int\phi dm=a$. Thus  $M_x (T)\subseteq \{\rho|\,\,\int \phi d\rho=a\}$ and so by (\ref{Regular-phi-equivalent})  $x\in R_{\phi,a} (T).$   If $x\in QW (T)$, by (\ref{weakperiodic-in-reccurent}) $x\in Rec (T).$ By (\ref{QW})   $x\in C_x=\omega_T (x).$ Then  by (\ref{C_x=UnionSupport}) $$x\in \omega_T (x)=C_x=\overline{\cup_{m\in M_x (T)}S_m}=\overline{\cup_{m\in K}S_m}=S_\omega \sqcup (S_{\omega'}\setminus  S_{\omega}).$$ So  $x \in S_\omega$ or $S_{\omega'}\setminus  S_{\omega}.$  Recall  that $S_\omega, S_{\omega'}\setminus  S_{\omega}$ are both  compact and invariant.   One has $Orb  (x)\subseteq S_\omega$ or $Orb  (x)\subseteq S_{\omega'}\setminus  S_{\omega}$ and thus   $\omega_T (x)\subseteq S_\omega\subsetneqq S_\omega \sqcup (S_{\omega'}\setminus  S_{\omega})=C_x=\omega_T (x)$ or $\omega_T (x)\subseteq (S_{\omega'}\setminus  S_{\omega})\subsetneqq S_\omega \sqcup (S_{\omega'}\setminus  S_{\omega})=C_x=\omega_T (x).$ That $\omega_T (x) \subsetneqq \omega_T (x)$ is a contradiction.    Hence, $x$ is not in $QW (T)$.
\qed



\bigskip

{\bf Proof of Theorem \ref{More-general-systems-3} } Since every periodic measure is  ergodic, by density of periodic measures, ergodic measures are   dense in the space of invariant measures.  Since $X$ is assumed infinite in this paper, by density of periodic points, $X$ is not a minimal set (in other words, $T$ is not minimal). By density of periodic measures, ergodic measures are dense in the space of invariant measures.  By Lemma \ref{lem-DGS-full-support-2015-22222}, there is some invariant measure such that it has  full support and in particular, the support is not minimal. By Lemma \ref{lem-PS}, $T$ is saturated. By Lemma \ref{lem-entropy-dense-Ps-2005}, $T$ has entropy-dense property (which also implies density of ergodic measures).   So Theorem  \ref{More-general-systems-3}  can be deduced from Propositions   \ref{3Main-Thm-2015-R-minus-A11111}, \ref{3Main-Thm-1}, \ref{3Main-Thm-2}, \ref{3Main-Thm-3}, \ref{3Main-Thm-4}. \qed

\bigskip

\begin{Rem}\label{Rem-phi-Birkhorf-regular}
 In Theorem  \ref{More-general-systems-3}, we require  $a$ to satisfy  $$\inf\{\int\phi d\mu|\,\mu\in M (T,X)\}<a<\sup\{\int\phi d\mu|\,\mu\in M (T,X)\}.$$
  For full shifts of finite type, let us give an example why we do not choose   $a$ from  extreme points.
   More precisely, firstly recall a result from   \cite{Gri,HaKa} that for any full shift  on  finite symbols, there exist uniquely ergodic minimal subshifts with any given entropy.
    Then we can take an ergodic measure $\mu$  with positive entropy such that $S_\mu$ is minimal, $S_\mu\neq X$ and $\mu$ is the unique invariant measure supported on $S_\mu.$  By density of periodic points, take a periodic measure $\nu$ such that $S_\nu\subseteq X\setminus S_\mu.$
 Since $S_\mu$ and $S_\nu$ are two disjoint closed sets, we can take   a continuous function  $\phi$ which restricts on the set $S_\mu$ with value $0$ and on the set $S_\nu$ with value $1$ respectively, and the values of other points are in the open interval of $ (0,1).$  Take  $a=0,$  let us show that $R_{\phi,a}(T)$ is nonempty and has positive entropy but $R_{\phi,a}(T)\cap I(T)=\emptyset.$  More precisely, on one hand,  
 by (\ref{Regular-equalto-a}),  $R_{\phi,a}(T)\supseteq G_\mu$. By ergodicity of $\mu$ and (\ref{eq-Bowen-ergodicentropy-equal-generic}), $h_{top}(T,G_\mu)=h_\mu(f)>0.$ So $h_{top}(T,R_{\phi,a}(T))\geq h_{top}(T,G_\mu)>0.$
On the other hand, observe that   $\mu$ is the unique invariant measure with integral by $\phi$ equal to $a$. So by weak$^*$ topology, $R_{\phi,a}(T)=G_\mu$  and thus $R_{\phi,a}(T)$ does not contain irregular point.

\end{Rem}

\bigskip

{\bf Proof of Theorem \ref{More-general-systems-2} (1) }

 {\bf Step 1.} Firstly we consider the case of  $I_\phi (T)\neq\emptyset.$  That is, we need to prove that
  $$\{A(T)\cup R(T), QR(T), W(T), V(T), QW(T), I(T)\}$$  has full entropy gaps with respect to $R_\phi(T)$.

 Fix $\varepsilon>0.$ By   Lemma \ref{Fact-concave}, we can take a number $a\in Int (L_\phi)$ such that $$h_{top} (T, R_{\phi,a} (T))> h_{top} (T)-\varepsilon.$$
 Recall that   $$R_\phi (T)=\bigsqcup_{b\in\mathbb{R}}R_{\phi,b} (T).$$ So by Theorem \ref{More-general-systems-3}, $$h_{top} (T,R_\phi (T)\cap \xi)\geq h_{top} (T, R_{\phi,a} (T)\cap \xi)=h_{top} (T, R_{\phi,a} (T)) > h_{top} (T)-\varepsilon,$$ where
$$
\xi = QR(T)\setminus (A(T)\cup R(T)),   W (T)\setminus QR(T), \,V(T)\setminus W (T) ,
  \,\  QW (T)\setminus V (T) ,
  \,I (T)\setminus QW (T).$$
By arbitrariness of $\varepsilon$, we complete the proof  of Step 1.

\bigskip

{\bf Step 2.}  The case that  $I_\phi (T)$ is empty. That is,  $R_\phi (T)=X$.

  Since $X$ is assumed infinite in this paper, by density of periodic measures, $M(T,X)$ is not a singleton. Take two invariant measures $\mu_1\neq \mu_2$. By weak$^*$
topology, there is some continuous function $\varphi:X\rightarrow \mathbb{R}$ such that $\int \varphi d \mu_1\neq\int \varphi d \mu_2. $ Since $T$ has   $g$-almost product  property (which is a little stronger  than  almost specification, see \cite{Tho2012}), by Lemma \ref{Lem-IC-notempty}, we have $I_\varphi(T)\neq \emptyset.$ By Step 1 for the function $\varphi,$  $$\{A(T)\cup R(T), QR(T), W(T), V(T), QW(T), I(T)\}$$ has full entropy gaps with respect to $R_\varphi(T)$. Since $R_\varphi(T)\subseteq X=R_\phi(T),$ then $$\{A(T)\cup R(T), QR(T), W(T), V(T), QW(T), I(T)\}$$ also has full entropy gaps with respect to $R_\phi(T)$.

Now we complete the proof of Theorem \ref{More-general-systems-2} (1). \qed

\bigskip

Remark that Step 2 also can be as a proof of  Theorem \ref{More-general-systems}.

\subsection{Time-$t$ map of hyperbolic flows}\label{section-partial-hyperbolic}
Let $f:X\rightarrow X$ be the  time-$t$ ($t\neq 0$) map of a topologically mixing  Anosov flow of  a compact Riemannian manifold $X$. In this case, $f$ is partially hyperbolic with one-dimension central bundle. Then $f$ is far from tangency so that  $f$ is entropy-expansive   (see \cite{LiaoVianaYang} or see   \cite{DFPV,PacVie}). Recall that from   \cite{Misiurewicz} entropy-expansive implies  asymptotically $h-$expansive and from   Theorem 3.1  of   \cite{PS} any expansive or asymptotically $h-$expansive system satisfies uniform separation property. Recall that the unique maximal entropy measure of the flow has full support and note that the invariant measure of the flow is also invariant  for the time-$t$  map.
  On the other hand, as said in Section 4.3 of \cite{To2010}, the time-$t$  map $f$ satisfies specification property (even though possible not Bowen's specification) which is stronger than $g$-almost product property. So $f$ satisfies entropy-dense property which implies density of ergodic measures in the space of invariant measures.
 Thus, $f$ satisfies the assumptions of Proposition \ref{Main-Thm-1}, \ref{Main-Thm-2},  \ref{3Main-Thm-2015-R-minus-A11111} (2) and  \ref{3Main-Thm-1} and so  we have a following  result.

\begin{Thm}\label{Thm-Application-timetmap-of-flow-Irregular-full-entropy}
   Let $f:X\rightarrow X$ be a  time-$t$ map($t\neq 0$) of a topological mixing  Anosov flow of  a compact Riemannian manifold $X$ (in this case, $f$ is partially hyperbolic whose central bundle only has zero Lyapunov exponents). Then  \\

  (A) $\{A(T)\cup R(T), QR(T), W(T),V(T) \}$  has full entropy gaps with respect to $X$;\\

  (B)  for any continuous function $\phi:X\rightarrow \mathbb{R}$,
    $\{A(T)\cup R(T), QR(T), W(T)    \}$ has full entropy gaps with respect to $  R_\phi(T)$;\\

 (C) for any continuous function $\phi:X\rightarrow \mathbb{R}$ satisfying
  $I_\phi(T)\neq \emptyset$,      $\{ QR(T), W(T),$ $ V(T)  \}$  has full entropy gaps with respect to $  I_\phi(T)$; \\

 (D)    for any continuous function $\phi:X\rightarrow \mathbb{R}$ satisfying $I_\phi(T)\neq \emptyset$ and for any $a\in Int (L_\phi),$
    $ \{ A(T)\cup R(T), QR(T), W(T) \}$  has full entropy gaps with respect to $   R_{\phi,a}(T)$.
 \end{Thm}

{\bf Proof.} Recall from above analysis, $f$ satisfies specification and uniform separation, and there is an invariant measure with full support. Now we start to prove.

(C)  Since $f$ satisfies the assumptions of Proposition \ref{Main-Thm-1}, \ref{Main-Thm-2}, then (C) is valid.

\bigskip

(D)  Since $f$ satisfies the assumptions of    \ref{3Main-Thm-2015-R-minus-A11111} (2) and  \ref{3Main-Thm-1}, then (D) is valid.

\bigskip

(B) One can follow  same idea as the proof of  Theorem \ref{More-general-systems-2} (1) to show that (D) implies (B).


\bigskip

(A) If take $\phi\equiv 1$, then $R_\phi(T)=X$ and thus (B) implies (A) except $V(T)\setminus W(T).$  Notice that Anosov flow carries  periodic orbit.  Take a periodic measure  for  flow, then it is also invariant (even though not ergodic) measure  for time-$t$ map $f$. This measure is just supported on one periodic orbit of flow and thus is different from  the maximal entropy measure which has full support. In other words, $f$ is not uniquely ergodic.
By weak$^*$ topology, there is a continuous function $\phi:X\rightarrow \mathbb{R}$ such that $$\inf_{\mu\in M(T,X)}\int  \phi(x)d\mu<\sup_{\mu\in  M(T,X)}\int \phi(x)d\mu$$ So by Lemma \ref{Lem-IC-notempty} $ I_\phi(T)\neq \emptyset.$ By (C), $ \{ W(T),V(T)\}$ has
full   entropy gaps with respect to $I_\phi(T)$. From Theorem  \ref{T-9} we know that $I_\phi(T)$ carries
full topological entropy.  So $V(T)\setminus W(T)$ also carries
full topological entropy and thus the proof of (A) is finished. \qed

\bigskip

  Theorems 17.6.2 and 18.3.6 in \cite{Katok1} (originally due to Anosov)  ensure that the geodesic flow of any compact connected Riemannian manifold of
  negative sectional curvature is topologically mixing and Anosov. So Theorem \ref{Thm-Application-timetmap-of-flow-Irregular-full-entropy} can be applicative to
 the  time-$t$ map of the geodesic flow of any compact connected Riemannian manifold of
  negative sectional curvature.

\section{Transitive Points}\label{section-transitive}

 Let $D (T)$ denote  the set of all transitive points (i.e. the points whose orbit is dense in the whole space). In other words, $D(T)=\{x\in X|\,\overline{Orb(x)}=X\}.$
 If $T$ is a minimal system (for example, irrational rotation on the circle), it is obvious that $D(T)=A(T)=X.$
 Now we consider non-minimal systems.

 Firstly we state  some simple observation. For any dynamical system which is not minimal, $D (T)$ does not contain any periodic points and almost periodic points. Here we admit $D(T)$ is empty set (for example, the system of  rational rotation on the circle).

\begin{Lem}\label{lem-transipoint-disjoint-AlmostPeriodc-non-minimalsystem}
For any  non-minimal  continuous   map $T:X\rightarrow X$    of a compact
metric space $X$,   $A(T)\cap D(T)=\emptyset.$
\end{Lem}
{\bf Proof.}  By contradiction, there is some  $x\in A (T)$ such that the orbit of $x$ is dense in   $X.$   $x\in A (T)$ implies that $x\in \omega_T(z)$ and   the closed invariant set $\omega_T (x)$ is minimal. $x\in \omega_T(x)$ implies that $\omega_T(x)=\overline{Orb(x)}$.  Then    $X= \overline{Orb(x)}= \omega_T (x)$  is minimal,    contradicting that $T$ is not minimal.\qed

\bigskip

 We will coordinate  $D(T)$  with $R(T)\setminus A(T)$ together.

\begin{Thm}\label{T11}  {\it For any non-minimal continuous   map $T:X\rightarrow X$    of a compact
metric space $X$, if  there is an ergodic measure $\mu$ with maximal entropy and full support, then    $D(T)\cap  R(T)\setminus A (T) $  carries  full topological entropy. In particular,   each one of $  R(T)\setminus A (T) $,  $ D (T)=D (T)\setminus A (T) $ and $ R(T)\cap D (T)$ carries  full topological entropy.}

   \end{Thm}

Before proof  let us first give some  simple observation as follows.

\begin{Lem}\label{lem-transipoint-contains-generic-pt-meausre-withfullsupp}
For any   continuous   map $T:X\rightarrow X$    of a compact
metric space $X$,
$$ \{x\in X\,|\, C_x=X\}\subseteq D(T).$$
In particular,\\
(1)  for any nonempty connected compact set $K\subseteq M(T,X)$, if   
$\overline{\cup_{\mu\in K}S_\mu}=X,$ 
then  $G_K\subseteq D(T),$
where $G_K=\{x\in X|\,M_x(T)=K\}.$  \\
(2) for any invariant measure $\mu\in M(T,X),$ if $\mu$ has full support, then $G_\mu\subseteq D(T).$
\end{Lem}

{\bf Proof.}  Fix a point $x\in X$ such that $  C_x=X.$ By (\ref{C_xinOmega}), $x\in X=C_x\subseteq \omega_T(x).$ So $ \overline{Orb(x)}=X,$ i.e., $x\in D(T).$

In particular, for any nonempty connected compact set $K\subseteq M(T,X)$, if
$\overline{\cup_{\mu\in K}S_\mu}=X,$  by (\ref{C_x=UnionSupport}) $G_K\subseteq \{x|\, C_x=X\}$ and thus 
$G_K\subseteq D(T).$ This implies (1). For the case (2), it is obvious from (1) by taking  singleton $K$. \qed


\bigskip


{\bf Proof
  of Theorem \ref{T11}.}
 Since $\mu$ is ergodic, by Birkhoff ergodic theorem  $G_\mu$ is of $\mu$ full measure.
  Since  $S_\mu=X$, from Lemma \ref{lem-transipoint-contains-generic-pt-meausre-withfullsupp}  we know that    $G_\mu\cap S_\mu= G_\mu\subseteq D (T)$.


   Remark that  $G_\mu= G_\mu\cap  S_\mu\cap D (T)\subseteq R(T)\cap D (T)=  R(T)\cap D(T)\setminus A (T),$ since by   non-minimal assumption  and Lemma \ref{lem-transipoint-disjoint-AlmostPeriodc-non-minimalsystem}   $A(T)\cap D(T)=\emptyset$.    Then by (\ref{Y1containY2})
      we obtain that $D(T)\cap  R(T)\setminus A (T) $ carries  full topological entropy. We complete the proof.
\qed

   \bigskip

From Theorem \ref{QW-W-in-Irregular},   $R(T) \subseteq W (T)$ and thus
$$D(T)\cap R(T)\setminus    A (T)\subseteq    D(T)\cap  W(T) \setminus A (T) .$$
So by Theorem \ref{T11} we have a following consequence.

\begin{Thm}\label{T12} {\it  For any non-minimal continuous   map $T:X\rightarrow X$    of a compact
metric space $X$, if  there is an ergodic measure $\mu$ with maximal entropy and full support, then    $  D(T)\cap (W (T)\setminus A (T))$ carries  full topological entropy
  and so does $W(T)\setminus A(T)$.}
\end{Thm}

\begin{Rem}\label{Remark-transitive-case-RminusA}
Remark that
Theorem \ref{T11} and Theorem \ref{T12}  are applicative to  all the examples in  Section \ref{section-examples},
 since each example    is non-minimal and  the  unique maximal entropy measure  is ergodic and has full support.
\end{Rem}



 \bigskip

Now we consider $D(T)$ under same assumption of our main theorems in first section.

 \begin{Thm}\label{Thm-transitive-More-general-systems}
 Let  $\,T$ be a      continuous map of a compact metric space $X$ with $g-$almost product property and uniform  separation property. If the periodic points are dense in $X$ (i.e., $\overline{Per(f)}=X$) and the periodic measures are dense in the space of invariant measures (i.e., $\overline{ M_{p} (T,X) }=M (T,X)),$ then \\

  (A) $\{A(T)\cup R(T), QR(T), W(T), V(T), QW(T)\}$ has full entropy gaps with respect to $D(T)$;\\

  (B)  for any continuous function $\phi:X\rightarrow \mathbb{R}$,
    $\{A(T)\cup R(T), QR(T), W(T), V(T),$ $ QW(T) \}$ has full entropy gaps with respect to $D(T)\cap R_\phi(T)$;\\

 (C) for any continuous function $\phi:X\rightarrow \mathbb{R}$ satisfying
  $I_\phi(T)\neq \emptyset$,      $\{ QR(T), W(T),$ $ V(T),$ $ QW(T) \}$  has full entropy gaps with respect to $D(T)\cap  I_\phi(T)$; \\

 (D)    for any continuous function $\phi:X\rightarrow \mathbb{R}$ satisfying $I_\phi(T)\neq \emptyset$ and for any $a\in Int (L_\phi),$
    $ \{ A(T)\cup R(T), QR(T), W(T), V(T), QW(T) \}$  has full entropy gaps with respect to $D(T)\cap  R_{\phi,a}(T)$.
 \end{Thm}

\begin{Rem}\label{Remark-transitive-case}
Remark that
Theorem \ref{Thm-transitive-More-general-systems} hold for all the examples in  Section \ref{section-examples}. However, for the results of
$QW(T)\setminus V(T)$ and $I(T)\setminus QW(T),$  it is still unknown.
\end{Rem}

{\bf Proof.}  Since $X$ is assumed infinite in this paper, by density of periodic measures, $X$ is not a minimal set (in other words, $T$ is not minimal).  By Lemma \ref{lem-DGS-full-support-2015-22222}, there is some invariant measure such that it has  full support and in particular, the support is not minimal. By Lemma \ref{lem-PS}, $T$ is saturated. By Lemma \ref{lem-entropy-dense-Ps-2005}, $T$ has entropy-dense property.

\bigskip

Case (C).   Observe that    from the proofs of  Proposition  \ref{Main-Thm-1}, \ref{Main-Thm-2} and \ref{Main-Thm-3},    each  constructed $K$
 satisfies $\overline{\cup_{\mu\in K}S_\mu}=X,$
and thus by Lemma \ref{lem-transipoint-contains-generic-pt-meausre-withfullsupp},  $G_K$ is contained in D(T). Then we can follow the proof of Proposition  \ref{Main-Thm-1}, \ref{Main-Thm-2} and \ref{Main-Thm-3}   to complete the proof of   (C).

\bigskip

Case (D).
 Observe that    from the proofs of  Proposition    \ref{3Main-Thm-2015-R-minus-A11111} (2), \ref{3Main-Thm-1}, \ref{3Main-Thm-2} and  \ref{3Main-Thm-3},   each  constructed $K$ (or a single measure) satisfies $\overline{\cup_{\mu\in K}S_\mu}=X,$  and thus by Lemma \ref{lem-transipoint-contains-generic-pt-meausre-withfullsupp},  $G_K$ is contained in D(T). Then we can follow the proof of Proposition   \ref{3Main-Thm-2015-R-minus-A11111} (2), \ref{3Main-Thm-1}, \ref{3Main-Thm-2} and  \ref{3Main-Thm-3} to complete the proof of  (D).

\bigskip

Case (B).    One can follow  same idea as the proof of  Theorem \ref{More-general-systems-2} (1) to show that (D) implies (B).

\bigskip

Case (A).   Take $\phi\equiv 1$ in (B), then  $R_\phi(T)=X$ and thus (B) implies  (A). \qed

\bigskip


In particular, for full shifts of finite symbols, we also have some more observation.

\begin{Thm}\label{Thm-fullshift-Irr-Transtive2015Jan12}

 Let  $\,T$ be a     full shift on $k$ symbols ($k\geq 2$).  Then   $\{Per(T), A(T), A(T)\cup R(T), QR(T),$ $ W(T), V(T), $ $ QW(T), I(T)\}$ has full entropy gaps with respect to $X\setminus D(T)$.

\end{Thm}

{\bf Proof.} By Theorem \ref{Thm-fullshifts},  $A(T)\setminus Per(T)$ has full entropy.  Since $T$ is not minimal, by Lemma \ref{lem-transipoint-disjoint-AlmostPeriodc-non-minimalsystem}  $A(T)\cap D(T)=\emptyset.$ Then $A(T)\setminus Per(T)$ also has full entropy
with respect to $X\setminus D(T)$. Now we start to consider other cases.

Recall a classical result in $ \S 7.3$  of   \cite{Walter} that  for  full shift $T$ of   $k$ symbols,  $T$  has proper  subshifts (that is,  subshifts not equal to the full shift) with topological entropy equal to any given positive real number less than the topological entropy of the shift itself.  The constructed subshift is in fact the well-known $\beta-$shift ($\beta>1$). Recall that  any $\beta-$shift is not minimal (which contains a fixed point) so that it is not uniquely ergodic. For convenience to explain, denote the  subshift by $T_\beta$ and the subspace by $\Sigma_\beta \subsetneqq  X.$


Note that  for $\beta<k,$  $\xi=I(T)\setminus QW(T),QW(T)\setminus V(T), V(T)\setminus W(T), W(T)\setminus QR(T),
QR(T)\setminus (R(T)\cup A(T)), R(T)\setminus A(T),$
 $$\xi\cap \Sigma_\beta 
 \subseteq
    \xi\cap (X \setminus D (T)).$$
By Theorem \ref{Thm-beta-shifts}, we have  $h_{top}(T,\xi\cap \Sigma_\beta)= \log\beta $ and thus  $h_{top}(T,\xi\cap (X \setminus D (T)))\geq \log\beta.$ Let $\beta \uparrow k$, then every set  $\xi\cap (X\setminus D (T))$ carries full  topological entropy of  $\log k.$ \qed

\bigskip

All in all, with the help to various periodic-like recurrence and regularity of points, for a certain class of dynamical systems,  we obtain a refined classification of transitive or non-transitive  points and each one carries full topological entropy.




\section{ Geometric characterization of  gap-sets}\label{section-geometric-Bowen-specification}

Under the assumption of Bowen's specification property,  we will show that the gap-sets of Theorem \ref{More-general-systems}, \ref{More-general-systems-2} and \ref{More-general-systems-3} are all  dense in the whole space.
For convenience of making these more precise, we introduce a concept as follows. Let $T:X\rightarrow X$  be a continuous map of a compact
metric space $X$.

\begin{Def}\label{def-full-entropygaps} For a collection of subsets $Z_1,Z_2,\cdots,Z_k\subseteq X$ ($k\geq 2$), we say $\{Z_i\}$ has {\it dense gaps} with respect to $Y\subseteq X$ if
$$  \overline{(Z_{i+1}\setminus Z_i)\cap Y}  = \overline{Y}  \,\,\,\text{ for all } 1\leq i<k.$$

\end{Def}
Often, but not always, the sets $Z_i$ are nested ($Z_i\subseteq Z_{i+1}$).

\bigskip

Now we state the geometric characterization as follows.

\begin{Prop}\label{Main-Thm-5} Let  $\,T$ be a      continuous map of a compact metric space $X$ with Bowen's specification property.
  Then \\

(A) $\{A(T), A(T)\cup R(T),  QR(T), W(T), V(T), QW(T), I(T)\}$  has dense gaps with respect to $X$;\\

(B) for any  $\phi\in C^0 (X)$,   $\{A(T), A(T)\cup R(T),  QR(T), W(T), V(T), QW(T), I(T)\}$  has dense gaps with respect to $R_\phi(T)$;\\

(C) for any  $\phi\in C^0 (X)$, if $I_\phi (T)\neq\emptyset$,  then  $\{  QR(T), W(T), V(T), QW(T), I(T)\}$  has dense gaps with respect to $I_\phi(T)$;\\

(D) for any  $\phi\in C^0 (X)$, if $I_\phi (T)\neq\emptyset$,  then  for any $a\in Int(L_\phi),$  $\{A(T)\cup R(T),  QR(T), $ $W(T), V(T), QW(T),$ $ I(T)\}$  has dense gaps with respect to $R_{\phi,a}(T)$.
\end{Prop}

Remark that Theorem \ref{Main-Thm-5} can be applicative to all topological mixing subshifts of finite type, systems restricted on topological mixing locally maximal
hyperbolic sets.

\bigskip


We need   some classical properties of Bowen's   specification property, see   \cite{DGS} (also see \cite{Sig,Bow,Bowen2}).

 \begin{Lem}\label{lem-DGS-periodic-point-measure-dense} (Proposition  21.3 and 21.8 in \cite{DGS})\\
Let  $\,T$ be a      continuous map of a compact metric space $X$ with Bowen's specification property. Then \\
  (1) the set of periodic points is dense in the whole space $X.$\\
  (2) the set of periodic measures is dense in the set of $T$-invariant measures.


 \end{Lem}

 \begin{Lem}\label{lem-DGS} (Proposition  21.9 and  21.12 in \cite{DGS})\\
Let  $\,T$ be a      continuous map of a compact metric space $X$ with Bowen's specification property.  Then there is a dense $G_\delta$ subset $\mathcal{R}$ of $T$-invariant measures such that for any $\mu\in\mathcal{R},$ $\mu$ is ergodic and $S_\mu=X$ (also saying $\mu$ has full support).
 \end{Lem}

\begin{Lem}\label{lem-saturatedset-G-K-dense} (Proposition 21.14 in \cite{DGS})\\
Let  $\,T$ be a      continuous map of a compact metric space $X$ with Bowen's specification property. Then for any  compact connected nonempty set $K \subseteq M (X, T ),$
$$G_K:=\{x\in X|\,M_x (T)=K\}$$ (called saturated set of $K$) is nonempty and  dense in $X.$ In particular, $$G_{max}:=\{x\in X|\,M_x (T)= M (X, T )\}$$ is nonempty and  contains a dense $G_\delta$ subset of $X.$

\end{Lem}


\bigskip

Now let us start to prove Proposition \ref{Main-Thm-5}.

{\bf Proof of Proposition \ref{Main-Thm-5}.} Notice that under the assumption of Bowen's specification, by Lemma \ref{lem-DGS-periodic-point-measure-dense} the periodic points are dense in $X$ (i.e., $\overline{Per(f)}=X$) and
 the periodic measures are dense in the space of invariant measures (i.e., $\overline{ M_{p} (T,X) }=M (T,X)).$ Moreover, by Lemma \ref{lem-entropy-dense-Ps-2005}, $T$ has entropy-dense property, since specification is stronger than  $g$-almost product  property.

(C).   Observe that the  constructions of $K$    in Propositions \ref{Main-Thm-1}-\ref{Main-Thm-4} all did not use the saturated property. So for such $K,$ by Lemma \ref{lem-saturatedset-G-K-dense}, we complete the proof of Proposition \ref{Main-Thm-5} (C).

\bigskip

(D).    Observe that the  constructions of $K$   in Propositions   \ref{3Main-Thm-2015-R-minus-A11111}, \ref{3Main-Thm-1}, \ref{3Main-Thm-2}, \ref{3Main-Thm-3}, \ref{3Main-Thm-4}  all did not use the saturated property. So for such $K,$ by Lemma \ref{lem-saturatedset-G-K-dense}, we complete the proof of Proposition \ref{Main-Thm-5} (D).

\bigskip

(B).  Firstly we prove that $\{ A(T)\cup R(T),  QR(T), W(T), V(T), $ $QW(T),$ $ I(T)\}$  has dense gaps with respect to $R_\phi(T)$. We need two steps, since (D) just holds for functions with $I_\phi(T)\neq \emptyset$.

 {\bf Step 1.} If $I_\phi(T)\neq \emptyset,$ then $Int(L_\phi)$ is nonempty and we can take $a\in Int(L_\phi)$. Observe that $ R_{\phi,a}(T) \subseteq R_\phi(T)$ and thus (D) implies that $\{ A(T)\cup R(T),$ $  QR(T), W(T), V(T), $ $QW(T),$ $ I(T)\}$  has dense gaps with respect to $R_\phi(T)$.

{\bf Step 2.}
If $I_\phi(T)=\emptyset,$ then $R_\phi(T)=X.$
Similar as Step 2 in the proof of Theorem \ref{More-general-systems-2} (1),     there is some continuous function $\varphi:X\rightarrow \mathbb{R}$ such that   $I_\varphi(T)\neq \emptyset.$  Then by Step 1 we have  $\{A(T)\cup R(T),  QR(T), W(T), V(T), QW(T), I(T)\}$  has dense gaps with respect to $R_\varphi(T)$.
Notice that $R_\varphi(T)\subseteq X=R_\phi(T),$  then   $\{A(T)\cup R(T),  QR(T), W(T),$ $ V(T), QW(T), I(T)\}$ also has dense gaps with respect to $R_\phi(T)$.
\bigskip

Secondly let us consider the gap-set $R(T)\setminus A(T)$.  By Lemma \ref{lem-DGS}, we can take an ergodic measure $\mu$ with full support.
Remark that $G_\mu=G_\mu\cap X=G_\mu\cap S_\mu\subseteq R(T).$ Notice that $T$ is not minimal so that  by Lemma \ref{lem-non-minimalsupport-disjoint-almostperiod}  $S_\mu=X$ implies that $G_\mu\cap A(T)=\emptyset.$ So $G_\mu\subseteq R(T)\setminus A(T).$   By Lemma \ref{lem-saturatedset-G-K-dense}, $G_\mu$ is dense in $X$ and thus $R(T)\setminus A(T)$ is also dense in $X$.
Now we complete the proof of Proposition \ref{Main-Thm-5} (B).

\bigskip

(A).   Take $\phi\equiv 1$ in (B), then  $R_\phi(T)=X$ and thus (B) implies  (A).
 Now we complete the proof of Proposition \ref{Main-Thm-5}. \qed

\bigskip

\bigskip

In particular, we have a better characterization for the set of $V(T)\setminus W(T)$, that is, $$ I_\phi (T)\cap\{x\in  QW (T)\setminus W (T)|\exists\, \omega\in M_x (T) \,s.t. \,S_\omega=C_x\},$$ which is to answer the open problem in Section \ref{zhoufeng-OpenProblem} from different sight.

\begin{Prop}\label{Main-Thm-6} Let  $\,T$ be a      continuous map of a compact metric space $X$ with Bowen's specification property.
Let $$IC^0 (X):=\{\phi\in C^0 (X)|\,\,I_\phi (T)\neq\emptyset\}.$$  Then the set $$\bigcap_{\phi\in IC^0 (X)}I_\phi (T)\cap\{x\in  QW (T)\setminus W (T)|\exists\, \omega\in M_x (T) \,s.t. \,S_\omega=C_x\}$$ contains a dense $G_\delta$ subset  of $X$ (called residual in $X$).

\end{Prop}

{\bf Proof} By Lemma \ref{lem-saturatedset-G-K-dense}, $G_{max}=\{x\in X|\,M_x (T)= M (X, T )\}$ is residual in $X$. So we only need to prove for any $\phi\in IC^0 (X),$ $$G_{max}\subseteq  I_\phi (T)\cap \{x\in QW (T)\setminus W (T)|\exists\, \omega\in M_x (T) \,s.t. \,S_\omega=C_x\}.$$

Since $\phi\in IC^0 (X)$, then $I_\phi(T)\neq \emptyset$ and thus by  (\ref{Ir-phi-nonempty-imply-differentmeasure})  there are two invariant measures $\omega,\mu\in M (X, T )$ such that they have different  integrals  for  $\phi.$ Fix  $x\in G_{max}.$ Remark that  $M (T, X )= M_x (T)$  and so $\mu,\omega\in M_x (T)$.  By (\ref{IrRegular-phi-equivalent})  $x\in I_\phi (T).$
By density of periodic points (Lemma \ref{lem-DGS-periodic-point-measure-dense}), $$C_x=\overline{\bigcup_{m\in M_x (T)}S_m}=\overline{\bigcup_{m\in M (T,X)}S_m}\supseteq\overline{\bigcup_{m\in M (T,X)\,\,is\,\,periodic \,\,measure}S_m}=X$$ (Remark that $C_x=X$ can be also obtained from the existence of invariant measure with full support by Lemma \ref{lem-DGS}).
By (\ref{C_x=X})    $  X=C_x $ implies   $x\in QW (T)$.  By (\ref{weakperiodic-in-reccurent}) $x\in Rec (T).$ By Lemma \ref{lem-DGS-periodic-point-measure-dense}, one can take a periodic measure $\nu\in M (T,X)=M_x (T)$ whose support $S_\mu\subsetneqq X=C_x$, then by (\ref{WT}) $x\in X\setminus W (T).$
Take an invariant measure $\mu$ with full support $X$ by Lemma \ref{lem-DGS}, then $\mu\in M (T,X)=M_x (T)$ and $S_\mu=X=C_x.$ We complete the proof.
 \qed

\bigskip

Remark that $G_{max}$ has zero topological entropy from (\ref{eq-generalentropyestimate-cptconnected-K}), since $K:=M (T,X)$ contains periodic measures which have zero entropy.
By (\ref{Bowen-disjointSet-entropy}) and Theorem \ref{Thm-specification1} we know that

\begin{Thm}\label{Thm-specification-minus-G-max} Let  $\,T$ be a      continuous map of a compact metric space $X$ with Bowen's specification property and uniform separation. The complementary set of $G_{max}$,  $$ \{x\in  I_\phi (T)\cap QW (T)\setminus W (T)|\exists\, \omega\in M_x (T) \,s.t. \,S_\omega=C_x\}\setminus G_{max},$$ has full topological entropy. This complementary set is just dense but not residual in $X$.

\end{Thm}

{\bf Proof.}  Full entropy can be deduced  from (\ref{Bowen-disjointSet-entropy}) and Theorem \ref{Thm-specification1}. Density can be seen from Lemma \ref{lem-saturatedset-G-K-dense}, since the choice  of $K$ in the proof of Proposition \ref{Main-Thm-2} is a proper subset of $M(T,X)$ and so $G_K\cap G_{max}=\emptyset$. \qed

\bigskip

Recall from \cite{Tho2012}   that for any system with almost specification (not necessarily satisfying uniform separation),
  every $\phi-$irregular set either is empty or carries full topological entropy. Together with the observation of Proposition \ref{Main-Thm-5}, a natural question aries for systems with Bowen's specification property (not necessarily satisfying uniform separation):

\begin{Que}\label{Que-Irregular-full-entropy}
Let $f$ be a     continuous map of a compact metric space $X$ with  Bowen's specification (or almost  specification).  Whether all the results in Theorem \ref{More-general-systems},  \ref{More-general-systems-2} and Theorem \ref{More-general-systems-3} hold? If not, which results hold, which results do not hold and how to construct such a counterexample?
\end{Que}

Let us give some simple observation for Question \ref{Que-Irregular-full-entropy} in the case of Bowen's specification.
   Since Bowen's specification is stronger than $g$-almost product property, then by Theorem \ref{lem-PS-etds2007-no-uniform-separation} $f$ is single-saturated. By Lemma \ref{lem-DGS-full-support-2015-22222}, there is an invariant measure with full support. So we can use
Proposition \ref{3Main-Thm-2015-R-minus-A11111} (2) to obtain $\{A(T)\cup R(T), QR(T)\}$ has full entropy gaps with respect to $R_{\phi,a}(T)$ (resp., $R_\phi(T)$ and $X$). For other results of Theorem \ref{More-general-systems},  \ref{More-general-systems-2} and Theorem \ref{More-general-systems-3}, from the proofs of such results the considered  $K\subseteq M(T,X)$ is not singleton. So single-saturated property is not enough.  Moreover, let us recall another possible  idea by  Thompson  \cite{Tho2012}, one needs to take two needed ergodic measures and then use these two measures to construct a set $F\subseteq  I_\phi (T)$ such that the topological entropy of $F$ is larger than $h_{top} (f)-\epsilon.$ In this process, Entropy Distribution Principle plays an important role. One can see   \cite{To2010,Tho2012} for more details. The constructed measures are required ergodic. However,   the constructed measures in present paper are not all ergodic so that we are not sure the idea of \cite{To2010,Tho2012} works. So Question \ref{Que-Irregular-full-entropy} is still open except the case for $\{A(T)\cup R(T), QR(T)\}$.

Remark that if the answer  of Question \ref{Que-Irregular-full-entropy}  is positive, then it can be as  a generalization of Theorem \ref{Thm-specification1}.  In particular, if Question \ref{Que-Irregular-full-entropy} is true,  it would  be applicative to all topological mixing interval maps, since it is known  from \cite{Bl,Buzzi} that any topologically mixing interval map satisfies Bowen's  specification. For example, Jakobson  \cite{Jakobson} showed that there exists a set of parameter values $\Lambda\subseteq [0,4]$ of positive Lebesgue measure such that if $\lambda\in \Lambda,$ then the logistic map $f_\lambda(x)=\lambda x(1-x)$ is topologically mixing. From Proposition 21.4 of \cite{DGS} we also know that the factor of a system with Bowen's specification   has Bowen's specification.



\bigskip

At the end of this section, we show that Proposition \ref{Main-Thm-6} holds for a larger class of systems, such as  $C^1$ generic (volume-preserving) diffeomorphisms, which means the open problem of   \cite{ZF} by Zhou and Feng is solved  for generic systems. Let $M$ be a compact Riemannian manifold and $m$ be a volume measure on $M.$ Let $\Diff^1 (M)$ and $\Diff^1_m (M)$ denote the space of all $C^1$ diffeomorphisms on $M$ and all volume-preserving $C^1$ diffeomorphisms on $M$ respectively. If $f:M\rightarrow M$  and $\Lambda$ is a compact subset of $M,$ let $$IC_f^0 (\Lambda):=\{\phi\in C^0 (\Lambda)|\,\,I_\phi (f)\neq\emptyset\}.$$

\begin{Thm}\label{Main-Thm-6-zhou-open-problem}
 (1) Let $\Lambda$ be an isolated non-trivial transitive set of a $C^1$ generic diffeomorphism $f\in \Diff^1 (M)$. Then the set $$\bigcap_{\phi\in IC_f^0 (\Lambda)} I_\phi (f)\cap\{x\in \Lambda \cap  QW (f)\setminus W (f)|\exists\, \omega\in M_x (f) \,s.t. \,S_\omega=C_x\}$$ contains a dense $G_\delta$ subset  of $\Lambda$ (called residual in $\Lambda$).

 (2) Let $f\in\Diff_m^1 (M)$ be a $C^1$ generic volume-preserving diffeomorphism. Then the set $$\bigcap_{\phi\in IC_f^0 (M)}I_\phi (f)\cap \{x\in   QW (f)\setminus W (f)|\exists\, \omega\in M_x (f) \,s.t. \,S_\omega=C_x\}$$ contains a dense $G_\delta$ subset  of $M$.
\end{Thm}

{\bf Proof.}  For the first case, by the main result of   \cite{suntian-2012} $$G_{max}=\{x\in \Lambda|\,M_x (f)= M (f,\Lambda  )\}$$ is residual in $X$. Recall  Theorem 3.5 of \cite{ABC} that generic invariant measures have full support. Thus, forward the proof of Proposition \ref{Main-Thm-6},  one can replace Lemma \ref{lem-DGS}  by Theorem 3.5 of   \cite{ABC} to prove.

For the second case, it is known that generic $f\in\Diff_m^1 (M)$ is transitive so that we can take $\Lambda=M$. Notice that Theorem 3.5 of   \cite{ABC}
and  the main result of   \cite{suntian-2012} also can be stated as the volume-preserving case. Then the proof is similar as above.
Here we omit the details.
\qed

\bigskip

 Inspired by Theorem \ref{Main-Thm-6-zhou-open-problem}, it is possible to ask

\begin{Que}\label{Que-Irregular-full-entropy-generic-diffeomorphism}
Let $\Lambda$ be an isolated non-trivial transitive set of a $C^1$ generic diffeomorphism $f\in \Diff^1 (M)$or let $f\in\Diff_m^1 (M)$ be a $C^1$ generic volume-preserving diffeomorphism. Then whether $\{W(f), V(f)\}$ (resp.,  $\{W(f), QW(f)\}$) has full entropy gaps with respect to $\Lambda$ or $M$?
\end{Que}

Recall that uniformly hyperbolic systems are not dense in the space  of all diffeomorphisms.
  So we can not use Theorem \ref{Thm-hyperbolic} to answer this question. 
  
Moreover, it is unknown whether we can use our main theorems (Theorem \ref{More-general-systems}-\ref{More-general-systems-3}) to answer this question.
From \cite{ABC} periodic points are all hyperbolic and dense in $\Lambda$ and periodic measures are dense in the space of all invariant measures supported on $\Lambda.$  However,   we do not know whether $g$-almost product property and uniform separation hold  for generic diffeomorphisms. 
So 
Question \ref{Que-Irregular-full-entropy-generic-diffeomorphism}  seems to be  nontrivial.



\bigskip

{\bf Acknowledgements.} The author  thanks  the anonymous referee very much  for his or her  constructive suggestions and careful reading. For  example, the referee introduces  the concept of {\it full entropy gaps} which makes the statement  of main results more easier to write and read.

\end{document}